\documentclass[12pt]{amsart}
\usepackage {amscd, amsthm}

\theoremstyle{plain} 
\newtheorem{thm}{Theorem}[subsection]

\newtheorem{lem}[thm]{Lemma}

 \theoremstyle{definition}
\newtheorem{defn}[thm]{Definition}

\theoremstyle{remark}
\newtheorem{exmpl}[thm]{Example}
\newtheorem{rem}[thm]{Remark}

\numberwithin{equation}{section}

\newcommand{\Ann}{\operatorname{Ann}}

\newcommand{\ber}{\operatorname{ber}}
\newcommand{\berdual}{\operatorname{ber}^*}
\newcommand{\Ber}{\operatorname{Ber}}
\newcommand{\Bersgr}{\operatorname{Ber(\Sgr)}}
\newcommand{\Bersgrdual}{\operatorname{Ber^*(\Sgr)}}
\newcommand{\Berx}{\mathcal B\text{er}_{X}}

\newcommand{\Berxsplit}
{\mathcal B\text{er}_{X}^{\text{split}}}
\newcommand{\Berxhat}{\mathcal B\text{{\^e}r}_{{X}}}
\newcommand{\Berxhatsplit}{\mathcal B\text{\^{e}r}_{{X}}\spl}

\newcommand{\Cox}{\mathcal{CO}_X}
\newcommand{ \datum}{\today}
\newcommand{\Dc}{D_{\mathcal C}}
\newcommand{\Dchat}{\hat{D}_{\mathcal C}}
\newcommand{\Div}{\operatorname{Div}}
\newcommand{\Divx}{\Div(X)}
\newcommand{\done}{D^{1|1}}

\newcommand{\halfz}{\frac12 \mathbb Z}
\newcommand{\hatxd}{\hat X^{(d)}}
\newcommand{\Jacx}{\operatorname{Jac}(X)}
\newcommand{\Jheis}{\mathcal J\text{Heis}}
\newcommand{ \inv}{^{-1}}
\newcommand{\Linfinf}{\mathbb \Lambda^{\infty\mid\infty}}%

\newcommand{\Ox}{\mathcal O_X}
\newcommand{\Oxev}{\mathcal O_{X,\text{ev} } }
\newcommand{\Oxevstar}{\mathcal O^\times_{X,\text{ev}}}
\newcommand{\Oxhat}{\hat{\mathcal O}_X}

\newcommand{\Oxhatsplit}{\hat{\mathcal O}_X^{\text{split}}}
\newcommand{\Oxred}{{\mathcal O}_{X}^{\text{red}}}
\newcommand{\Oxsplit}{{\mathcal O}_{X}\spl}

\newcommand{\per}{\operatorname{per}}
\newcommand{\Pic}{\operatorname{Pic}}
\newcommand{\Picxz}{\Pic^0(X)}

\newcommand{\Ratx}{\mathcal Rat(X)}
\newcommand{\Ratxevstar}{\mathcal Rat_{\text{ev}}^\times(X)}
\newcommand{ \red}{^{\,\text{red}}}
\newcommand{\Sfr}{\mathcal S\text{fr}}
\newcommand{\Sgr}{\mathcal S\text{gr}}
\newcommand{\Sheis}{\mathcal S\text{Heis}}
\newcommand{\Spec}{\operatorname{Spec}}
\newcommand{\Str}{\operatorname{Str}}
\newcommand{ \spl}{^{\,\text{split}}}

\begin{document}
\title[Super Curves]
{Super Curves, their Jacobians, and super KP equations}

    \author{M.J. Bergvelt}
    \address[Bergvelt] {Department of Mathematics,
University of Illinois, Urbana, IL 61801, USA }
    \email{bergv@math.uiuc.edu}
    \author{J.M. Rabin}
    \address[Rabin] {Department of Mathematics, University
of California at San Diego,
    	La Jolla, CA, 92093, USA}
    \email{jrabin@euclid.ucsd.edu }
    \date{\datum}
\begin{abstract}
We study the geometry and cohomology of algebraic

super curves, using a new
contour integral for holomorphic differentials. For a class of super
curves
(``generic SKP curves'') we define
a period matrix. We show that the odd part of the period matrix
controls the cohomology of the dual curve.
 The Jacobian of a generic SKP curve is a
smooth supermanifold; it is principally polarized, hence projective,
if the
even
part of the period matrix is symmetric.
 In general symmetry is not guaranteed by the
Riemann bilinear equations for our contour integration, so it remains
open
whether  Jacobians are always
projective  or carry theta functions.

\hskip.3truein
These results on generic SKP curves are applied to the study of
algebro-geometric  solutions of the super KP hierarchy.
The tau function is shown to be, essentially, a meromorphic section
of
a
line bundle with trivial Chern class on the Jacobian, rationally
expressible in
terms of super theta functions when these exist. Also we relate the
tau function and the Baker function for this hierarchy,
 using a generalization of Cramer's rule to the supercase.

\end{abstract}

\maketitle

\tableofcontents

\section{Introduction.}
In this paper we study algebraic super curves with a view towards
applications
to super Kadomtsev-Petviashvili hierarchies (SKP). We deal from the
start with
super curves $X$ over a nonreduced ground ring $\Lambda$, i.e., our
curves
carry global nilpotent constants. This has as an advantage, compared
to super
curves over
the complex numbers $\mathbb C$, that our curves can be nonsplit, but
this
comes at the price of some technical complications. The main problem
is that
the cohomology groups of coherent sheaves on our curves should be
thought of as
finitely generated modules over the ground ring $\Lambda$, instead of
vector
spaces over $\mathbb C$. In general these modules are of course not
free.
Still we have in this situation Serre duality, as explained in
Appendix
\ref{app:dualSerredual}, the dualizing sheaf being the (relative)
berezinian
sheaf $\Berx$. In applications to SKP there occurs a natural class of
super
curves that we call generic SKP curves. For these curves the most
important
sheaves, the structure sheaf and the dualizing sheaf, have free
cohomology. In
the later part of the paper we concentrate on these curves.

Super curves exhibit a remarkable duality uncovered in
\cite{DoRoSc:SupModSpaces}.
The projectivized cotangent bundle of any $N=1$ super curve has the
structure of an $N=2$ super Riemann surface (SRS), and super curves
come in
dual pairs $X,\hat{X}$ whose associated $N=2$ SRSs coincide.
Further, the ($\Lambda$-valued) points of a super curve can
be identified with the effective divisors of degree 1 on its dual.
Ordinary
$N=1$ SRSs, widely studied in the context of super string theory, are
self dual
under this duality. By the resulting identification of points with
divisors they enjoy many of the properties that distinguish Riemann
surfaces from higher-dimensional varieties.
By exploiting the duality we extend this good behaviour to all
super curves.

In particular we define for all super curves a contour integration
for sections
of $\Berx$, the holomorphic differentials in this situation. The
endpoints of a
super contour turn out to be not $\Lambda$-points of  our super
curve, but
rather irreducible divisors, i.e., $\Lambda$-points
on the dual curve! For SRSs these notions are the same and our
integration is a
generalization of the procedure already known for SRSs. We use this
to prove
Riemann bilinear relations, connecting in this situation periods of
holomorphic
differentials on our curve $X$ with those on its dual curve.

In case the cohomology of the structure sheaf is free, e.g. if $X$ is
a generic
SKP curve, we can define a period matrix and use this to define the
Jacobian
of $X$ as the quotient of a super vector space by a lattice generated
by the
period matrix. In this case the Jacobian is a smooth supermanifold. A
key
question is whether the Jacobian of a generic SKP curve admits ample
line
bundles (and hence embeddings in projective super space),
whose sections could serve as the super analog of theta functions.
We show that the symmetry of the even part of the period matrix
(together with the automatic positivity of the imaginary part of the
reduced matrix) is sufficient for this, and construct the super
theta functions in this case.
We derive some geometric necessary and sufficient conditions for
this symmetry to hold, but it is not an automatic consequence of the
Riemann bilinear period relations in this super context.
Neither do we know an explicit example in which the symmetry fails.
The usual proof that symmetry of the period matrix is necessary for
existence of a (principal) polarization also fails because crucial
aspects of Hodge theory, particularly the Hodge decomposition of
cohomology, do not hold for supertori.

The motivation for writing this paper was our wish to generalize
the theory of the
algebro-geometric solutions to the KP hierarchy of nonlinear PDEs,
as described in \cite{SeWi:LpGrpKdV} and references
therein, to the closest supersymmetric analog, the ``Jacobian"
super KP  hierarchy of Mulase and Rabin \cite{Mu:Jac,Ra:GeomSKP}.

In the super KP case the geometric data leading to a solution include
a super
curve $X$ and a line bundle $\mathcal L$  with vanishing cohomology
groups over
$X$. For such a line bundle to exist the super curve $X$ must have a
structure
sheaf $\Ox$ such that the associated split sheaf $\Oxsplit$, obtained
by
putting the global nilpotent constants in $\Lambda$ equal to zero, is
a direct
sum $\Oxsplit=\Oxred\oplus \mathcal N$, where $\Oxred$ is the
structure sheaf
of the underlying classical curve $X^{\text{red}}$ and $\mathcal N$
is an
invertible $\Oxred$-sheaf of degree zero. We call such an $X$ an SKP
curve, and
if moreover $\mathcal N$ is not isomorphic to $\Oxred$ we call $X$ a
generic
SKP curve.

The Jacobian SKP hierarchy describes linear
flows $\mathcal L(t_i)$ on the Jacobian of $X$ (with even and odd
flow
parameters).
The other known SKP hierarchies, of Manin--Radul
\cite{ManRad:SusyextKP} and
Kac--van de
Leur \cite{KavdL:SuperBoson} , describe flows on the universal
Jacobian over
the moduli space
of super curves, in which $X$ as well as $\mathcal L$ vary with the
$t_i$
\cite{Ra:GeomSKP}.
These are outside the scope of this paper, although we hope to
return to them elsewhere.
As in the non-super case, the basic objects in the theory are the
(even and odd) Baker functions, which are sections of $\mathcal
L(t_i)$
holomorphic except for a single simple pole, and a tau
function which is a section of the super determinant (Berezinian)
bundle over a super Grassmannian $\Sgr$.
In contrast to the non-super case, we show that the Berezinian
bundle has trivial Chern class, reflecting the fact that the
Berezinian is a ratio of ordinary determinants.
The super tau function descends, essentially, to $\Jacx$ as a section
of a
bundle with trivial Chern class also, and can be rationally expressed
in
terms of super theta functions when these exist (its reduced part is
a ratio of ordinary tau functions).
We also obtain a formula for the even and odd Baker functions in
terms of the tau function, confirming that one must know the tau
function for the more general Kac--van de Leur flows to compute the
Baker functions for even the Jacobian flows in this way, cf.
\cite{DoSc:SuperGrass,Takama:GrassmannSKP}.
For this we need a slight extension of Cramer's rule for solving
linear equations in even and odd variables, which is developed in an
Appendix via the theory of quasideterminants.
In another Appendix we use the Baker functions found in
\cite{Ra:SupElliptic} for Jacobian flow in the case of super elliptic
curves to
compute the corresponding tau function.

Among the problems remaining open we mention the following.
First, to obtain a sharp criterion for when a super Jacobian admits
ample line bundles --- perhaps always?
Second, the fact that generic SKP curves have free cohomology is a
helpful simplification which allows us to represent their period
maps by matrices and results in their Jacobians being smooth
supermanifolds.
However, our results should generalize to arbitrary super curves with
more singular Jacobians.
Finally, one should study the geometry of the universal Jacobian and
extend our analysis to the SKP system of Kac--van de Leur.

\section{Super curves and their Jacobians.}

\subsection{Super curves.}
Fix a Grassmann algebra $\Lambda$ over $\mathbb C$; for instance we
could take
$\Lambda=\mathbb C\,[\beta_1,\beta_2,\dots,\beta_n]$, the polynomial
algebra
generated by $n$ odd indeterminates. Let $(\bullet,\Lambda)$
be the super scheme $\Spec \Lambda $, with underlying topological
space a
single point $\bullet$.

A smooth compact connected complex super curve over $\Lambda$ of
dimension
$(1|N)$ is a pair $(X,\Ox)$, where $X$ is a topological space and
$\Ox$ is a
sheaf of super commutative $\Lambda$-algebras
over $X$, equipped with a structure morphism $(X,\Ox)\to
(\bullet,\Lambda)$,
such that
\begin{enumerate}
\item $(X,\Oxred)$ is a smooth  compact connected complex  curve,
algebraic or
holomorphic, depending on the category one is working in. Here
$\Oxred$ is the
reduced sheaf of $\mathbb C$-algebras on $X$ obtained by quotienting
out the
nilpotents in the structure sheaf $\Ox$,
\item For suitable open sets $U_\alpha\subset X$ and suitable
linearly
independent odd elements $\theta_\alpha^i$ of $\Ox(U_\alpha)$ we have
$$
\Ox(U_\alpha)=\Oxred\otimes
\Lambda[\theta_\alpha^1,\theta_\alpha^2,\dots,\theta_\alpha^N].
$$
\end{enumerate}
The $U_\alpha$\rq s above are called coordinate neighborhoods of
$(X,\Ox)$ and
$(z_\alpha,\theta_\alpha^1,\theta_\alpha^2,\dots,\theta_\alpha^N)$
are called
local coordinates  for $(X,\Ox)$, if $z_\alpha$ (mod nilpotents) is a
local
coordinate for $(X,\Oxred)$. On overlaps of coordinate neighborhoods
$U_\alpha\cap U_\beta$ we have
\begin{equation}\label{eq:coordchange}
\begin{split}
z_\beta 	  &= F_{\beta\alpha}(z_\alpha,\theta_\alpha^j),\\
\theta_\beta^i  &= \Psi_{\beta\alpha}^i(z_\alpha,\theta_\alpha^j).
\end{split}
\end{equation}
Here  the $F_{\beta\alpha}$ are even functions and
$\Psi_{\beta\alpha}^i$  odd
ones, holomorphic or algebraic depending on the category we are
using.

\begin{exmpl}\label{exmpl:split} A special case is formed by the {\it
split}
super curves. For $N=1$ they are given by transition functions
\begin{equation}\label{eq:splitcoordchange}
\begin{split}
z_\beta 	  &= f_{\beta\alpha}(z_\alpha),\\
\theta_\beta  &= \theta_\alpha B_{\beta\alpha}(z_\alpha),
\end{split}
\end{equation}
with $f_{\beta\alpha}(z_\alpha),B_{\beta\alpha}(z_\alpha)$ even
holomorphic (or
algebraic) functions that are independent of the nilpotent constants
in
$\Lambda$.
So in this case the $f_{\beta\alpha}$ are the transition functions
for $\Oxred$
and $\Ox=\Oxred\otimes\Lambda \mid \mathcal N\otimes \Lambda$, where
$\mathcal
N$ is the $\Oxred$-module with transition functions
$B_{\beta\alpha}(z_\alpha)$. Here and henceforth we denote by a
vertical $\mid$
a direct sum of  free $\Lambda$-modules, with on the left an evenly
generated
summand and on the right an odd one.

To any super curve $(X,\Ox)$ there is canonically associated a split
curve
$(X,\Oxsplit)$ over $\mathbb C$: just take
$\Oxsplit=\Ox\otimes_{\Lambda}
\Lambda/\mathfrak{m}=\Ox/\mathfrak{m}\Ox$, with
$\mathfrak m=\langle \beta_1,\dots,\beta_n\rangle$  the maximal ideal
of
nilpotents in $\Lambda$. There is a functor from the category of
$\Ox$-modules
to the category of $\Oxsplit$-modules that associates to a sheaf
$\mathcal F$
the {\it associated split sheaf} $\mathcal F^{\text{split}}=\mathcal
F/\mathfrak {m}\mathcal F$.\qed\end{exmpl}
\smallskip
A $\Lambda$-point of a super curve $(X,\Ox)$ is a morphism
$\phi:(\bullet,\Lambda)\to  (X,\Ox)$ such that the composition with
the structural morphism \linebreak $(X,\Ox)\to(\bullet,\Lambda)$ is
the
identity  (of $(\bullet,\Lambda)$). Locally, in an open set
$U_\alpha$
containing $\phi(\bullet)$,  a $\Lambda$-point is given by specifying
the
images under the even $\Lambda$-homo\-mor\-phism
$\phi^\sharp:\Ox(U_\alpha)\to
\Lambda$
of the local coordinates:
$p_\alpha=\phi^\sharp(z_\alpha),\pi^i_\alpha
=\phi^\sharp( \theta^i_\alpha)$ . The local parameters
$(p_\alpha,\pi^i_\alpha)$ of a $\Lambda$-point transform precisely as
the
coordinates do, see (\ref{eq:coordchange}). By quotienting out
nilpotents in a
$\Lambda$-point $(p_\alpha,\pi^i_\alpha)$ we obtain a complex number
$p_\alpha^{\text{red}}$, the coordinate of the reduced point of
$(X,\Oxred)$
corresponding to the $\Lambda$-point $(p_\alpha,\pi^i_\alpha)$.

\subsection{Duality and $N=2$ curves.}\label{ss:DualN=2curves}
Our main interest is the theory of $N=1$ super curves but  as a
valuable tool
for the study of these curves we make use of $N=2$ curves as well in
this
paper. Indeed, as is well known,
\cite{DoRoSc:SupModSpaces,Schwarz:SuperanalogsSYMPLCONT}, one can
associate in
a canonical way to an $N=1$ curve an (untwisted) super conformal
$N=2$ curve,
as we  will now recall. The introduction of the super conformal $N=2$
curve
clarifies the whole theory of $N=1$ super curves.

Let from now on $(X,\Ox)$ be an $N=1$ super curve. Any invertible
sheaf
$\mathcal E$ for $(X,\Ox)$ and any extension
of $\mathcal E$ by the structure sheaf:
\begin{equation*}
0\rightarrow \Ox\rightarrow \hat{\mathcal E}\rightarrow \mathcal
E\rightarrow
0,\label{eq:extension}
\end{equation*}
defines in the obvious way an $N=2$ super curve $(X,\hat{\mathcal
E})$. It has
local coordinates $(z_\alpha,\theta_\alpha,\rho_\alpha)$, where
$(z_\alpha,\theta_\alpha)$ are local coordinates for $(X,\Ox)$. On
overlaps we
will have
\begin{equation}\label{eq:coordchangen=2}
\begin{split}
z_\beta 	&= F_{\beta\alpha}(z_\alpha,\theta_\alpha),\\
\theta_\beta  	&= \Psi_{\beta\alpha}(z_\alpha,\theta_\alpha),\\
\rho_\beta 	&= H_{\beta\alpha}(z_\alpha,\theta_\alpha)
		\rho_\alpha + \phi_{\beta\alpha}
			(z_\alpha,\theta_\alpha).
\end{split}
\end{equation}
(So $H_{\beta\alpha}(z_\alpha,\theta_\alpha)$ is the transition
function for
the generators of the invertible sheaf $\mathcal E$.) We want to
choose
the extension (\ref{eq:extension}) such that $(X,\hat{\mathcal E})$
is {\it
super conformal}, in the sense that the local differential form
$\omega_\alpha=dz_\alpha- d\theta_\alpha \rho_\alpha$   is globally
defined up
to a scale factor. Now
$$
\omega_\beta	=dz_\beta- d\theta_\beta \rho_\beta
		=dz_\alpha
	(\frac{\partial F}{\partial z_\alpha} - \frac{\partial
\Psi}{\partial
z_\alpha}\rho_\beta) -
	d\theta_\alpha
	(-\frac{\partial F}{\partial \theta_\alpha} + \frac{\partial
\Psi}{\partial
\theta_\alpha}\rho_\beta).
$$
(Here we suppress the subscripts on $F$ and $\Psi$, as we will  do
below.) We
see that for $\hat{\mathcal E}$ to be super conformal we need
$$
\rho_\alpha=\frac{(-\frac{\partial F}{\partial \theta_\alpha} +
\frac{\partial
\Psi}{\partial \theta_\alpha}\rho_\beta)}{
(\frac{\partial F}{\partial z_\alpha} - \frac{\partial \Psi}{\partial
z_\alpha}\rho_\beta)},
$$
or
\begin{equation}\label{eq:transformrho}
\rho_\beta=\frac{
(\frac{\partial F}{\partial \theta_\alpha} + \frac{\partial
F}{\partial
z_\alpha}\rho_\alpha)}
{(\frac{\partial \Psi}{\partial \theta_\alpha} - \frac{\partial
\Psi}{\partial
z_\alpha}\rho_\alpha)}.
\end{equation}
Conversely one checks that if (\ref{eq:transformrho}) holds for all
overlaps
the cocycle condition is satisfied and that we obtain in this manner
an $N=2$
super curve. To show that this super curve is an extension as in
(\ref{eq:extension}), it is useful to note that
(\ref{eq:transformrho}) can
also be written as
\begin{equation}\label{eq:transformrho2}
\rho_\beta=\ber\begin{pmatrix}\partial_z F&\partial_z \Psi\\
			\partial_\theta F&\partial_\theta\Psi
			\end{pmatrix}
	\rho_\alpha+\frac{\partial_\theta F}{\partial_\theta\Psi}.
\end{equation}
The homomorphism $\ber$ is defined in Appendix
\ref{app:Lineqsupercat}, see
\eqref{eq:defberber*}.
Recall that the local generators $f_\alpha$ of the dualizing sheaf
(see
Appendix \ref{app:dualSerredual}) $\mathcal{B}\text{er}_X$ of
$(X,\Ox)$
transform as
\begin{equation}\label{eq:transfBer}
f_\beta=\ber\begin{pmatrix}\partial_z F&\partial_z \Psi\\
			\partial_\theta F&\partial_\theta\Psi
			\end{pmatrix}f_\alpha.
\end{equation}

If we denote by $\Cox$ the structure sheaf of the super conformal
$N=2$ super
curve just constructed, we see that we have an exact sequence
\begin{equation}
0\rightarrow\Ox\rightarrow\Cox\rightarrow
\mathcal{B}\text{er}_X\rightarrow0.
\label{eq:extberbystruct}\end{equation}
$\Cox$ is the only extension of $\Berx$ by the structure sheaf that
is super
conformal. This sequence is {\it trivial} if it isomorphic
(\cite{HiSt:HomAlg})
to a split sequence.

\begin{defn}\label{def:projected}
A super curve is called {\it{projected} }if there is a cover of $X$
such that
the transition functions $F_{\beta\alpha}$ in (\ref{eq:coordchange})
are
independent of the odd coordinates $\theta_\alpha^j$.
\end{defn}

For projected curves we have a
projection morphism $(X,\Ox)\to (X,\Oxred\otimes \Lambda)$
corresponding to the
sheaf inclusion $\Oxred\otimes \Lambda\to \Ox$.
This inclusion can be defined only for projected curves.

A projected super curve  has a $\Cox$ that is a trivial extension but
the
converse is not true, as we will see when we discuss super Riemann
surfaces in
subsection
\ref{ss:SRS}.
The relation between projectedness of $(X,\Ox)$ and the triviality of
the
extension defining $(X,\Cox)$ is discussed in detail in subsection
\ref{ss:Symmperiodmatrices}.

\begin{exmpl}\label{exmpl:splitber}
If $(X,\Ox)$ is split, (\ref{eq:transfBer}) becomes
\begin{equation}\label{eq:transfBersplit}
f_\beta=\frac{\partial_z f_{\beta\alpha}}{B_{\beta\alpha}} f_\alpha.
\end{equation}
This means that in this case $\Berx=\mathcal K\mathcal N\inv\otimes
\Ox=\mathcal K\mathcal N\inv\otimes\Lambda\mid \mathcal K\otimes
\Lambda$,
where $\mathcal K$ is the canonical sheaf for $\Oxred$.

Split curves are projected and the sequence (\ref{eq:extberbystruct})
becomes
trivial. As an $\Oxred$-module we have $\Cox=(\Oxred \oplus \mathcal
K)\otimes\Lambda\mid (\mathcal N\oplus \mathcal K\mathcal
N\inv)\otimes
\Lambda$.
\qed\end{exmpl}
\smallskip

The map $\Cox\to \Berx$ is locally described by the differential
operator
$\Dc^\alpha=\partial_{\rho_\alpha}$.
Indeed, the operator $\Dc^\alpha$ transforms homogeneously,
$\Dc^\beta=\ber\begin{pmatrix}\partial_z F&\partial_z \Psi\\
			\partial_\theta F&\partial_\theta\Psi
			\end{pmatrix}\inv \Dc^\alpha$, so this
defines a
global $(0\mid 1)$ dimensional distribution $\Dc$ and the quotient of
$(X,
\mathcal {CO}_X)$ by this distribution is precisely $(X,\Ox)$.

Now the distribution $\Dc$ annihilates the 1-form $\omega$ used to
find $\Cox$.
This form  locally looks like $\omega_\alpha=dz_\alpha
-d\theta_\alpha
\rho_\alpha$ and its kernel is generated by  $\Dc^\alpha$ and  a
second
operator
$\Dchat^\alpha=\partial_{\theta_\alpha}+
\rho_\alpha\partial_{z_\alpha}$. (The
operators that we call $\Dc$ and $\Dchat$ are in the literature
also denoted by $D^+$ and $D^-$, cf. \cite{DoRoSc:SupModSpaces})
To study the result of ``quotienting by the distribution $\Dchat$''
we
introduce in each coordinate neighborhood $U_\alpha$ new coordinates:
\begin{equation*}
\begin{split}
\hat z_\alpha 		&= z_\alpha-\theta_\alpha\rho_\alpha,\\
\hat \theta_\alpha 	&=\theta_\alpha ,\\
\hat \rho_\alpha 	&= \rho_\alpha.
\end{split}
\end{equation*}
In the sequel we will drop the hats $\hat {}$ on $\theta$ and $\rho$,
hopefully
not causing too much confusion.

In these new coordinates we have
$$
\Dchat^\alpha= \partial_{\theta_\alpha},\quad
\Dc^\alpha=\partial_{\rho_\alpha}+ \theta_\alpha\partial_{\hat
z_\alpha}.
$$
So the kernel of  $\Dchat$ consists locally of functions of
$\hat z_\alpha,\rho_\alpha$. To see that this makes global sense we
observe
that
\begin{equation}\label{eq:coordtransfdualcurve}
\begin{split}
\hat z_\beta 	&= F(\hat z_\alpha,\rho_\alpha) + \frac{DF
	(\hat z_\alpha,\rho_\alpha)}{D\Psi(\hat
z_\alpha,\rho_\alpha)}\Psi(\hat
z_\alpha,\rho_\alpha),\\
	\rho_\beta 	&=\frac{DF(\hat
	z_\alpha,\rho_\alpha)}{D\Psi(\hat z_\alpha,\rho_\alpha)},
\end{split}
\end{equation}
where $D=\partial_\theta +\theta\partial_z$. The details of this
somewhat
unpleasant calculation are left to the reader.  From
(\ref{eq:coordtransfdualcurve}) we see that $\Cox$ contains the
structure sheaf
$\Oxhat$ of another $N=1$ super curve:  $\Oxhat$ is the sheaf of
$\Lambda$-algebras locally generated by $\hat z_\alpha,\rho_\alpha$.
We call
$\hat X=(X,\Oxhat)$ the {\em{dual curve}} of $(X,\Ox)$. We have
\begin{equation}\label{eq:extberbystructdual}
0\rightarrow\Oxhat\rightarrow\Cox
\overset{\Dchat}\rightarrow\Berxhat\rightarrow0,
\end{equation}
where $\Berxhat$ is the dualizing sheaf of the dual curve.  One can
show that
the dual curve of the dual curve is the original curve, thereby
justifying the
terminology.

\begin{exmpl}\label{exmpl:dualsplit}
We continue the discussion of split curves. In this case
(\ref{eq:coordtransfdualcurve}) becomes
\begin{equation}\label{eq:coordtransfdualcurvesplit}
\begin{split}
\hat z_\beta 	&= f(\hat z_\alpha),\\
	\rho_\beta 	&=\frac{\partial_{\hat z}f(\hat
	z_\alpha)}{B(\hat z_\alpha)}\rho_\alpha,
\end{split}
\end{equation}
So the dual split curve is $\Oxhatsplit=\Oxred\otimes\Lambda\mid
\mathcal
K\mathcal N\inv \otimes \Lambda$. The Berezinian sheaf for the dual
split curve
has generators that satisfy
\begin{equation}\label{eq:transfBersplitdual}
\hat f_\beta=B(z_\alpha) \hat f_\alpha.
\end{equation}
This means that $\Berxhat=\mathcal N\otimes \Oxhat=\mathcal N\otimes
\Lambda\mid \mathcal K\otimes \Lambda$. \qed
\end{exmpl}
\smallskip
A very useful geometric interpretation of the dual curve exists, cf.
\cite{DoRoSc:SupModSpaces,Schwarz:SuperanalogsSYMPLCONT}: the points
(i.e., the
$\Lambda$-points) of the dual curve correspond precisely to the
irreducible
divisors of the original curve and vice versa, as we will presently
discuss. In
 subsection \ref{ss:IntegraOnSupCurve} we will see that irreducible
divisors
are the limits that occur in contour integration on a super curve.

An irreducible divisor (for $\Ox$) is locally given by an even
function
$P_\alpha=z_\alpha-\hat z_\alpha -
\theta_\alpha \rho_\alpha\in \Ox(U_\alpha)$, where $\hat z_\alpha$
and
$\rho_\alpha$ are now respectively even and odd constants, i.e.,
elements of
$\Lambda$. Two divisors $P_\alpha,P_\beta$  defined on coordinate
neighborhoods
$U_\alpha$ and $U_\beta$, respectively,
are said to correspond to each other on the overlap if
\begin{equation}
P_\beta(z_\beta,\theta_\beta)=P_\alpha(z_\alpha,\theta_\alpha)
g(z_\alpha,\theta_\alpha), \quad g(z_\alpha,\theta_\alpha)\in
\Oxevstar(U_\alpha\cap U_\beta).\label{eq:corresponddiv}
\end{equation}
(If $R$ is a ring (or sheaf of rings)  ${R}^\times$ is the set of
invertible
elements.)

\begin{lem}\label{lem:roots} Let $(U, \mathcal O(U))$ be a $(1\mid
1)$
dimensional super domain with coordinates $(z,\theta)$ and let
$f(z,\theta)\in
\mathcal O(U)$. Then, with $D=\partial_\theta+\theta\partial_z$,
\begin{equation*}
f(z,\theta)=(z-\hat z -\theta
\rho)g(z,\theta)\quad\Leftrightarrow\quad f(\hat
z,\rho)=0,Df(\hat z,\rho)=0,
\end{equation*}
for $g(z,\theta)$  in $\mathcal O(U)$.
\end{lem}
Applying Lemma \ref{lem:roots} to (\ref{eq:corresponddiv}) we find
\begin{equation*}
\begin{split}
P_\beta(F(\hat z_\alpha,\rho_\alpha),
	\Psi(\hat z_\alpha,\rho_\alpha))
	&= F(\hat z_\alpha,\rho_\alpha)-\hat z_\beta-
		\Psi(\hat z_\alpha,\rho_\alpha)\rho_\beta=0,\\
 DP_\beta(F(\hat z_\alpha,\rho_\alpha),
	\Psi(\hat z_\alpha,\rho_\alpha))
&=DF(\hat z_\alpha,\rho_\alpha)-
		D\Psi(\hat z_\alpha,\rho_\alpha)\rho_\beta=0.
\end{split}
\end{equation*}
{}From this one sees that the parameters $(\hat
z_\alpha,\rho_\alpha)$ in the
local expression for an irreducible divisor transform as in
(\ref{eq:coordtransfdualcurve}), so they are $\Lambda$-points of the
dual
curve.

The $N=2$ super conformal super curve canonically associated to a
super curve
has a structure sheaf $\Cox$ that comes equipped with two sheaf maps
$\Dc$ and
$\Dchat$ with kernels the structure sheaves $\Ox$ and $\Oxhat$ of the
original
super curve and its dual. The intersection of the kernels is the
constant sheaf
$\Lambda$. The images of these maps are the dualizing sheaves
$\Berx$ and $\Berxhat$. In fact we can restrict $\Dc,\Dchat $ to the
subsheaves
$\Oxhat$ and $\Ox$, respectively, without changing the images. This
gives us
exact sequences
\begin{equation}\label{eq:DandDhatseq}
\begin{split}
0\rightarrow \Lambda\rightarrow &\Ox\overset{\hat
D}\rightarrow\Berxhat\rightarrow0,\\
0\rightarrow \Lambda\rightarrow
&\Oxhat\overset{D}\rightarrow\Berx\rightarrow0,
\end{split}
\end{equation}
with $D=\Dc |_{\Oxhat}$ and $\hat D=\Dchat|_{\Ox}$. Just as the sheaf
maps
$\Dc,\Dchat$ have local expressions as differential operators, also
their
restrictions are locally expressible in terms of differential
operators: if
$\{f_\alpha(z_\alpha,\theta_\alpha)\}$ is a section of $\Ox$ then the
corresponding section $\{(\Dchat f_\alpha)(\hat
z_\alpha,\rho_\alpha)\}$ of
$\Berxhat$ is given by
$$
\hat D f_\alpha(\hat
z_\alpha,\rho_\alpha)=
[(\partial_\theta+\theta\partial_z)f_\alpha]|_{z_\alpha=\hat
z_\alpha,\theta_\alpha=\rho_\alpha}.
$$
Similarly, if $\{\hat f_\alpha(\hat z_\alpha,\rho_\alpha)\}$ is a
section of
$\Oxhat$ then the corresponding section of $\Berx$ is
$$
{D}\hat{f}_\alpha(
z_\alpha,\theta_\alpha)=[(\partial_\rho+\rho\partial_{\hat
z})\hat f_\alpha]|_{\hat z_\alpha=
z_\alpha,\rho_\alpha=\theta_\alpha}.
$$
We summarize the relationships between the various sheaves and sheaf
maps in
the following commutative diagram (of sheaves of $\Lambda$-algebras):
\begin{equation}
\begin{CD}\label{eq:bigcd}
{}	@.	0	@.	0	@.	{}	@.	{}
\\
@.		@VVV		@VVV		@.		@.
\\
0  @>>>	\Lambda@>>>\Oxhat@>{ D}>>\Berx 	@>>>0	\\
@.		@VVV		@VVV		@\vert		@.
\\
0  @>>>	\Ox@>>>\Cox@>{\Dc}>>\Berx 	@>>>0	\\
@.		@V\hat{D}VV		@V\Dchat VV		@.
@.	\\
{}	@.  \Berxhat	@=	\Berxhat	@.{}
@.{}\\
@.		@VVV		@VVV		@.		@.
\\
{}	@.	0	@. 	0	@.	{}	@.	{}
\\
\end{CD}
\end{equation}
\medskip
We conclude this subsection with the remark that the dualizing sheaf
${\mathcal B\text{er}(\Cox)}$ of the super conformal super curve
$(X,\Cox)$
associated to a super curve $(X,\Ox)$ is trivial, making $(X,\Cox)$ a
super
analog of an elliptic curve or a Calabi-Yau manifold, cf.
\cite{DistNelson:SemiRigidSGra}. In fact, this statement is true for
any $N=2$
super curve $(X,\mathcal E)$ where $\mathcal E$ is an extension of
$\Berx$ by
the structure sheaf: if we have
\begin{equation*}
0 \to \Ox \to{\mathcal E} \to \Berx \to 0,
\end{equation*}
then $\mathcal E$ has local generators
$(z_\alpha,\theta_\alpha,\rho_\alpha)$
on $U_\alpha$, and on overlaps we get
\begin{equation} \label{eq:generalextensionberbystruct}
\begin{split}
z_\beta 	&= F_{\beta\alpha}(z_\alpha,\theta_\alpha),\\
\theta_\beta  	&= \Psi_{\beta\alpha}(z_\alpha,\theta_\alpha),\\
\rho_\beta 	&=
\Phi_{\beta\alpha}(z_\alpha,\theta_\alpha,\rho_\alpha)=
		\ber(J(z,\theta))\rho_\alpha + \phi_{\beta\alpha}
			(z_\alpha,\theta_\alpha),
\end{split}
\end{equation}
where $\ber(J(z,\theta))$ is the Berezinian of the super Jacobian
matrix of the
change of $(z,\theta)$ coordinates; this is precisely the transition
function
for $\Berx$, see \eqref{eq:transfBer}. Then  the super Jacobian
matrix
\begin{multline*}
J(z,\theta,\rho)=\ber
	\begin{pmatrix}
	\partial_z F 		&\partial_z \Psi 	&\partial_z
\Phi\\
	\partial_\theta F	&\partial_\theta \Psi
&\partial_\theta \Phi\\
	\partial_\rho F 	&\partial_\rho \Psi
&\partial_\rho \Phi\\
	\end{pmatrix}=
\ber
	\begin{pmatrix}
	\partial_z F 		&\partial_z \Psi 	&\partial_z
\Phi\\
	\partial_\theta F	&\partial_\theta \Psi
&\partial_\theta \Phi\\
	0		 	&0
&\partial_\rho \Phi\\
	\end{pmatrix}=    \\
	=\ber( J(z,\theta))/\partial_\rho \Phi=1,\quad
\end{multline*}
for all overlaps $U_\alpha\cap U_\beta$, and therefore $(X,\mathcal
E)$ has
trivial dualizing sheaf.

\subsection{Super Riemann surfaces.}\label{ss:SRS}
In this subsection we briefly discuss a special class of $N=1$ super
curves,
the super Riemann surfaces (SRS). This class of curves is studied
widely in the
literature because of its applications in super string theory, see
e.g.,
\cite{Fried:NoteString2DCFT,GidNelson:GeomSRS,LebrRoth:ModuliSRS,%
CraneRabin:SRSuniTeichm}. (Also the term
$\text{SUSY}_1$ curve is used,
\cite{Manin:GaugeFieldTheoryComplexGeom,Manin:Topicsnoncomgeom}, or
super
conformal manifold, \cite{RoSchVor:GeomSupConf}.)
{}From our point of view super Riemann surfaces are special because
irreducible
divisors and $\Lambda$-points can be identified and because there is
a
differential operator taking functions to sections of the dualizing
sheaf. Both
facts simplify the theory considerably. However, by systematically
using the
duality of
the $N=2$ super conformal curve one can extend results previously
obtained
solely for super Riemann surfaces to arbitrary super curves.

In the previous subsection we have seen that every $N=1$ super curve
$(X,\Ox)$
has
a dual curve $(X,\Oxhat)$. Of course it can happen that the
transition
functions of
$(X,\Ox)$ are identical to those of the dual curve $(X,\Oxhat)$. This
occurs if
the transition functions satisfy
\begin{equation}\label{eq:SRScondition}
DF(z_\alpha,\theta_\alpha)=
\Psi(z_\alpha,\theta_\alpha)D\Psi(z_\alpha,\theta_\alpha).
\end{equation}
If (\ref{eq:SRScondition}) holds then the operator
$D_\alpha=\partial_{\theta_\alpha}+\theta_\alpha\partial_{z_\alpha}$
transforms
as
\begin{equation}\label{eq:transfDSRS}
D_\beta=(D\Psi)\inv D_\alpha
\end{equation}
So in the situation of (\ref{eq:SRScondition}) the super curve
$(X,\Ox)$ carries a $(0\mid 1)$ dimensional distribution $D$ such
that $D^2$ is
nowhere vanishing (in fact $D^2=\partial_z$). A super curve carrying
such a
distribution is called a ($N=1$) super Riemann surface.  Equivalently
an $N=1$
super Riemann surface is a  ($N=1$) super curve that carries an odd
global
differential operator with nowhere vanishing square that takes values
in some
invertible sheaf.

Recall the Berezinian that occurs in the transformation law for
generators of
$\Berx$, (\ref{eq:transfBer}). It can be written in general as
$$
\ber\begin{pmatrix}	\partial_z F		&\partial_z \Psi\\
			\partial_\theta F	&\partial_\theta\Psi
			\end{pmatrix}=D(\frac{DF}{D\Psi}).
$$
Therefore if (\ref{eq:SRScondition}) holds we have
$\ber\begin{pmatrix}\partial_z F&\partial_z \Psi\\
			\partial_\theta F&\partial_\theta\Psi
			\end{pmatrix}=D\Psi$ so (\ref{eq:transfDSRS})
tells us that $D$ takes values in the dualizing sheaf $\Berx$.

So  super Riemann surfaces are self dual, as probably first noted in
\cite{DoRoSc:SupModSpaces}. More generally, the question then arises
what
happens  if the curves $(X,\Ox)$ and $(X,\Oxhat)$ are isomorphic, but
a priori
not with identical transition functions. We claim that also in this
case
the curve $(X,\Ox)$ is  a super Riemann surface.
 Indeed, the operator $\Dchat$ restricted to $\Ox$ takes values in
the
dualizing sheaf $\Berxhat$ of $\Oxhat$, as we have seen above. Using
the
isomorphism we can think of $\Dchat$ as a differential operator
taking values
in a sheaf isomorphic to the dualizing sheaf $\Berx$ on $\Ox$. Since
$\Dchat^2$
does not vanish we see that $(X,\Ox)$ is a super Riemann surface. Now
it is
known (and easy to see) that for any super Riemann surface there are
coordinates such that (\ref{eq:SRScondition}) holds. In these
coordinates the
transition functions of $(X,\Ox)$ and $(X,\Oxhat)$ are in fact equal.

The $N=2$ super conformal curve $(X,\Cox)$ associated to a SRS
$(X,\Ox)$ is
very
simple. Recall that $\Cox$ is an extension
\begin{equation*}
0\to\Ox\to\Cox\overset\epsilon\to\Berx \to 0.
\end{equation*}
where locally $\epsilon(z)=\epsilon(\theta)=0$ and
$\epsilon(\rho)=f$, with
$f$ a local generator of $\Berx$. For SRS there is a splitting
$e:\Berx\to\Cox$, given locally by $e(f)=\rho-\theta$. One needs to
use the
definition of a SRS to check that this definition makes global sense,
i.e.,
that $\rho-\theta$ transforms as a section of $\Berx$; for this see
\cite{Ra:oldnew}. In other words for a SRS the associated $N=2$ curve
has a
split structure sheaf:
$$
\Cox=\Ox\oplus \Berx.
$$
Note that not all SRS's are projected, so there are examples where
$\Cox$ is a
trivial extension but where $(X,\Ox)$ is not projected.

\subsection{Integration on super curves.}\label{ss:IntegraOnSupCurve}

Let us first recall the classical situation. On an ordinary Riemann
surface
$(X,\Oxred)$ we can integrate a holomorphic 1-form $\omega$ along a
contour
connecting two points $p$ and $q$ on $X$.
If the contour connecting $p$ and $q$
lies in a single, simply connected, coordinate neighborhood
$U_\alpha$ with
local coordinate $z_\alpha$ we can write
$\omega=d f_\alpha$, with $f_\alpha\in \Oxred(U_\alpha)$ determined
up to a
constant. The points $p,q$  are described by the irreducible divisors
$z_\alpha-p_\alpha$ and $z_\alpha-q_\alpha$. Then we calculate the
integral of
$\omega$ along the contour by
$\int_p^q\omega=f_\alpha(q_\alpha)-f_\alpha(p_\alpha)$.
Suppose next that $p$ and $q$ are in different coordinate
neighborhoods
$U_\alpha$ and $U_\beta$, with coordinates $z_\alpha,z_\beta$ related
by
$z_\beta=F(z_\alpha)$ on overlaps. Assume furthermore that the
contour
connecting them contains a point $r \in U_\alpha\cap U_\beta$. Then
we can
write $\omega=df_\alpha$ on $U_\alpha$, and $\omega=df_\beta$ on
$U_\beta$,
with $f_\alpha(z_\alpha)=f_\beta(F(z_\alpha))+ c_{\alpha\beta}$ on
overlaps,
where $c_{\alpha\beta}$ is locally constant on $U_\alpha\cap
U_\beta$. The
intermediate point $r$ can be described by two (corresponding)
irreducible
divisors $z_\alpha-r_\alpha$ and $z_\beta-r_\beta$. Then
$\int_p^q\omega=\int_p^r\omega+\int_r^q\omega=f_\beta(q_\beta)-
 f_\beta(r_\beta)+f_\alpha(r_\alpha) -f_\alpha(p_\alpha)$. This is
independent
of the intermediate point because the parameter $r_\alpha$ in the
irreducible
divisor $z_\alpha-r_\alpha$ transforms as a $\mathbb C\,$-point of
the curve:
we have $r_\beta=F(r_\alpha)$, and $f_\alpha(r_\alpha)-
f_\beta(r_\beta)=c_{\alpha\beta}$; therefore we can replace $r$ by
any other
intermediate point in the same connected component of $U_\alpha\cap
U_\beta$.
If $p$ and $q$ are not in adjacent coordinate neighborhoods we need
to
introduce more intermediate points.

So there are three crucial facts in the construction of the contour
integral of
holomorphic 1-forms on a Riemann surface: the parameter in  an
irreducible
divisor transforms as a point,  $d$ is an operator that produces from
a
function on $X$ a section of the dualizing sheaf  on $X$, and the
kernel of the
operator $d$ consists of the constants. We will find analogs of all
three facts
for super curves.

We have seen that for an $N=1$ super curve in general the parameters
in an
irreducible divisor correspond to a $\Lambda$-point of the dual
curve. Also the
sheaf map $D$ acting on the dual curve maps sections of $\Oxhat$ to
sections of
$\Berx$, see (\ref{eq:bigcd}).

This suggests that we define a {\it (super) contour}
$\Gamma=(\gamma,P,Q)$
 on $(X,\Ox)$ as an ordinary contour $\gamma$ on the underlying
topological
space $X$, together with two  irreducible divisors $P$ and $Q$ for
$(X,\Ox)$
such that the reduced divisors of $P$ and $Q$ are the endpoints of
$\gamma$. So
if
$$
P=z-\hat p-\theta\hat \pi,\quad Q=z-\hat q-\theta \hat\chi,
$$
then the corresponding $\Lambda$-points on the dual curve
$(X,\Oxhat)$ are
$(\hat p,\hat \pi)$, $(\hat q,\hat \chi)$, and $z=\hat
p^{\text{red}}$ and
$z=\hat q^{\text{red}}$ are the equations for the endpoints of the
curve
$\gamma$.  Then we define the integral of a section
$\{\omega_\alpha={D}\hat
f_\alpha\}$ of the dualizing sheaf on $(X,\Ox)$ along $\Gamma$  by
$$
\int_P^Q \omega=\int_P^Q  D \hat f=\hat f(\hat q,\hat \chi)-\hat
f(\hat p,\hat
\pi).
$$
Here we assume that the contour connecting $P$ and $Q$ lies in a
single simply
connected open set. If the contour traverses various  open sets we
need to
choose intermediate divisors on the contour,  as before.

A super contour $\Gamma$ is called {\it closed} if it is of
 the form $\Gamma=(\gamma,P,P)$, with the underlying
 contour $\gamma $ closed in the usual sense. Observe
 that the integral over $\Gamma$ is independent of the choice of $P$,
so
 we will omit reference to it.

The contour integration on $N=1$ super curves introduced here seems
to be new;
it is a nontrivial generalization of the contour integral on super
Riemann
surfaces, as described for instance in
\cite{Fried:NoteString2DCFT,McArthur:LineintegralsSRS,Rog:ContourSRS}.

 For
closed contours it agrees with the
integration theory described in
\cite{GaiKhudShvar:IntegrationSurSuperSpce,Khud:BVformoddsympl}.

We can also understand this integration procedure in terms of  the
contour
integral on the
$N=2$ super conformal super curve $(X,\mathcal{CO}_X)$, introduced by
Cohn,
\cite{Cohn:N=2SRS}. To this end
define on
$\Cox(U)\oplus\Cox(U)$ the sheaf map $(\Dc,\Dchat )$ by the local
componentwise
action of the differential operators $\Dc^\alpha$ and $\Dchat^\alpha$
as
before. Then the square of the operator $(\Dc,\Dchat )$ vanishes and
the
Poincar\'e Lemma holds for  $(\Dc,\Dchat)$:

\begin{lem}\label{lem:poincare}
Let $U$ be a simply connected open set on $X$ and let $(f,g)\in
\Cox(U)\oplus\Cox(U)$ such that $(\Dc,\Dchat )(f,g)=0$. Then there is
an
element $H\in \Cox(U)$, unique up to an additive constant, such that
$$(f,g)=(\Dc  H,\Dchat  H).$$
\end{lem}

Let then $\mathcal M(U)\subset \Cox(U)\oplus\Cox(U)$ be the subsheaf
of
$(\Dc,\Dchat )$-closed sections. Note that a section of
$\mathcal M$ looks in $U_\alpha$ like
$(f_\alpha,g_\alpha)=(f(z_\alpha,\theta_\alpha),g(\hat
z_\alpha,\rho_\alpha))$
and furthermore $f$ is a section of $\Berx$ and $g$ is a section of
$\Berxhat$.
This means that $\mathcal M$ globalizes to the direct sum
$\Berx\oplus
\Berxhat$.

So we get an  exact sequence of sheaves:
$$
0\to \Lambda\to \mathcal{CO}_X\overset(\Dc,\Dchat)\to
\mathcal M\rightarrow0.
$$
Now the sections of $\mathcal M$ are the objects on $\Cox$ that can
be
integrated. A contour for $\Cox$ is a triple $(\gamma, \mathcal
CP,\mathcal
CQ)$ where $\mathcal CP,\mathcal CQ$ are  two $\Lambda$-points of
$(X,\mathcal{CO}_X)$ with as reduced points the endpoints of the
contour
$\gamma$. Assume that  the contour lies in a single simply connected
open set
$U$. If  $\omega\in \mathcal M(U)$ then we can write $\omega=(\Dc
H,\Dchat
H)$ for some $H\in \Cox(U)$ and we put $\int_{\mathcal CP}^{\mathcal
CQ}\omega=H(\mathcal CQ)-H(\mathcal CP)$. Extension to more
complicated
contours as before.

Now start with a section $\{s_\alpha\}$ of $\Berx$ on $(X,\Ox)$. We
can lift it
to the section $\{(s_\alpha,0)\}$ of $\mathcal M $.  In particular
there is a
section $\{H_\alpha\}$ of
$\mathcal{CO}_X$  such that $s_\alpha=\Dc  H_\alpha$, $\Dchat
H_\alpha=0$. This
means that $\{H_\alpha\}$ is in fact a section of the subsheaf
$\Oxhat$. So in
specifying the $\Lambda$-points of  $(X,\mathcal{CO}_X)$ on the ends
of the
contour we have the freedom to shift along the fiber of the
projection $\hat
\pi:(X,\Cox)\to (X,\Oxhat)$. In other words we only need  to specify
$\Lambda$-points of the dual curve, or, equivalently, irreducible
divisors on
the original curve.

Therefore we can define the integral of a section $s=\{s_\alpha\}$ of
$\mathcal
 B\text{er}_X$ along a contour with $P,Q$ two irreducible divisors
for
$(X,\Ox)$ at the end point as follows. We choose two $\Lambda$-points
$\mathcal
CP$ and  $\mathcal CQ$ of $(X,\Cox)$ that project to the
$\Lambda$-points of
$(X,\hat{\mathcal{O}}_X)$ corresponding to $P,Q$. Then $\int_P^Q
s=H(\mathcal
CQ)-H(\mathcal CP)$
if $s_\alpha=\Dc H$ and $\Dchat H=0$. Again we are assuming here that
the
contour lies in a simply connected region and extend for the general
case using
intermediate points. One checks that this procedure of integrating a
section of
the dualizing sheaf on $(X,\Ox)$ using integration on $\Cox$ is the
same as we
had defined before.

\subsection{Integration on the universal cover.}

We consider from now on only holomorphic (compact, connected, $N=1$)
super
curves $(X,\Ox)$ of genus $g>1$. We fix a point $x_0\in X$ and
1-cycles $A_i,B_i, i=1,\dots ,g$ through $x_0$ with intersection
$A_i\cdot B_j=\delta_{ij}$, $A_i\cdot A_j=B_i\cdot B_j=0$ as usual.
Then the
fundamental group $\pi_1(X,x_0)$ is generated by the classes
$a_i,b_i$
corresponding to the loops $A_i,B_i$, subject solely to the relation
$a_1b_1a_1\inv b_1\inv a_2b_2a_2\inv b_2\inv \dots a_gb_ga_g\inv
b_g\inv=e$.

The universal cover of the super curve $(X,\Ox)$ is the open
superdisk
$\done=(D,\mathcal O_{\done})$ of dimension $(1\mid1)$, where
$D=\{z\in \mathbb
C\mid |z|<1\}$ and
$\mathcal O_{\done}=\mathcal O_D\otimes_{\mathbb C}\Lambda[\theta]$,
with
$\mathcal O_D$ the usual sheaf of holomorphic functions on the unit
disk. The
group $G$ of covering transformations of $(D,\mathcal O_{\done})\to
(X,\Ox)$ is
isomorphic to $\Pi_1(X,x_0)$ and each covering transformation $g$ is
determined
by its action on the global coordinates $(z,\theta)$ of $\done$.
Introduce
super holomorphic functions by
$$
F_g(z,\theta):=g\inv\cdot z , \quad\Psi_g(z,\theta):=g\inv\cdot
\theta.
$$
If $P_p$ is a $\Lambda$-point of $\done$, i.e., a homomorphism
$\mathcal O_{\done}\to \Lambda$, determined by $z\mapsto
 z_P\in \Lambda_0,\theta \mapsto \theta_P\in \Lambda_1$, then
 the action of $g$ in the covering group is defined by
$g\cdot P_p(f)=P_p(g\inv\cdot f)$. Then $z_P\mapsto
F_g(z_P,\theta_P)$ and
$\theta_P\mapsto \Psi_g(z_P,\theta_P)$. So
$\Lambda$-points transform as the coordinates under the covering
group.

Next consider irreducible divisors $P_d=z-\hat z_1
-\theta\hat\theta_1$,
$Q_d=z-\hat z_2 -\theta\hat\theta_2$.  We say that $g\cdot  P_d=Q_d$
as
divisors if we have the identity
$g\inv Q_d=P_d h(z,\theta)$ as holomorphic functions for some
invertible $h(z,\theta)$. By the same calculation as the one
following Lemma
\ref{lem:roots} we find that
\begin{equation*}
\begin{split}
\hat z_2&= F_g(\hat z_1, \hat \theta_1) + \frac {DF_g(\hat z_1, \hat
\theta_1)}{D\Psi_g(\hat z_1, \hat \theta_1)}\Psi_g(\hat z_1, \hat
\theta_1),\\
\hat \theta_2&= \frac {DF_g(\hat z_1, \hat \theta_1)}{D\Psi_g(\hat
z_1, \hat
\theta_1)}.
\end{split}
\end{equation*}
So irreducible divisors transform with the dual action, compare with
(\ref{eq:coordtransfdualcurve}).

There is a  parallel theory for the dual curve: we have a covering
$(D,\mathcal O_{\done})\to (X,\Oxhat)$, with covering group $\hat G$.
The dual covering group $\hat G$ is isomorphic to $G$ by a
distinguished
isomorphism: $g$ and $\hat g$ are identified if they give the same
transformation of the reduced disk. Their action on functions is in
general
different, however, unless we are dealing with a super Riemann
surface. In
fact, since duality interchanges irreducible divisors and
$\Lambda$-points on
the curve and its dual we see that the action
of $\hat g$ on the coordinates is dual to the transformation of
$g$:
\begin{equation*}
\begin{split}
\hat g\inv\cdot z&= F_g(z,\theta) + \frac {DF_g(z,\theta )}{D\Psi_g(
z ,
\theta )}\Psi_g( z ,  \theta ),\\
 \hat g\inv\cdot\theta &= \frac {DF_g(  z , \theta)}{D\Psi_g( z ,
\theta )}.
\end{split}
\end{equation*}

A function $f$ on $(X,\Ox)$ lifts to a function that is invariant
under the
covering group $G$ and similarly $\hat f$, a function on
$(X,\Oxhat)$, lifts to a function that is invariant under the dual
covering
group $ \hat G$. An irreducible divisor or a $\Lambda$-point on
$(X,\Ox)$ lifts
to an infinite set of divisors or points, one for each point on the
underlying
disk above the corresponding reduced point of $X$.

Let as before $x_0$ be a point on $X$ and $d_0$ a point on the disk
lying over
$x_0$.
Let $\gamma$ be a contour for integration on $(X,\Ox)$, so $\gamma$
consists of
a contour on
$X$ and two irreducible divisors at the endpoints.
The contour lifts to a unique contour on the disk starting at $d_0$
and the
irreducible divisors lift to unique
irreducible divisors for $(D,\mathcal O_{\done})$ that reduce to
$d_0$ and the
endpoint on the disk, respectively.
Also we can pull back sections of $\Berx$ to $(D,\mathcal O_{\done})$
and
calculate integrals on $(X,\Ox)$ by
lifting to $(D,\mathcal O_{\done})$.
Since $D$ is simply connected this is a great simplification.
For instance any integral over a closed contour is zero.

Similar considerations apply to the $N=2$ curve $(X,\Cox)$ and its
universal covering space $D^{1|2}$ and covering group $\mathcal G$.
Of course $D^{1|2}$ is the $N=2$ curve canonically associated to the
$N=1$
curve
$D^{1|1}$ as in subsection \ref{ss:DualN=2curves}, and the lifts of
$f \in \Ox$
to $D^{1|2}$ via
either $(X,\Cox)$ or $D^{1|1}$ as intermediate space coincide.

\subsection{Sheaf cohomology for super
curves.}\label{ss:Sheafcohomology}

Our super curves are in fact families of curves over the base scheme
$(\bullet,\Lambda)$, with $\Lambda$ the Grassmann algebra of
nilpotent
constants. This means that for any coherent sheaf the cohomology
groups are
finitely generated $\Lambda$-modules, but they are not necessarily
free. This
means in particular that standard classical theorems, like the
Riemann-Roch
theorem, do not hold in general in our situation. (See for instance
\cite{Hodg:ProblFieldSRS}.)

The basic facts about  sheaf cohomology of families of super curves
are
completely parallel to the classical theory (explained for instance
in
\cite{Kempf:AbInt}). For a coherent locally free sheaf  $\mathcal L$
there
exist $\Lambda$-homo\-mor\-phisms $\alpha:F\to G$, with $F, G$ free
finite rank
$\Lambda$-modules, that calculate the cohomology. More precisely, for
every
$\Lambda$-module $M$ we have an exact sequence
\begin{equation}\label{eq:calccohom}
0\to H^0(X,\mathcal L\otimes M)\to F\otimes M\overset{\alpha\otimes
1_{M}}
\to G\otimes M\to H^1(X,\mathcal L\otimes M)\to0.
 \end{equation}
Recall from Example \ref {exmpl:split} that for any  sheaf of
$\Ox$-modules
$\mathcal F$ we have  an associated split sheaf  $\mathcal
F^{\text{split}}=\mathcal F\otimes_\Lambda \Lambda/\mathfrak m$.
Therefore, if we choose $M=\Lambda/\mathfrak m$, the sequence
(\ref{eq:calccohom}) calculates the cohomology groups  of the split
sheaf
$\mathcal L\spl$. (These cohomology groups are
$\mathbb Z_2$-graded vector spaces over $\Lambda/\mathfrak m=\mathbb
C$.)
Without loss of generality one can choose the homomorphism
$\alpha:F\to G$ such
that $\alpha\spl=\alpha\otimes 1_{\Lambda/\mathfrak m}$ is
identically zero.
This means that
$H^0(X,\mathcal L)$ (respectively $H^1(X,\mathcal L)$) is a submodule
(resp. a
quotient module) of a free $\Lambda$-module of rank $\dim
H^0(X,\mathcal
L\spl)$ (resp. of rank $\dim H^1(X,\mathcal L\spl)$).

We are interested in the question when the $H^i(X,\mathcal L)$ are
free. The
idea is to check this by an inductive procedure, starting with the
free
cohomology of $\mathcal L\spl$.
We have for every $j=1,\dots,n-1$ the split exact sequence
\begin{equation}\label{eq:seqlambda}
0\to \mathfrak m^j/\mathfrak m^{j+1}\to \Lambda/\mathfrak m^{j+1}\to
\Lambda/\mathfrak m^j\to 0.
\end{equation}
 Since $\mathfrak  m^j/\mathfrak m^{j+1}\otimes_\Lambda \mathcal
L=\mathfrak
m^j/\mathfrak m^{j+1}\otimes_{\mathbb C} \mathcal L\spl$,
$\Lambda/\mathfrak
m^i\otimes_\Lambda \mathcal L=\mathcal L/\mathfrak m^i \mathcal L$
and
$\mathcal L$ is flat over $\Lambda$ we obtain by tensoring with
$\mathcal L$
and taking cohomology the exact sequence ($\Lambda^j=\mathfrak
m^j/\mathfrak
m^{j+1}$)
\begin{equation}\label{eq:longexactcohom}
\begin{aligned}
0 &\to\Lambda^j\otimes_{\mathbb C}H^0(X,\mathcal L\spl)&
  &\to H^0(X, \mathcal L/\mathfrak m^{j+1}\mathcal L) &
  &\to  H^0(X, \mathcal L/\mathfrak m^{j}\mathcal L) &
  &\overset{q^j}\to \\
 &\overset{q^j}\to \Lambda^j\otimes_{\mathbb C}H^1(X,\mathcal L\spl)
&
 &\to H^1(X, \mathcal L/\mathfrak m^{j+1}\mathcal L) &
 &\to  H^1(X, \mathcal L/\mathfrak m^{j}\mathcal L) &
 &\to 0.
\end{aligned}
\end{equation}
If $H^0(X, \mathcal L/\mathfrak m^{j}\mathcal L)$ and $H^1(X,
\mathcal
L/\mathfrak m^{j}\mathcal L)$ are free $\Lambda/\mathfrak
m^j$-modules, then
the module
$H^0(X, \mathcal L/\mathfrak m^{j+1}\mathcal L)$ is free over
$\Lambda/\mathfrak m^{j+1}$ iff the connecting map $q^j$ in
(\ref{eq:longexactcohom}) is zero iff $H^1(X, \mathcal L/\mathfrak
m^{j+1}\mathcal L)$ is free as $\Lambda/\mathfrak m^{j+1}$-module
(see \cite{Kempf:AbInt}, Lemma 10.4).
The relation between the connecting homomorphisms $q^j$ and the
homomorphism $\alpha$ that calculates cohomology is as follows: if we
assume as
above $\alpha\spl$ is zero then $q^1=\alpha\otimes
1_{\Lambda/\mathfrak m^2}$.
More generally, if $q^1=q^2=\dots=q^{j-1}=0$
then $q^j=\alpha\otimes 1_{\Lambda/{\mathfrak m}^{j+1}}$.

More concretely, we can assume that
$\alpha$ is a matrix of size rank $G\times \text{rank }F$ and the
$q^j$ are
quotients of this matrix by $\mathfrak m^{j+1}$. Then the cohomology
of
$\mathcal L$ is the kernel and cokernel of the matrix $\alpha$, and
the
cohomology is free iff $\alpha$ is identically zero.

If now $\mathcal L$ is an invertible sheaf, $\mathcal{L}\spl$ obeys a
super
Riemann-Roch relation and in case of free cohomology we get ($h^i=
\text{rank } H^i$):
\begin{equation} \label{superRR}
h^0(X,\mathcal{L}) - h^1(X,\mathcal{L}) = (\deg \mathcal{L} + 1-g\mid
\deg \mathcal{L} + \deg \mathcal{N} + 1-g),
\end{equation}
where  $\Ox\spl=\Oxred\mid \mathcal N$. We can relate by Serre
duality the
cohomology groups of $\mathcal L$ and $\mathcal L^*\otimes \Berx$,
see Appendix
\ref{app:dualSerredual}. In particular, $H^0(X,\mathcal L^*\otimes
\Berx)$ is
free iff $H^1(X,\mathcal L)$ is.

We summarize the discussion in this subsection in the following
theorem.

\begin{thm}\label{thm:freeness}
Let $\mathcal L$ be an invertible $\Ox$-sheaf. Then $H^0(X,\mathcal
L)$
(respectively $H^1(X,\mathcal L)$) is a submodule (respectively a
quotient
module) of a free $\Lambda$-module of rank $\dim H^0(X,\mathcal
L\spl)$
(respectively of rank $\dim H^1(X,\mathcal L\spl)$). Furthermore
\begin{align*}
H^0(X,\mathcal L)\text{ is a free $\Lambda$-module}
&\Longleftrightarrow
H^1(X,\mathcal L)\text{ is free},\\
&\Longleftrightarrow H^0(X,\mathcal L^*\otimes\Berx)\text{ is
free},\\
&\Longleftrightarrow H^1(X,\mathcal L^*\otimes\Berx)\text{ is free},
\end{align*}
in which case the rank of $H^i(X,\mathcal L)$ is equal to $\dim
H^i(X,\mathcal
L\spl)$.
\end{thm}

\subsection{Generic SKP curves.}

\begin{defn} An { \it  SKP curve} is a super curve $(X,\Ox)$ such
that the
split sheaf $\Oxsplit$ is of the form
$$
\Oxsplit=\Oxred\mid \mathcal N,
$$
where $\mathcal N$ is an invertible $\Oxred$-module of degree zero.
If
$\mathcal N\ne\Oxred$ then $(X,\Ox)$ is called a {\it generic SKP
curve}.
\qed
\end{defn}

We will discuss in subsection \ref{ss:Krichever}  a Krichever map
that
associates to an invertible sheaf on a super curve $(X,\Ox)$ (and
additional
data) a point $W$ of an infinite super Grassmannian. If this point
$W$ belongs
to the {\it big cell} (to be defined below) we obtain a solution of
the super
KP hierarchy. For $W$ to belong to the big cell it is necessary that
$(X,\Ox)$
is an SKP curve. The generic SKP curves enjoy simple cohomological
properties.

\begin{thm}\label{thm:cohomcurve}
Let $(X,\Ox)$ be  a generic SKP curve. Then the cohomology groups of
the
sheaves
$\Ox, \Berx$  are free $\Lambda$-modules. More precisely:
\begin{alignat*}{2}
H^0(X,\Ox) &= \Lambda \mid 0,
	&\qquad H^1(X,\Ox)&= \Lambda^g\mid \Lambda ^{g-1},\\
H^0(X,\Berx) &= \Lambda^{g-1} \mid \Lambda^g,
	&\qquad H^1(X,\Berx)&= 0\mid\Lambda.
\end{alignat*}
\end{thm}

\begin{proof}
Since $\mathcal N$ has no global sections,  $H^0(X, \Oxsplit)=\mathbb
C\mid 0$
consists of the constants only.
Now by definition of a curve over $(\bullet,\Lambda)$ we have an
inclusion
$0\to\Lambda\to \Ox$, so $H^0(X,\Ox)$ contains at least the constants
$\Lambda$.  By Theorem \ref{thm:freeness} then $H^0(X,\Ox)$ must be
equal to
$\Lambda\mid 0$. Again using Theorem \ref{thm:freeness} then also
$H^1(X, \Ox)$
and the cohomology of $\Berx$ will be free, and the rest of the
theorem follows
from the properties of the split sheaves, see Examples
\ref{exmpl:split},
\ref{exmpl:splitber}.
\end{proof}

\begin{rem}
It is not true that all invertible $\Ox$-sheaves for a generic SKP
curve have
free cohomology. For instance, consider a sheaf $\mathcal L$ with
$\mathcal
L\spl=\Oxsplit$, but $\mathcal L\neq \Ox$. Then, for a covering
$\{U_\alpha\}$
of $X$, the transition functions of $\mathcal L$ will have the form
$g_{\alpha\beta}=1 + f_{\alpha\beta}(z,\theta)$, with
$f_{\alpha\beta}(z,\theta)=0$ modulo the maximal ideal $\mathfrak m$
of
$\Lambda$. Let then $I\subset\Lambda$ be the ideal of elements that
annihilate
all $f_{\alpha\beta}$. Then we have $H^0(X,\mathcal L)=I$ and is in
particular
not free.
\qed
\end{rem}

\subsection{Riemann bilinear relations.}
Let us call sections of $\Berx$ and $\Berxhat$ holomorphic
differentials (on
$(X,\Ox)$ and $(X,\Oxhat)$ respectively).
We will in this subsection introduce analogs of the classical
bilinear
relations for holomorphic differentials.
\begin{thm}
\label{thm:bilinear}
 Let $(X,\Ox)$ be a super curve and let $\omega$, $\hat \omega$ be
holomorphic
differentials on $(X,\Ox)$ and $(X,\Oxhat)$ respectively. Let
$\{a_i,b_i\}$ be a standard symplectic basis for $H_1(X,\mathbb Z)$.
Then
$$
\sum_{i=1}^g\oint_{a_i}\omega\oint_{b_i} \hat
\omega=\sum_{i=1}^g\oint_{a_i}\hat\omega\oint_{b_i}  \omega.
$$
\end{thm}

Note that we  think here of closed contours on the underlying
topological space
$X$ as closed super contours on either $(X,\Ox)$ or on $(X,\Oxhat)$.

\begin{proof}
The argument is clearest using the $N=2$ curve $(X,\Cox)$ and its
universal
covering superdisk $D^{1\mid2}$;
this way only one universal covering group $\mathcal G$ appears
instead of both
$G$ and $\hat{G}$.
Choose any holomorphic differentials  $\omega$ on $X$ and
$\hat{\omega}$ on
$\hat{X}$, and lift them
to sections $(\omega,0)$ and $(0,\hat{\omega})$ of $\mathcal M$ on
$(X,\Cox)$.
Lifting further to $D^{1\mid2}$, let $\Omega$ be an antiderivative of
$(0,\hat{\omega})$, so that
$(\Dc\Omega,\Dchat\Omega) = (0,\hat{\omega})$.
The crucial point is that $(\Omega\omega,0)$ is itself a section of
$\mathcal
M$, because $\Dc(\Omega\omega)=0$.
This could not have been achieved using only differentials from $X$.
As per the standard argument, we integrate this object around the
polygon
obtained by cutting open $(X,\Cox)$.
To form this polygon, fix arbitrarily one vertex $P$ (a
$\Lambda$-point of the
$N=2$ disk $D^{1\mid2}$)  and let the other vertices be
$a_1^{-1}P,\;b_1^{-1}a_1^{-1}P,\ldots,a_gb_g^{-1}a_g^{-1}\cdots
b_1a_1b_1^{-1}a_1^{-1}P$, where $a_i,b_i$ are
the generating elements of $\mathcal G$.
The vertices are the endpoints of super contours whose reduced
contours are the
sides of the usual polygon bounding a
fundamental region for $\mathcal G$.

These contours project down to any of $X,\hat{X},(X,\Cox)$ as closed
loops
generating the homology; integrating
a differential lifted from any of these spaces along a side of our
polygon will
yield the corresponding period.
Labeling the sides of the polygon with generators of $\mathcal G$ as
usual,
neighborhoods of the sides labeled $a_i$
are identified with each other by $b_i$ and vice versa.
Then we have
\begin{multline*}
0 = \oint \Omega(\omega,0) =
\sum_{i=1}^g \left[ \int_{a_i} \Omega(\omega,0) - \int_{a'_i}
\Omega(\omega,0)
\right]
+ \\ +\sum_{i=1}^g \left[ \int_{b_i} \Omega(\omega,0) - \int_{b'_i}
\Omega(\omega,0) \right].
\end{multline*}
In the first sum, the two integrals are related by the change of
variables
given by $b_i$; the differential
$(\omega,0)$ is invariant under this covering transformation while
$\Omega$
changes by the
$b_i$-period of $\hat{\omega}$.
The second sum is simplified in the same manner, with the result
\begin{equation*}
\sum_{i=1}^g \left[ \int_{ {a}_i}\omega \int_{b_i} \hat{\omega} -
\int_{a_i}
\hat{\omega} \int_{{b}_i} \omega \right] = 0.
\end{equation*}
\end{proof}

\subsection{The period map and cohomology.}

The commutative diagram \eqref{eq:bigcd} gives a commutative diagram
in
cohomology that partly reads:
\begin{equation}
\begin{CD}\label{eq:bigcdcohom}
{}			@.	H^0(X,\Berxhat)	@=
H^0(X,\Berxhat)	\\
@.		@V\operatorname{p\hat er}VV		@V{\hat q}VV
\\
H^0(X,\Berx)  	@>{\operatorname{per}}>> H^1(X,\Lambda)
		@>{\operatorname{r\hat ep}}>> H^1(X,\Oxhat)	\\
@|		@V\operatorname{rep}VV		@.	\\
H^0(X,\Berx) @>{q}>>H^1(X,\Ox)@.
\end{CD}
\end{equation}

Let  $\{a_i,b_i; i=1,\dots, g\}$ be a symplectic basis for
$H_1(X,\mathbb Z)$
and let $\{a_i^*,b_i^*; i=1,\dots, g\}$ be a dual basis for
$H^1(X,\mathbb Z)$
and also for $H^1(X,\Lambda)$. We will use Serre
duality (see Appendix \ref{appss:SerredualSupermanifold}) to identify
$H^1(X,\Ox)$ and $H^1(X,\Oxhat)$ with the duals of $H^0(X, \Berx)$
and $H^0(X,
\Berxhat)$.

\begin{lem}\label{lem:perrep}
The maps $\operatorname{per}$, $\operatorname{p\hat er}$,
$\operatorname{rep}$
and $\operatorname{r\hat ep}$ are explicitly given by
\begin{align*}
\operatorname{per}(\omega)	&=\sum_{i=1}^g(\oint_{a_i} \omega)
a_i^*
		+\sum_{i=1}^g(\oint_{b_i} \omega) b_i^*,\\
\operatorname{p \hat er}(\hat\omega)	&=\sum_{i=1}^g(\oint_{a_i}
\hat\omega)
	a_i^* +\sum_{i=1}^g(\oint_{b_i} \hat\omega) b_i^*,\\
\operatorname{rep}(\sigma)(\omega)
&=\sum_{i=1}^g\alpha_i(\oint_{b_i}
\omega)-\sum_{i=1}^g \beta_i(\oint_{a_i} \omega),\\
\operatorname{r\hat ep} (\sigma)(\hat\omega)
&=\sum_{i=1}^g\alpha_i(\oint_{b_i}
\hat\omega) -\sum_{i=1}^g \beta_i(\oint_{a_i} \hat\omega),
\end{align*}
where $\omega\in H^0(X,\Berx)$, $\hat\omega\in H^0(X,\Berxhat)$ and
$\sigma=
\sum_{i=1}^g \alpha_i a_i^*+\beta_i b_i^*\in H^1(X,\Lambda)$.
\end{lem}

\noindent If we introduce
a basis $\{\omega_\alpha,\alpha=1,\dots,g-1\mid w_j,j=1,\dots,g\}$ of
$H^0(X,\Berx)$ we obtain the {\it period matrix} associated to
$\per$:
$$
\Pi=\begin{pmatrix}
\oint_{a_i}\omega_\alpha&\oint_{a_i}w_j\\
\oint_{b_i}\omega_\alpha&\oint_{b_i}w_j
\end{pmatrix},
$$
where $i,j$ run from $1$ to $g$ and $\alpha$ runs from $1$ to $g-1$.

For the split curve  we have $H^0(X,\Oxhatsplit)=\mathbb C\mid\mathbb
C^{g-1}$
and the map
$$
D:H^0(X,\Oxhatsplit)\to H^0(X,\Berxsplit)
$$
has as image a $g-1$ dimensional, even subspace of exact
differentials. For
these elements the periods vanish, and one finds the reduction mod
$\mathfrak
m$ of $\Pi$ is given by
$$
\Pi\spl=\begin{pmatrix}
0&\Pi\red(a)\\
0&\Pi\red(b)
\end{pmatrix},
$$
where $\Pi\red=\begin{pmatrix}\Pi\red(a)\\ \Pi\red(b)\end{pmatrix}$
is the
classical period matrix of the underlying curve $(X,\Oxred)$. By
classical
results we can choose the basis of holomorphic differentials on the
reduced
curve so that $\Pi\red(a)=1_g$. From this it follows that we can also
choose in
$H^0(X,\Berx)$ a basis such that $\oint_{a_i}w_j=\delta_{ij}$ and so
that the
period matrix takes the form
\begin{equation}\label{eq:normperiodmatrix}
\Pi=\begin{pmatrix}
0&1_g\\
Z_{o}&Z_e
\end{pmatrix}.
\end{equation}
Note that $\Pi$ is not uniquely determined by the conditions we have
imposed:
we are still allowed to change $\Pi\mapsto \Pi^\prime=
\begin{pmatrix}
0&1_g\\
Z_{o}G&Z_e + Z_o \Gamma
\end{pmatrix}
$, corresponding to a change of basis of $H^0(X,\Berx)$ by an even
invertible
matrix
$\begin{pmatrix}
G&\Gamma\\
0&1_g
\end{pmatrix}$ of size $g-1\mid g\times g-1\mid g$.

Using the same basis we see that $\operatorname{rep}$ has matrix
$$
\begin{pmatrix}
\oint_{b_i}\omega_\alpha&-\oint_{a_i}\omega_\alpha\\
\oint_{b_i}w_j&-\oint_{a_i}w_j
\end{pmatrix}=\Pi^tI=\begin{pmatrix}
Z_o^t&0\\
Z_e^t&-I_g
\end{pmatrix}, \quad I=\begin{pmatrix}
0&-1_g\\
1_g&0
\end{pmatrix}.
$$
Again this matrix is not entirely determined by our choices.

{}From the commutativity of the diagram \eqref{eq:bigcdcohom} we see
that
the matrix of the map $q$ is given by
\begin{equation}\label{eq:matrixconnectinghom}
Q=\Pi^t I \Pi=\begin{pmatrix}
0&Z_o^t\\
-Z_o&Z_e^t-Z_e
\end{pmatrix}.
\end{equation}

In general, the structure sheaf $\Oxhat$ and dualizing sheaf
$\Berxhat$ of the
dual curve will not have free cohomology, so that we cannot represent
the maps
$\operatorname{rep}$, $\operatorname{r\hat ep}$ and $\hat q$ by
explicit
matrices.

The nonfreeness of the cohomology of $\Oxhat$ and $\Berxhat$
is determined by the odd component $Z_o$ of the period matrix, see
\eqref{eq:normperiodmatrix}.
Recall that $\Oxhatsplit=\Oxred \mid \mathcal K\mathcal N\inv$ and
$\Berxhatsplit=\mathcal N\mid \mathcal K$
(see Example \ref{exmpl:dualsplit}) and hence
\begin{alignat*}{2}
H^0(X,\Oxhatsplit)&= \mathbb C\mid \mathbb C^{g-1}, &
H^1(X,\Oxhatsplit)&=
\mathbb C^g\mid 0,\\
H^0(X,\Berxhatsplit)&= 0\mid \mathbb C^{g}, & H^1(X,\Berxhatsplit)&=
\mathbb
C^{g-1}\mid \mathbb C.
\end{alignat*}
{}From the diagram \eqref{eq:bigcd} we extract in cohomology, using
that
the map $H^1(X, \Berx) \to H^2(X,\Lambda)=\Lambda\mid 0$ is an (odd!)
isomorphism and $H^0(X, \Lambda)=\Lambda$,
\begin{equation}\label{eq:period-seq}
0 \to   \frac{H^0(X, \Oxhat)}{\Lambda}\overset{D}\to  H^0(X, \Berx)
\overset
{\per}\to  H^1(X,\Lambda) \to H^1(X, \Oxhat) \to  0
\end{equation}
so that the period map
has as kernel $H^0(X, \Oxhat)$ mod constants and as cokernel $H^1(X,
\Oxhat)$.
Therefore $\per$ is essentially one of the
 homomorphisms that  calculate cohomology introduced in
 subsection \ref{ss:Sheafcohomology}. We can even be more
 explicit: if  $\{\omega_\alpha\mid w_j\}$ is the (partially)
normalized basis
of
holomorphic differentials as above the homomorphism $\per$ maps
 the submodule generated by the $w_j$ isomorphically to a free rank
 $g$ summand of $H^1(X,\Lambda)$.
 This is irrelevant for the calculation of cohomology,
so we can replace the sequence \eqref{eq:period-seq} by
\begin{equation}\label{eq:period-seqsimple}
0 \to   H^0(X, \Oxhat)/\Lambda \overset{D}\to  \Lambda^{g-1}
\overset{Z_o}\to
  \Lambda^g \to H^1(X, \Oxhat) \to  0,
\end{equation}
and $H^0(X, \Oxhat)$ mod constants is the kernel of $Z_o$, whereas
the cokernel
of $Z_o$ is $H^1(X, \Oxhat)$.

Similarly, the cohomology of $\Berxhat$ is calculated  by the
sequence
\begin{multline}\label{eq:repiod-seq}
0 \to   H^0(X, \Berxhat) \overset{\operatorname{p\hat er}}\to  H^1(X,
\Lambda)
\overset{\operatorname{rep}}\to  H^0(X,\Berx)^* \\
\to H^1(X, \Berxhat) \to  \Lambda\to 0
\end{multline}
The image of a holomorphic differential  $\hat\omega$ in
$H^1(X,\Lambda)$ is
then a vector
$\operatorname{p\hat er}(\hat\omega)=\begin{pmatrix} a(\hat\omega)\\
b(\hat\omega)\end{pmatrix}$, where $a(\hat\omega)$ and
$b(\hat\omega)$ are the
vectors of $a$ respectively $b$ periods of $\hat\omega$.  By
exactness of
\eqref{eq:repiod-seq} we have $\operatorname{rep}\circ
\operatorname{p\hat
er}=0$, or, using bases,
$$
\begin{pmatrix}
Z_o^t&0\\
Z_e^t&-I_g
\end{pmatrix}\begin{pmatrix} a(\hat\omega)\\
b(\hat\omega)\end{pmatrix}=0
$$
This means that  the vector $b(\hat\omega)$ of $b$ periods is
(uniquely)
determined by the $a$ periods: $b(\hat\omega)= Z_e^t a(\hat\omega)$,
and the
vector of $a$ periods is constrained by the equation $Z_o^t
a(\hat\omega)=0$.
The submodule of $H^1(X, \Lambda)$ generated by the elements $b_i^*$
maps under
$\operatorname {rep}$ isomorphically to a free rank $0\mid g$ summand
of
$H^1(X,\Berx)^*$, so that for the calculation of cohomology we can
simplify
\eqref{eq:repiod-seq} to
\begin{equation}\label{eq:repiod-seqsimple}
0 \to   H^0(X, \Berxhat) \to  \Lambda^{g}  \overset{Z_o^t}\to
\Lambda^{g-1}
\to H^1(X, \Berxhat) \to  \Lambda\to 0.
\end{equation}
We summarize the results on the cohomology of the dual curve in the
following
Theorem.
\begin{thm}
Let $(X,\Ox)$ be a generic SKP curve with odd period matrix $Z_o$.
Then
$$
H^0(X,\Oxhat)/\Lambda\simeq \operatorname{Ker}(Z_o),\quad
H^1(X,\Oxhat)\simeq\operatorname{Coker}(Z_o).
$$
Furthermore  $H^0(X,\Berxhat)\simeq \operatorname{Ker}(Z_o^t)$ and
$\operatorname{Coker}(Z_o^t)$ is a submodule of $H^1(X,\Berxhat) $
such that
$$H^1(X,\Berxhat)/\operatorname{Coker}(Z_o^t)\simeq\Lambda.
$$
\end{thm}
\subsection{$\Cox$ as extension of $\Berx$.}
We discuss in this subsection, for generic SKP curves, the structure
of $\Cox$
as extension of $\Berx$ and the relation with free cohomology and the
projectedness of the curve $(X,\Ox)$.

{}From the sequence (\ref{eq:extberbystruct}) that defines $\Cox$ we
obtain in
cohomology
\begin{gather}\label{eq:cohomn=2}
\begin{aligned}
0 &\to	H^0(X,\Ox)& 	&\to H^0(X, \Cox) & &\to H^0(X, \Berx) &
&\overset{q}\to
\\
{}&\overset{q}\to H^1(X,\Ox)&&\to H^1(X, \Cox) &&\to  H^1(X,
\Berx)&&\to 0
\end{aligned}
\end{gather}
The cohomology of the sheaves $\Ox, \Berx$ is given by Theorem
\ref{thm:cohomcurve}. By Theorem \ref{thm:freeness} (or its extension
to rank
two sheaves) $H^0(X,\Cox)$ is a submodule of
a $\Lambda^{g+1}\mid\Lambda^{g-1}$ and $H^1(X,\Cox)$ is a quotient of
a $\Lambda^{g+1}\mid \Lambda^{g-1}$. We see from this that the
cohomology of
$\Cox$ is free if and only if $q$ is the zero map.

To describe the map $q$ in more detail we need to recall some facts
about
principal parts and extensions, (see e.g., \cite{Kempf:AbInt}). For
any
invertible sheaf $\mathcal L$ let $\underline { \mathcal
Rat}(\mathcal L)$ and
$\underline { \mathcal Prin}(\mathcal L)$ denote the sheaves of
rational
sections and principal parts for $\mathcal L$ and denote by $\mathcal
Rat(\mathcal L)$ and
$ \mathcal Prin(\mathcal L)$ their $\Lambda$-modules of  global
sections. Then
the cohomology of $\mathcal L$ is calculated by
$$
0\to H^0(X,\mathcal L)\to \mathcal Rat(\mathcal L)\to
\mathcal Prin(\mathcal L)\to H^1(X,\mathcal L)\to 0.
$$
In particular we can represent a class $\alpha\in H^1(X,\mathcal L)$
as a
principal part $p=\sum p_{x}$, where $p_{x}\in \mathcal Rat(\mathcal
L)/\mathcal L_{x}$, for $x\in X$.

If $\alpha\in H^1(X,\mathcal L)$ and $\omega\in H^0(X,\mathcal M)$,
for some
other invertible sheaf $\mathcal M$, then we can define the {\it cup
product}
$\omega\cup \alpha$ by representing $\alpha$ by a principal part $p$
and
calculating the principal part $\omega p= \sum \omega_{x}p_{x}$ in
$\mathcal
Prin(\mathcal M\otimes \mathcal L)$; the image of $\omega p$ in
$H^1(X,
\mathcal M\otimes \mathcal L)$ is then by definition $\omega\cup
\alpha$.

We want to understand the kernel of the cup product with $\omega \in
H^0(X,
\mathcal M)$ in case $\omega$ is odd and free (i.e., linearly
independent over
$\Lambda$). In this case there will be for any invertible sheaf
$\mathcal L$
sections that are immediately annihilated by $\omega$; let therefore
$\Ann(\mathcal L,\omega)\subset \mathcal L$ be the subsheaf of such
sections.
Putting $\mathcal L_\omega=\mathcal L/\Ann(\mathcal L,\omega)$, we
get, because
$\omega^2=0$, the exact sequence
\begin
{equation}
\label{eq:cupsequence}
0\to \mathcal L_\omega\overset{\omega}\to \Ann(\mathcal L\otimes
\mathcal
M,\omega)\to Q\to 0
\end
{equation}
Locally, in an open set $U_\alpha\subset X$, we have $\mathcal
L(U_\alpha)=\Ox(U_\alpha)l_\alpha$, $\mathcal
M(U_\alpha)=\Ox(U_\alpha)m_\alpha$ and we write
$\omega=\omega_\alpha(z,\theta)m_\alpha$, with
$\omega_\alpha=\phi_\alpha+\theta_\alpha f_\alpha$. Then
$f_{\alpha}^{\text{red}}$ is a regular function on $U_\alpha$ with
some divisor
of zeros $D_f=\sum n_i q_i$. Some of the $q_i$ may also be zeros of
(the lowest
order part of) $\phi_\alpha$ and there will be a maximal
$g_\alpha(z,\theta) \in \Ox(U_\alpha)_{\bar 0}$ (here
$\Ox(U_\alpha)_{\bar 0}$
is the module of even sections) such that
$$
\omega_\alpha(z,\theta)=\omega_\alpha(z,\theta)^\prime
g_\alpha(z,\theta)
$$
with $g_\alpha^{\text{red}}$ a regular function with divisor of zeros
$D_g$ (on
$U_\alpha$) satisfying $0\le D_g\le D_f$.
Then $\Ann(\mathcal L\otimes \mathcal M,\omega)(U_\alpha)$ is
generated by
$\omega_\alpha(z,\theta)^\prime l_\alpha\otimes m_\alpha$ and we see
that
$Q$ is a torsion sheaf: $Q$ is killed by the invertible sheaf
generated locally
by the even invertible rational function $g_\alpha(z,\theta)$.  Let
$D_\omega=\{(g_\alpha(z,\theta),U_\alpha)\}$ be the corresponding
Cartier
divisor. Then we have an isomorphism
$$
\Ann(\mathcal L\otimes \mathcal M,\omega)\to \mathcal L(D_\omega),
\quad
\omega_\alpha(z,\theta)^\prime l_\alpha\otimes m_\alpha\mapsto
l_\alpha\otimes
1/g_\alpha(z,\theta)
$$
The sequence \eqref{eq:cupsequence} is equivalent to
$$
0\to \mathcal L\to \mathcal L(D_\omega)\to \mathcal
L(D_\omega)|_{D_\omega}\to
0.
$$
Now the cup product with $\omega$ gives a map $H^1(\mathcal
L_\omega)\to H^1(X,
\Ann(\mathcal L\otimes \mathcal M,\omega))$ with kernel the image of
the
natural map $\phi:H^0(X,Q)\to H^1(\mathcal L_\omega)$. Identifying
$H^0(X,Q)$
with $H^0(X, \mathcal L(D_\omega)|_{D_\omega})$, we see that $\phi$
is the
composition
$$
H^0(X,\mathcal L(D_\omega)|_{D_\omega})\to \mathcal Prin (\mathcal
L)\to
H^1(X,\mathcal L).
$$
Therefore the kernel of $\omega\cup$ consists of those $\alpha\in
H^1(X,\mathcal L)$ that have a representative $p\in
\mathcal{P}{rin}(\mathcal
L)$ such that $\omega p$ has zero principal part, i.e., the poles in
$p$ are
compensated by the zeros in $\omega$.

Extensions of the form (\ref{eq:extberbystruct})
are classified by  $\delta\in H^1(X,\Berx^*)$: we think of
$\Cox$ as a subsheaf of $\underline {\mathcal Rat}(\Ox\oplus \Berx)$
consisting
on an open set $U$  of pairs $(f, \omega)$ where $\omega\in \Berx(U)$
and $f$ a
rational function such that the principal part $\bar f$ is equal to
$\omega p$,
for $p\in \mathcal Prin(\Berx^*)$; then $\delta\in H^1(X,\Berx^*)$ is
the class
of $p$. It is then easy to see that the connecting map $q:
H^0(X,\Berx)\to
H^1(X,\Ox)$   is cup product by the extension class $\delta$:
$q(\omega)=\omega\cup \delta$. The class $\delta$ is  ly represented
by the
{\it \v Cech} cocycle
\begin{equation}\label{eq:cechcocycle}
\phi_{\beta\alpha} = \partial_\theta F_{\beta\alpha} /
\partial_\theta \Psi_{\beta\alpha}\in \Berx^*(U_\alpha\cap U_\beta),
\end{equation}
 from \eqref{eq:coordchangen=2}, \eqref{eq:transformrho2}.

 \begin{lem}\label{lem:q=0iffexttriv}
 Let $(X,\Ox)$ be a generic SKP curve. Then the connecting
homomorphism $q:
H^0(X,\Berx)\to H^1(X,\Ox)$ in \eqref{eq:cohomn=2} is the zero map
iff the
extension \eqref{eq:extberbystruct} is trivial. In particular the
cohomology of
$\Cox$ is free iff the extension is trivial.
 \end{lem}

\begin{proof}It is clear that if the extension is trivial the
connecting map
$q$ is trivial.

 From the explicit form, in particular the $\theta$ independence, of
the
cocycle we see that it is not immediately killed by multiplication by
an odd
free section $\omega$ of $\Berx$, i.e., the cocycle is not zero in
the
cohomology group $H^1(X,\Berx^*)/\Ann(\Berx^*,\omega)$ if it is
nonzero in
$H^1(X, \Berx^*)$. The split sheaf $\Berxsplit$ is $\mathcal
K\mathcal
N\inv\mid \mathcal K$. An odd free section $\omega$ of $\Berx$
therefore has an
associated divisor $D_\omega$ as constructed above with reduced
support
included in the divisor of  a section of $\mathcal K$ on the
underlying curve.
Now $q(\omega)=\omega\cup \delta$ is zero if the zeros of $\omega$
cancel the
poles occuring in the principal part $p$ representing $\delta$. But
by
classical results the complete linear system of $\mathcal K$ has no
base
points, i.e.,
there is no point on $X$ where where all global sections of $\mathcal
K$
vanish. This means that wherever the poles of $p$ occur, there will
be a
section $\omega$ of $\Berx$ that does not vanish there. So $q$ being
zero on
all odd generators of $H^0(X,\Berx)$ implies that the extension is
trivial. A
fortiori if $q$ is the zero map the extension will also be trivial.
\end{proof}

The extension given by the cocycle \eqref{eq:cechcocycle} is trivial
if
\begin{equation} \label{eq:phitriv}
\phi_{\beta\alpha}(z_\alpha) = \sigma_\beta(z_\beta,\theta_\beta) -
H_{\beta\alpha}(z_\alpha,\theta_\alpha)
\sigma_\alpha(z_\alpha,\theta_\alpha)
\end{equation}
for some 1-cochain $\sigma_\alpha \in\Berx^*(U_\alpha)$.
In that case, a splitting $e: \Berx \rightarrow \Cox$ is obtained by
$e(f_\alpha) = \rho_\alpha - \sigma_\alpha(z_\alpha,\theta_\alpha)$.

\begin{thm} \label{thm:split=proj}
For a generic SKP curve $(X,\Ox)$, $\Cox$ is a trivial extension of
$\Berx$ iff $(X,\Ox)$ is projected.
\end{thm}

\begin{proof}
We have already observed (in subsection \ref{ss:DualN=2curves}) that
$X$
projected implies
$\phi_{\beta\alpha} = 0$ in a projected atlas, making the extension
trivial.

Now suppose, if possible, that the extension is trivial but that $X$
is not projected.
Write the transition functions of $X$ in the form
\begin{equation*}
z_\beta = f_{\beta\alpha}(z_\alpha) + \theta_\alpha
\eta_{\beta\alpha}(z_\alpha), \;\;\;\;
\theta_\beta = \psi_{\beta\alpha}(z_\alpha) + \theta_\alpha
B_{\beta\alpha}(z_\alpha)
\end{equation*}
and assume that the atlas has been chosen so that
$\eta_{\beta\alpha}$ vanishes to the highest possible (odd) order
$n$ in nilpotents.
That is, $\eta_{\beta\alpha} = 0$ mod $\mathfrak m^n$, but not mod
$\mathfrak m^{n+2}$.
Writing also $\sigma_\alpha(z_\alpha,\theta_\alpha) =
\chi_\alpha(z_\alpha) + \theta_\alpha h_\alpha(z_\alpha)$ and
substituting in (\ref{eq:phitriv}) yields two conditions.
 From the $\theta_\alpha$-independence of $\phi_{\beta\alpha}$ one
finds that $h_\alpha$ mod $\mathfrak m^{n+1}$ is a global section
of ${\mathcal K}^{-1}{\mathcal N}^2$.
Since $X$ is a generic SKP curve, there are no such sections and
$h_\alpha=0$ to this order.
Using this, the second condition becomes,
\begin{equation*}
\eta_{\beta\alpha} = B_{\beta\alpha} \chi_\beta(f_{\beta\alpha}) -
f'_{\beta\alpha} \chi_\alpha \;\;\;
{\text{mod }} {\mathfrak m}^{n+2}.
\end{equation*}
This condition implies that the coordinate change
$\tilde{z}_\alpha = z_\alpha - \theta_\alpha \chi_\alpha$
will make $\eta_{\beta\alpha}$ vanish to higher order than $n$, a
contradiction.
\end{proof}

To lowest order in nilpotents, the cocycle conditions for the
transition functions of $X$ imply that
$\eta_{\beta\alpha}/B_{\beta\alpha}$ is a cocycle for
$H^1(X,\mathcal{NK}^{-1})$, while $\psi_{\beta\alpha}$ is a
cocycle for $H^1(X,\mathcal N^{-1})$.
This implies that the projected $X$'s have codimension
$(0 \mid 3g-3)$ in the moduli space of generic SKP curves, which
has dimension $(4g-3 \mid 4g-4)$ (see
\cite{Vain:DeformSupSpacShe}).
The proof of Theorem \ref{thm:split=proj} generalizes to higher order
in
nilpotents the fact that at lowest order $\phi_{\beta\alpha}$ is a
cocycle in
$H^1(X,\mathcal{NK}^{-1} \mid \mathcal N^2 \mathcal K^{-1})$.

\subsection{$\Cox$ as extension of $\Berxhat$ and symmetric period
matrices}
\label{ss:Symmperiodmatrices}
One can equally view $\Cox$ as an extension of $\Berxhat$ by
$\Oxhat$.
Obviously, if $(X,\Oxhat)$ is projected this extension is trivial,
but the converse no longer holds.
In the proof of  Theorem \ref{thm:split=proj} there is now the
possibility that
$h_\alpha \neq 0$.
(Recall from subsection \ref {ss:SRS} that for $X$ a SRS, a splitting
of the
extension was
universally given by $\chi_\alpha=0,h_\alpha=-1$.)
One can see that this extension is not always trivial, however, by
constructing examples with $\psi_{\beta\alpha} = 0$ and
$\phi_{\beta\alpha}$ a nontrivial class.
(We are now refering to an atlas for $(X,\Oxhat)$.)
In this subsection we will exhibit a connection between the structure
of $\Cox$
as extension of $\Berxhat$ and the symmetry
of the component $Z_e$ of the period matrix, see
\eqref{eq:normperiodmatrix}.

By classical results $Z_e^{\text{red}}$ is symmetric.
However, there seems to be no reason that $Z_e$ is  symmetric in
general.

\begin{thm} \label{thm:Zsym&noZo=proj}
Let $(X,\Ox)$ be a generic SKP curve and $Z_e,Z_o$ its (partially)
normalized
period matrices (as in \eqref{eq:normperiodmatrix}).
Then we have $Z_e$ symmetric and $Z_o=0$
iff $(X,\Ox)$ is projected.
\end{thm}

\begin{proof}
This follows immediately from Theorem \ref{thm:split=proj}, Lemma
\ref{lem:q=0iffexttriv} and the explicit form
\eqref{eq:matrixconnectinghom} of
the connecting homomorphism $q$.
\end{proof}

Recall the exact sequence
$$
0\to \Lambda\to \Cox \overset{(\Dc,\Dchat )}\to \mathcal M
\rightarrow 0,
$$
where $\mathcal M=\Berx\oplus\Berxhat$ is the sheaf of objects that
can be
integrated on $\Cox$.
The corresponding cohomology sequence is in part
\begin{equation*}
0 \to \Lambda \to H^0(X,\Cox)
\stackrel{(\Dc,\Dchat)}{\longrightarrow}
H^0(X,{\mathcal M}) \stackrel{\operatorname{cper}}{\longrightarrow}
H^1(X,\Lambda)
\end{equation*}
where $\operatorname{cper}(\omega,\hat\omega)=\{\sigma\mapsto
\int_\sigma
[\omega+\hat\omega]\}$. So we see that we can identify
$H^0(X,\Cox)/\Lambda$
with
pairs $(\omega,\hat{\omega})$ of differentials with opposite periods.

Now let $(\omega,\hat{\omega})$ be such a pair.
$\omega$ can be written in terms of the basis of $H^0(X,\Berx)$ in
the form
\begin{equation*}
\omega = \sum a_i(\omega) w_i + \sum A_\alpha \omega_\alpha,
\end{equation*}
where $a_i(\omega)$ denote the a-periods and $A_\alpha$ are
other constants uniquely determined by $\omega$.
Then the vector of b-periods of $\omega$ will be
\begin{equation*}
b(\omega) = Z_e a(\omega) + Z_o A.
\end{equation*}
Since these coincide with minus the b-periods of $\hat{\omega}$,
which are
$b(\hat{\omega}) = Z_e^t a(\hat\omega)=-Z_e^t a(\omega)$, we obtain
for each
such pair of
differentials a relation
\begin{equation} \label{eq:basicZrelation}
(Z_e - Z_e^t) a(\omega) + Z_o A = 0.
\end{equation}
We have a sequence analogous to \eqref{eq:cohomn=2} for  $\Cox$ as
extension
of $\Berxhat$ and a connecting map $\hat q$ for this situation.
\begin{thm} \label{thm:Zsym->dualexttriv}
Let $(X,\Ox)$ be a generic SKP curve and $Z_e,Z_o$ its normalized
period matrices.
If $Z_e$ is symmetric, then
$\hat q$ is the zero map.
\end{thm}

\begin{proof}
Assuming that $Z_e=Z_e^t$, we determine the set of pairs
$(\omega,\hat{\omega})$ with opposite periods.
The a-periods of $\hat{\omega}$ can be chosen freely from the kernel
of $Z_o^t$.
According to (\ref{eq:basicZrelation}), any $\omega$ chosen to match
these a-periods will also have matching b-periods iff $A_\alpha$
belongs to the kernel of $Z_o$.
Therefore, $H^0(X,\Cox)$ mod constants can be identified with
$\text{Ker} Z_o \oplus \text{Ker} Z_o^t$, which is precisely
$H^0(X,\Oxhat)/\Lambda \oplus H^0(X,\Berxhat)$.
In this case $\hat q$ is the zero map.
\end{proof}

In general it seems that $\hat q=0$ will not imply that the extension
$\Cox$ of
$\Berxhat$ is trivial, as in Lemma \ref{lem:q=0iffexttriv} for the
extension
of $\Berx$ by $\Ox$. Also it seems  that $Z_e=Z_e^t$ cannot be
deduced from
(\ref{eq:basicZrelation}) as long as
the a-periods are constrained to the kernel of $Z_o^t$.

\subsection{Moduli of invertible sheaves.}\label{ss:modinvertible
sheaves}

In this subsection we will discuss some facts about invertible
sheaves on super
curves and their moduli spaces, see also
\cite{RoSchVor:GeomSupConf,GidNelson:LinebSRS}.

An invertible sheaf on $(X,\Ox)$ is determined by transition
functions
$g_{\alpha\beta}$  on overlaps $U_\alpha\cap U_\beta$, and so
isomorphism
classes of invertible sheaves are classified by the cohomology group
$H^1(X,\Oxevstar)$.

The degree of an invertible sheaf $\mathcal L$ is the degree of the
underlying
reduced sheaf $\mathcal L^{\text{red}}$, with transition functions
$g_{\alpha\beta}^{\text{red}}$. Let $\Picxz$ denote the group of
degree zero
invertible sheaves on $(X,\Ox)$.
The exponential sheaf sequence
\begin{equation}\label{expsequence}
0\to \mathbb Z\to \Oxev\overset{\exp(2\pi i \times \cdot )}\to
\Oxevstar\to 0
\end{equation}
reduces mod nilpotents to the usual exponential sequence for $\Oxred$
and we
see that $\Picxz=H^1(X,\Oxev)/H^1(X,\mathbb Z)$.

If $(X,\Ox)$ is a generic SKP curve $H^1(X,\Ox)$ is a free rank
$g\mid g-1$
$\Lambda$-module and the map $H^1(X,\mathbb Z)\to H^1(X,\Ox)$ is the
restriction of the map $H^1(X,\Lambda)\to  H^1(X,\Ox)$, which is dual
to
the map $\per$ of Lemma \ref{lem:perrep}. So with respect to a
suitable basis
$H^1(X,\mathbb Z)\to  H^1(X,\Ox)$ is described by the transpose of
the period
matrix \eqref{eq:normperiodmatrix}. This implies that the image of
$H^1(X,\mathbb Z)$ is generated by $2g$ elements that are linearly
independent
over the real part $\Lambda_\Re$ of $\Lambda$ (see Appendix
\ref{ss:realstrconj} for the definition of $\Lambda_\Re$). The
elements of the
quotient
$\Picxz=H^1(X,\Oxev)/H^1(X,\mathbb Z)$ are the
$\Lambda$-points of a super torus of dimension $(g\mid g-1)$. Each
component of
$\Pic(X)$ is then isomorphic as a supermanifold to this
supertorus.

 In general, however, $H^1(X,\Ox)$ is not free, nor is the image of
$H^1(X,\mathbb Z)$ generated by $2g$ independent vectors. It seems an
interesting question to understand $\Picxz$ in this generality.

For any supercurve $(X,\Ox)$ we define the Jacobian by
$$\Jacx=H^0(X,\Berx)^*_{\text{odd}}/H_1(X,\mathbb Z),$$
 where elements of
$H_1(X,\mathbb Z)$ act by odd linear functionals on holomorphic
differentials
from $H^0(X,\Berx)$ by integration over 1-cycles.

We have, as discussed in Appendix \ref{app:dualSerredual}, a pairing
of
$\Lambda$-modules
\begin{equation}\label{eq:lambdapairing}
H^1(X,\Ox)\times H^0(X,\Berx)\to \Lambda.
\end{equation}
As we will discuss in more detail in subsection
\ref{ss:effdivisorpoinc}
invertible sheaves are also described by divisor classes. We use this
in the
following Theorem.

\begin{thm}
The pairing (\ref{eq:lambdapairing}) induces an isomorphism of the
identity
component $\Picxz$ with  the Jacobian $\Jacx$ given by the usual
Abel map: a bundle $\mathcal L \in \Picxz$ with divisor $P-Q$
corresponds to
the class of
linear functionals $\int_Q^P$, modulo the action of $H_1(X,\mathbb
Z)$ by
addition of cycles to
the path from $Q$ to $P$.
\end{thm}
\begin{proof}
Let $\mathcal L \in \Picxz$ have divisor $P-Q$, with the reduced
points
$P^{\text{red}}$ and $Q^{\text{red}}$ contained
in a single chart $U_0$ of a good cover of $X$.
If $P=z-p-\theta\pi$ and $Q=z-q-\theta\xi$,
this bundle has a canonical section equal to unity in every other
chart, and
equal to
$$
\frac{z-p-\theta\pi}{z-q-\theta\xi} = \frac{z-p}{z-q} -
\frac{\theta\pi}{z-q} +
\theta\xi \frac{z-p}{(z-q)^2}
$$
in $U_0$.
In the covering space $H^1(X,\Oxev)$ of $\Picxz$, with covering group
$H^1(X,\mathbb Z)$, $\mathcal L$ lifts to a discrete
set of cocycles given by the logarithms of the transition functions
of
$\mathcal L$ in the chart overlaps, namely
$$
a_{0i} = \frac{1}{2\pi i}[\log (z-p) - \log (z-q) -
\frac{\theta\pi}{z-p} +
\frac{\theta\xi}{z-q}]
$$
in $U_0 \cap U_i$, and zero in other overlaps.
The covering group acts by changing the choice of branches for the
logarithms.
We now fix the particular choice for which the branch cut $C$ from
$Q$ to $P$
lies entirely in $U_0$ and meets no other
$U_i$.
Under the Dolbeault isomorphism, this cocycle corresponds to a
$(0,1)$ form
most conveniently represented by the current
$\bar{\partial}a_i$ in $U_i$, where $a_{ij} = a_i - a_j$ and
$\bar{\partial} = d\bar{z} \partial_{\bar{z}} + d\bar{\theta}
\partial_{\bar{\theta}}$.
It is supported on the branch cut $C$, and we can take $a_i = 0$ for
$i \ne 0$.
The pairing (\ref{eq:lambdapairing}) now associates to this the
linear
functional on $H^0(X,\Berx)$
which sends $\omega \in H^0(X,\Berx)$, written as $f(z) + \theta
\phi(z)$ in
$U_0$, to \cite{HaskeWells:Serreduality}
$$
\int_X i(\partial_{\bar{z}}) \bar{\partial}a_0 \, \omega \bar{\theta}
\,[dz\,
d\bar{z}\, d\bar{\theta} \, d\theta] =
\int_X (\partial_{\bar{z}}a_0) \omega dz \,d\bar{z}\, d\theta.
$$
By the definition of the derivative of a current
\cite{GrHa:PrincAlgGeom} and
Stokes' theorem this can be rewritten
\begin{multline*}
- \int_{\partial(X-C)} dz \int d\theta \, a_0 (f+\theta\phi)= \\
= - \frac{1}{2\pi i} \int_{\partial(X-C)} dz \{[\log(z-p) -
\log(z-q)]\phi +
[\frac{\xi}{z-q} - \frac{\pi}{z-p}]f\},
\end{multline*}
where $\partial(X-C)$ denotes the limit of a small contour enclosing
$C$.
Using the residue theorem and the discontinuity of the logarithms
across the
cut, this evaluates to
$$ \int_C \phi \, dz + \pi f(p) - \xi f(q) = \int_Q^P \omega .$$
By linearity of the pairing (\ref{eq:lambdapairing}),
we can extend this correspondence to arbitrary bundles of degree zero
by taking
sums of divisors of the form $P_i-Q_i$.
In particular, the divisor $(P-Q) + (P_1-P) + (P_2-P_1) + \cdots +
(P_n-P_{n-1}) + (P-P_n)$ is equivalent to $P-Q$,
but if the contour $PP_1P_2 \cdots P_nP$ represents a nontrivial
homology class
then the corresponding linear
functionals $\int_Q^P$ differ by addition of this cycle to the
integration
contour.
This shows that the action of $H_1(X,\mathbb Z)$ specified in the
definition of
$\Jacx$ is the correct one.
\end{proof}

\subsection{Effective divisors and Poincar\'e sheaf for generic SKP
curves.}
\label{ss:effdivisorpoinc}
Another description of invertible sheaves is given by divisor
classes.

Recall that a divisor $D\in \Divx$ is a global section of
the sheaf $\Ratxevstar/\Oxevstar$, so $D$, up to equivalence, is
given by a
collection
$(f_\alpha,U_\alpha)$ where the $f_\alpha$ are even
invertible rational functions that are on overlaps related by
an element of $\Oxevstar(U_\alpha\cap U_\beta)$. Each
$f_\alpha$ reduces mod nilpotents to a nonzero rational
function $f_\alpha^{\text{red}}$ on the reduced curve, so
that $D$ determines a divisor $D^{\text{red}}$. Then the {\it
degree} of $D$ is the usual degree of its reduction
$D^{\text{red}}$. We have a mapping $\Ratxevstar\to \Divx$,
$f\mapsto (f)$, and elements $(f)$ of the image are called
{\it principal}. Two even invertible rational functions $f_1,
f_2$ give rise to the same divisor iff $f_1=kf_2$ where $k\in
H^0(X,\Oxevstar)$. So if $(X,\Ox)$ is a generic
SKP curve $k$ is just an even invertible element of $\Lambda$
but in general more exotic possibilities for $k$ exist. A
divisor $D$ is {\it effective}, notation $D\ge 0$, if all
$f_\alpha\in \Oxev(U_\alpha)$. An invertible $\Ox$-module
$\mathcal L$ can be thought of as a submodule of rank $1|
0$ of the constant sheaf $\Ratx$. If $\mathcal
L(U_\alpha)=\Ox(U_\alpha)e_\alpha$, then
$e_\alpha\in \Ratxevstar$ and $\mathcal L$ determines the
divisor $D=\{(f_\alpha=e_\alpha\inv,U_\alpha)\}$. Conversely
any divisor $D$ determines an invertible sheaf $\Ox(D)$ (in
$\Ratx$) with local generators $e_\alpha=f_\alpha\inv$. Two
divisors $D_1=\{(f^{(1)}_\alpha,U_\alpha)\}$ and
$D_2=\{(f^{(2)}_\alpha,U_\alpha)\}$ give rise to  equivalent
invertible sheaves iff they are {\it linearly equivalent},
i.e., $D_1=D_2+(f)$ for some element $f$ of $\Ratxevstar$, or
more explicitly iff $f^{(1)}_\alpha=ff^{(2)}_\alpha$ for all
$\alpha$. If $f\in \Ratxevstar$ is a global section of an
invertible sheaf $\mathcal L=\Ox(D)$ then $D+(f)\ge 0$ and
vice versa. The {\it complete linear system} $|D|=|\Ox(D)|$
of a divisor (or of the corresponding invertible sheaf) is
the set of all effective divisors linearly equivalent to $D$. So
we see that if $\mathcal L=\Ox(D)$ then
$$
|D|\simeq H^0(X,\mathcal
L)_{\text{ev}}^\times/H^0(X,\Ox)_{\text{ev}}^\times.
$$
In case the cohomology of $\mathcal L$ is free of rank
$p+1\mid q$ and $H^0(X,\Ox)$ is just the constants
$\Lambda\mid 0$, the complete linear system $|D|$ is (the set
of $\Lambda$-points of) a super projective space $\mathbb
P^{p\mid q}_\Lambda$. In particular, if $(X,\Ox)$ is a
generic SKP curve and the degree $d$  of $\mathcal L$ is $\ge
2g-1$ the first cohomology of $\mathcal L$ vanishes, the
zeroth cohomology is free of rank $d +1-g\mid d+1-g$ and
$|D|\simeq \mathbb P_\Lambda^{d-g\mid d+1-g}$.

Let $\hat X=(X,\Oxhat)$ be the dual curve and denote by
$\hatxd$ the $d$-fold symmetric product of $\hat X$, see
\cite{DomPerHerRuiSanchSal:Superdiv}. This smooth
supermanifold of dimension $(d\mid d)$ parametrizes
effective divisors of degree $d$ on $(X,\Ox)$. We have a map
(called {\it Abelian sum}) $A:\hatxd\to \Pic^d(X)$ sending an
effective divisor $D$ to the corresponding invertible sheaf
$\Ox(D)$. An invertible sheaf $\mathcal L$ is in the image of
$A$ iff $\mathcal L$ has a even invertible global section: if
$D\in \hatxd$ and $\mathcal L=A(D)$ then the fiber of $A$ at
$\mathcal L$ is the complete linear system $|D|$. If the
degree $d$ of $\mathcal L$ is at least $2g-1$ $H^1(X,\mathcal L)$ is
zero and
hence the cohomology
of $\mathcal L$ is free. So in that case  $A$ is
surjective and  the fibers of $A$ are all projective spaces $\mathbb
P^{d-g\mid d+1-g}$ and $A$ is in fact a fibration.

The symmetric product $\hatxd$ is a universal parameter space
for effective divisors of degree $d$. This  is studied in
detail by Dom\'\i nguez P\'erez et al.
\cite{DomPerHerRuiSanchSal:Superdiv}; we will summarize some
of their results and refer to their paper for more details.
(In fact they consider curves over a field, but the theory is
not significantly different for curves over $\Lambda$.) A
{\it family of effective divisors } of degree $d$ on
$X$ parametrized by a super scheme $S$ is a pair $(S,D_S)$,
where $D_S$ is a Cartier divisor on $X\times_\Lambda S$ such
that for any morphism $\phi:T\to S$ the induced map
$(1\times\phi)^*\mathcal O_{X\times S}(-D_S)\to
(1\times\phi)^*\mathcal O_{X\times S}$ is injective and such
that for any $s\in S$ the restriction  of $D_S$ to $X\times
\{s\}\simeq X$ is an effective divisor of degree $d$. For example, in
$X\times \hatxd$ there is a canonical divisor $\Delta^{(d)}$
such if $p_D$ is any $\Lambda$-point of $\hatxd$
corresponding to a divisor $D$ then the restriction of
$\Delta^{(d)}$ to $X\times \{p_D\}\simeq X$ is just $D$. Then
$(\hatxd, \Delta^{(d)}) $ is universal in the sense that for
any family $(S,D_S)$ there is a unique morphism $\Psi:S\to
\hatxd$ such that $D_S=\Psi^*\Delta^{(d)}$.

A {\it family of invertible sheaves  } of degree $d$ on $X$
parametrized by a super scheme $S$ is a pair $(S,\mathcal
L_S)$, where $\mathcal L_S$ is an invertible sheaf on
$X\times_\Lambda S$ such that  for any $s\in S$ the
restriction  of $\mathcal L_S$ to $X\times \{s\}$ is a sheaf
of degree $d$ on $X$. For example, $(\hatxd,\mathcal
O_{X\times \hatxd}(\Delta^{(d)})$ is a family of invertible
sheaves of degree $d$. Two families $(S,\mathcal L_1)$,
$(S,\mathcal L_2)$ are equivalent if $\mathcal L_1=\mathcal
L_2\otimes \pi_S^* \mathcal N$, where $\pi_S:X\times S\to S$
is the canonical projection and $\mathcal N$ is an invertible
sheaf on $S$. For example, fix a point $x$ of $X$;
then$(\hatxd,\mathcal O_{X\times \hatxd}(\Delta^{(d)}))$ is
equivalent to $(\hatxd,\mathcal R_{x}))$, where $\mathcal
R_x=\mathcal O_{X\times \hatxd}(\Delta^{(d)})\otimes
\pi_{\hatxd}^*[\mathcal O_{X\times
\hatxd}(\Delta^{(d)})|_{\{x\}\times\hatxd})]\inv$. The family
$(\hatxd,\mathcal R_x)$ is normalized: it has the property
that $\mathcal R_x$ restricted to $\{x\}\times \hatxd$ is
canonically trivial. Now consider the mapping $(1\times
A):X\times \hatxd\to X\times\Pic^d(X)$ and the direct image
$\mathcal P^{(d)}_x=(1\times A)_* \mathcal R_x$.

\begin{thm}
Let $(X,\Ox)$ be a generic SKP curve. Let $d\ge
2g-1$. Then $\mathcal P^{(d)}_x$ is a Poincar\'e sheaf on
$X\times\Pic^d(X)$, i.e., $(\Pic^d(X), \mathcal P^{(d)}_x)$
is a family of invertible sheaves of degree $d$ that is
universal in the sense that for any family $(S,\mathcal L)$ of
degree $d$ invertible sheaves there is a unique morphism
$\phi:S\to \Pic^d(X)$ so that $\mathcal L=\phi^*\mathcal
P^{(d)}_x$. Furthermore $\mathcal P^{(d)}_x$ is normalized so
that the restriction to $\{x\}\times \Pic^d(X)$ is
canonically trivial.
\end{thm}

\subsection{Berezinian bundles.}\label{subs:Berbundles}
We continue with the study of a generic SKP  curve
$(X,\Ox)$; we fix an integer $n$ and write $\mathcal P$ for
$\mathcal P^n_X$, the Poincar\'e sheaf on $X\times
\Pic^n(X)$. Let $\mathcal L_s$ be an  invertible sheaf corresponding
to $s\in \Pic^n(X)$. The cohomology groups $H^i(X,\mathcal
L_s)$  will vary as $s$ varies over $\Pic^n(X)$ and can in
general be nonfree, as we have seen. Even if the cohomology
groups are free $\Lambda$-modules their ranks will jump.
Still it is possible to define  an invertible sheaf $\Ber$ over
$\Pic^n(X)$ with fiber at $s$ the line
$$
\ber(H^0(X,\mathcal L_s))\otimes \ber^*(H^1(X,\mathcal L_s)),
$$
in case $\mathcal L_s$ has free cohomology. Here $\ber(M)$ for
a free rank $d\mid \delta$ $\Lambda$-module with basis
$\{f_1,\dots,f_d,\phi_1,\dots,\phi_\delta\}$ is the rank 1
$\Lambda$-module with generator
$B[f_1,\dots,f_d,\phi_1,\dots,\phi_\delta]$. If we are given
another basis
$\{f^\prime_1,\dots,f^\prime_d,\phi^\prime_1,\dots,\phi^\prime
_\delta\}=g\cdot \{f_1,\dots,f_d,\phi_1,\dots,\phi_\delta\}$,
with $g\in Gl(d\mid\delta,\Lambda)$, we have the relation
$$
B[f^\prime_1,\dots,f^\prime_d,\phi^\prime_1,\dots,\phi^\prime_\delta]
=\ber(g)B[f_1,\dots,f_d,\phi_1,\dots,\phi_\delta].
$$
Similarly  $\berdual(M)$ is defined using the inverse
homomorphism $\berdual$. Here $\ber$ and $\berdual$ are the
group homomorphisms defined in (\ref{eq:defberber*}).

The invertible sheaf $\mathcal L_s$ is obtained from the
Poincar\'e sheaf via $i_s^*\mathcal P$. We can reformulate
this somewhat differently: $\mathcal P$ is an $\mathcal
O_{\Pic^n(X)}$-module and for every $\Lambda$-point $s$ of
$\Pic^n(X)$, via the homomorphism $s^\sharp:\mathcal
O_{\Pic^n(X)}\to \Lambda$, also $\Lambda$ becomes an
$\mathcal O_{\Pic^n(X)}$-module, denoted by $\Lambda_s$. Then
$\mathcal L_s=i_s^*\mathcal P=\mathcal P\otimes_{\mathcal
O_{\Pic^n(X)}} \Lambda_s$. It was Grothendieck's idea to
study the cohomology of $\mathcal P\otimes M$ for arbitrary
$\mathcal O_{\Pic^n(X)}$-modules $M$. We refer to Kempf
(\cite{Kempf:AbInt}) for an excellent discussion and more
details on these matters.

The basic fact is that, given the Poincar\'e bundle $\mathcal
P$ on  $X\times \Pic^n(X)$, there is a homomorphism
$\alpha:\mathcal F\to \mathcal G$ of locally free coherent
sheaves on $\Pic^n(X)$ such that we get for any sheaf of
$\mathcal O_{\Pic^n(X)}$-modules $M$ an exact sequence
\begin{gather*}
\begin{aligned}
 0 \to H^0(X\times \Pic^n(X),\mathcal P\otimes M)\to \mathcal
F\otimes
M\overset
{\alpha\times
1_M}\to \mathcal G\otimes M\to \\
\to H^1(X\times \Pic^n(X),
\mathcal P\otimes M)\to 0.
\end{aligned}
\end{gather*}
The proof of this is the same as for the analogous statement in the
classical case, see \cite{Kempf:AbInt}.

Now $\mathcal F$ and $\mathcal G$ are locally free, so for small
enough open
sets $U$ on $\Pic^n(X)$ one can define
$\ber(\mathcal F(U))$ and $\ber^*(\mathcal G(U))$. This
globalizes to invertible sheaves $\ber(\mathcal F)$ and
$\ber^*(\mathcal G)$. Next we  form the ``Berezinian
of the cohomology of $\mathcal P$'' by defining
$\Ber=\ber(\mathcal F)\otimes \ber^*(\mathcal G)$. Finally one
proves, as in  Soul\'e, \cite{Soul:Arakelov}, VI.2, Lemma 1,
that $\Ber$ does not depend, up to isomorphism, on the choice
of homomorphism $\alpha:\mathcal F\to \mathcal G$.

\begin{thm}\label{thm:1ChernBertriv} The first Chern class of
the $\Ber$ bundle is zero. \end{thm}
We will prove this theorem
in subsection \ref{ss:ChernclassBeronPic}, after the
introduction of the infinite super Grassmannian and the
Krichever map. The topological triviality of the $\Ber$ bundle
is a fundamental difference from the situation of classical
curves: there the determinant bundle on $\Pic$ is ample.

Next we consider the special case of $n=g-1$. In this case,
because of Riemann-Roch (\ref{superRR}) $\mathcal F$ and
$\mathcal G$ have the same rank. Indeed, locally
$\alpha(U):\mathcal F(U)\to \mathcal G(U)$ is given, after
choosing bases, by a matrix over $\mathcal O_{\Pic^n(X)}(U)$ of
size $d\mid\delta \times e\mid \epsilon$, say. If we fix a
$\Lambda$-point $s$ in $U$ we get a homomorphism
$\alpha(U)_s:\mathcal F(U)\otimes \Lambda_s\to \mathcal
G(U)\otimes\Lambda_s$ represented by a matrix over $\Lambda$.
The kernel and cokernel are the cohomology groups of
$\mathcal L_s$ and these have the same rank by Riemann-Roch. On the
other hand if the kernel and cokernel of a matrix over
$\Lambda$ are free we have rank of kernel $-$ rank of
cokernel=  $d-e\mid\delta- \epsilon=0\mid 0$. So $\alpha(U)$ is a
square
matrix. This allows us to
define a map
$$
\ber(\alpha):\ber(\mathcal F)\to \ber (\mathcal G).
$$
But this is a (non-holomorphic!) section of $\ber^*(\mathcal
F)\otimes \ber(\mathcal G)$, i.e., of the dual Berezinian
bundle $\mathcal P^*$ on  $\Pic^{g-1}$, because of the non-polynomial
(rational) character of the Berezinian.  This section $\ber(\alpha)$
is
essential for the definition of the $\tau$-function in subsection
\ref{ss:Bakerf-fullsuperH-tau}.

\subsection{Bundles on the Jacobian; theta functions}

We continue with $X$ being a generic SKP curve.
Super theta functions will be defined as holomorphic sections
of certain ample bundles on $J=\text{Jac}(X)$, when such bundles
exist. (As usual, the existence of ample invertible sheaves is
necessary and
sufficient for projective embeddability.)
Given one such bundle, all others with the same Chern
class $c_1$ are obtained by tensor product with bundles having
trivial Chern class, so we begin by determining these, that
is, computing $\text{Pic}^0(J)$.
As we briefly discussed in subsection \ref{ss:modinvertible sheaves}
$J$ is the
quotient of the affine super space
$V={\mathbb A}^{g|g-1} = \Spec
\Lambda[z_1,\ldots,z_g,\eta_1,\ldots,\eta_{g-1}]$
by the lattice $L$ generated by the columns of the transposed
period matrix:
\begin{equation}\label{eq:latticegen}
\begin{aligned}
\lambda_i: & \quad z_j \rightarrow z_j + \delta_{ij}, &\quad
		&\eta_\alpha \rightarrow\eta_\alpha,\\
\lambda_{i+g}: &\quad  z_j \rightarrow z_j+(Z_e)_{ij}, &\quad
		&\eta_\alpha \rightarrow  \eta_\alpha + (Z_o)_{i
\alpha}, \quad
i=1,2,\ldots,g.
\end{aligned}
\end{equation}

We will often omit the parity labels $e,o$ on $Z$, since the index
structure makes clear which is meant.

Any line bundle $\mathcal L$ on such a supertorus $J$ lifts to a
trivial bundle on the covering space $V$.
A section of $\mathcal L$
lifts to a function on which the translations $\lambda_i$ act
by multiplication by certain invertible holomorphic functions,
the {\it multipliers} of $\mathcal L$. We can factor the quotient map
$V \rightarrow J$ through the cylinder $V/L_0$, where $L_0$ is
the subgroup of $L$ generated by the first $g$ $\lambda_i$
only. Since holomorphic line bundles on a cylinder are
trivial, this means that the multipliers for $L_0$ can always
be taken as unity. We have $\text{Pic}^0(J) \cong
H^1(J,\mathcal{O}_{\text{ev}})/H^1(J,\mathbb{Z})$. It is very
convenient to compute the numerator as the group cohomology
$H^1(L,\mathcal{O}_{\text{ev}})$ of $L$ acting on the even
functions on the covering space $V$, in part because the
cocycles for this complex are precisely (the logarithms of)
the multipliers. For the basics of group cohomology, see for
example \cite{Si:ArithEllipticCurves,Mum:AbelVar}. In
particular, factoring out the subgroup $L_0$ reduces our
problem to computing $H^1(L/L_0,\mathcal{O}^{L_0})$, the
cohomology of the quotient group acting on the $L_0$-invariant
functions.

A 1-cochain for this complex assigns to each generator of
$L/L_0$ an even function (log of the multiplier) invariant
under each shift $z_j \rightarrow z_j + 1$,
\begin{equation*}
\lambda_{i+g} \mapsto F^i(z,\eta) = \sum_{\vec{n}}
F_{\vec{n}}^i(\eta) e^{2 \pi i \vec{n} \cdot \vec{z}}.
\end{equation*}
It is a cocycle if the multiplier induced for
every sum $\lambda_{i+g} + \lambda_{j+g}$ is independent of the
order of addition,
which amounts to the symmetry of the matrix $\Delta_i F^j$ giving
the change in $F^j$ under the action of $\lambda_{i+g}$:
\begin{multline*}
F^i(z_k + Z_{jk},\eta_\alpha + Z_{j \alpha}) - F^i(z_k,\eta_\alpha)
=\\
F^j(z_k + Z_{ik},\eta_\alpha + Z_{i \alpha}) - F^j(z_k,\eta_\alpha),
\end{multline*}
or, in terms of Fourier coefficients,
\begin{equation*}
F_{\vec{n}}^i(\eta_\alpha + Z_{j \alpha})
e^{2\pi i\sum_k n_k Z_{jk}} - F_{\vec{n}}^i(\eta) =
F_{\vec{n}}^j(\eta_\alpha + Z_{i \alpha}) e^{2\pi i\sum_k
n_k Z_{ik}} - F_{\vec{n}}^j(\eta).
\end{equation*}
One does
not have to allow for an integer ambiguity in the logarithms
of the multipliers in these equations, precisely because we
are considering bundles with vanishing Chern class. The
coboundaries are of the form,
\begin{equation*}
\lambda_{i+g}
\mapsto A(z,\eta) - A(z_k + Z_{ik},
\eta_\alpha + Z_{i \alpha})
\end{equation*}
for a single function $A$, that is,
those cocycles for which
\begin{equation*}
F_{\vec{n}}^i(\eta)
= A_{\vec{n}}(\eta) - A_{\vec{n}}(\eta_\alpha + Z_{i
\alpha}) e^{2\pi i\sum_k n_k Z_{ik}}.
\end{equation*}
This equation has the form,
\begin{equation*}
F_{\vec{n}}^i(\eta)
= A_{\vec{n}}(\eta)
( 1 - e^{2\pi i\sum_k n_k Z_{ik}}) + O(Z_o).
\end{equation*}

The point now is that, by the linear independence of the columns of
$Z_e^{\text{red}}$, for any $\vec{n} \ne \vec{0}$ there is some
choice of $i$
for which the
reduced part of the exponential in the last
equation differs from unity.
This ensures that, for this $i$, the equation is solvable for
$A_{\vec{n}}$, first to zeroth order in $Z_o$ and then to all orders
by iteration.
Adding this coboundary to the cocycle produces one for which
$F_{\vec{n}}^i = 0$, whereupon the cocycle conditions imply
$F_{\vec{n}}^j = 0$ for all $\vec{n} \ne \vec{0}$ and all $j$ as
well.

Thus the only potentially nontrivial cocycles are independent
of $z_i$.
In the simplest case, when the odd period matrix $Z_o=0$,
all such cocycles are indeed nontrivial, and we have an analog of
the classical fact that bundles of trivial Chern class are specified
by $g$ constant multipliers.
Here a cocycle is specified by giving $g$
even elements $F^i_{\vec{0}}(\eta)$ in the exterior algebra
$\Lambda[\eta_{\alpha}]$ (elements of $H^0(J,\mathcal O_J$)),
leading to $\dim \text{Pic}^0(J) =
g^{2^{g-2}} \mid g^{2^{g-2}}$ (the number of $\eta_{\alpha}$
is $g-1$).
In general, when $Z_o \neq 0$, not all cochains
specified in this way will be cocycles, and some cocycles
will be trivial:
$\text{Pic}^0(J)$ will be smaller, and in general not a
supermanifold.

As to the existence of ample line bundles, let us examine in
the super case the classical arguments
leading to the necessary and sufficient Riemann conditions
\cite{GrHa:PrincAlgGeom,LaBirk:ComplAbVar}.
The Chern class of a very ample bundle is represented in de Rham
cohomology by a $(1,1)$ form obtained as the pullback of the
Chern class of the hyperplane bundle via a projective embedding.
We can introduce real even coordinates $x_i,i=1,\ldots,2g$ for $J$
dual to
the basis $\lambda_i$ of the lattice $L$, meaning that
$x_j \rightarrow x_j + \delta_{ij}$ under the action of
$\lambda_i$.
The associated real odd coordinates $\xi_\alpha,\alpha=1,\ldots,2g-2$
can be taken
to be globally defined because every real supermanifold is split.
The relation between the real and complex coordinates can be taken
to be
\begin{eqnarray*}
z_j & = & x_j + \sum_{i=1}^g Z_{ij} x_{i+g}, \; j=1,\ldots,g, \\
\eta_\alpha & = & \xi_\alpha + i \xi_{\alpha + g-1} +
\sum_{i=1}^g Z_{i \alpha} x_{i+g}, \; \alpha = 1,\ldots,g-1.
\end{eqnarray*}
The de Rham cohomology is isomorphic to that of the reduced torus
and can be represented by translation-invariant forms in the
$dx_i$.
The Chern class represented by a form
$\sum_{i=1}^g \delta_i\, dx_i\, dx_{i+g}$
is called a polarization of type $\Delta =
\text{diag}(\delta_1,\ldots,\delta_g)$ with elementary divisors the
positive integers $\delta_i$.
We consider principal polarizations $\delta_i=1$ only, because
nontrivial nonprincipal polarizations generically do not exist,
even on the reduced torus \cite{Lef:ThmCorrAlgCurv}.
Furthermore, a nonprincipal polarization is always obtained by
pullback
of a principal one from another supertorus whose lattice $L'$
contains $L$ as a sublattice of finite index
\cite{GrHa:PrincAlgGeom}.
Reexpressing the Chern form in complex coordinates, the standard
calculations lead to the usual Riemann condition $Z_e = Z_e^t$ to
obtain a $(1,1)$ form.
Together with the positivity of the imaginary part of the reduced
matrix, the symmetry of $Z_e$ (in some basis) is necessary and
sufficient
for the existence of a $(1,1)$ form
with constant coefficients representing the Chern class.
This can be viewed as the cocycle condition, symmetry of
$\Delta_iF^j$, for the usual multipliers of a theta bundle,
$F^j = -2\pi i z_j$.

The usual argument that the $(1,1)$ form representing the Chern
class can always be taken to have constant coefficients depends
on Hodge theory, particularly the Hodge decomposition of
cohomology, for a K\"ahler manifold such as a torus.
This does not hold in general for a supertorus with $Z_o \neq 0$.
For example, $H_{\text{dR}}^1(J)$ is generated by the $2g$ 1-forms
$dx_i$,
whereas $H^{1,0}(J)$ contains the $g \mid g-1$ nontrivial forms
$dz_i,\,d\eta_\alpha$, with certain nilpotent multiples of the
latter being trivial.
Indeed, since by (\ref{eq:latticegen}) $\eta_\alpha$ is defined
modulo entries of column $\alpha$ of $Z_o$, $\epsilon\eta_\alpha$
is a global function and $\epsilon d \eta_\alpha$ is exact when
$\epsilon \in \Lambda$ annihilates these entries.
Thus, $H^{1,0}(J)$ cannot be a direct summand in
$H^1_{\text{dR}}(J)$.
Correspondingly, some $\eta$-dependent multipliers
$F^j = -2\pi i z_j + \cdots$ may satisfy the cocycle condition
and give ample line bundles.
We do not know a simple necessary condition for a Jacobian to
admit such polarizations.

When $Z_e$ is symmetric, we can construct theta functions
explicitly.
Consider first the trivial case with $Z_o=0$ as well.
Then the standard Riemann theta function $\Theta(z;Z_e)$
gives a super theta function on $\Jacx$, where
$\Theta(z;Z_e)$ is defined by Taylor expansion in the nilpotent part
of
$Z_e$ as usual.
It has of course the usual multipliers,
\begin{equation} \label{eq:thetafactors}
\Theta(z_j + \delta_{ij};Z_e) = \Theta(z_j;Z_e),\;\;\;
\Theta(z_j + Z_{ij};Z_e) = e^{-\pi i (2z_i + Z_{ii})}
\Theta(z_j;Z_e).
\end{equation}
Multiplication of $\Theta(z;Z_e)$ by any monomial in the odd
coordinates
$\eta_{\alpha}$ gives another, even or odd,
theta function having the same multipliers, whereas
translation of the argument $z$ by polynomials in the
$\eta_{\alpha}$ leads to the multipliers for another bundle
with the same
Chern class.

In the general case with $Z_o \neq 0$, theta functions with the
standard multipliers can be constructed as follows.
Such functions must obey
\begin{align*}
H(z_j + \delta_{ij},\eta_{\alpha};Z) &= H(z_j,\eta_{\alpha};Z),\\
H(z_j + Z_{ij},\eta_{\alpha} + Z_{i \alpha};Z) &=
e^{-\pi i (2z_i + Z_{ii})} H(z_j,\eta_{\alpha};Z).
\end{align*}
The function $\Theta(z;Z_e)$ is a trivial example independent of
$\eta$;
to obtain others one checks that
when $H$ satisfies these relations then so does
$$H_\alpha = \left( \eta_{\alpha} + \frac{1}{2\pi i} \sum_k Z_{k
\alpha}
\frac{\partial}{\partial z_k} \right) H.$$
Applying this operator repeatedly one constructs super theta
functions
$\Theta_{\alpha \cdots \gamma}$ reducing to
$\eta_\alpha \cdots \eta_\gamma \Theta(z;Z_e)$
when $Z_o=0$.

``Translated" theta functions which are sections of other bundles
having the same Chern class can be obtained by literally
translating the arguments of these only in the simplest cases.
Constant shifts in the multiplier exponents $F^j$ can be achieved by
constant shifts of the arguments $z_j$.
Shifts linear in the $\eta_\alpha$ are obtained by
$z_j \rightarrow z_j + \eta_\alpha \Gamma_{\alpha j}$,
which is a change in the chosen basis of holomorphic differentials
on $X$, see the discussion after \eqref{eq:normperiodmatrix}.
The resulting theta functions have the new period matrix
$Z_e + Z_o\Gamma$.
More generally, translated theta functions
can be obtained by the usual method of determining their
Fourier coefficients from the recursion relations following from the
desired multipliers.
We do not know an explicit expression for them in terms of
conventional
theta functions.

It is easy to see that any meromorphic function $F$ on the Jacobian
can be rationally expressed in terms of the theta functions we have
defined.
Expand $F(z,\eta) = \sum_{IJ} \beta_I \eta_J F_{IJ}(z)$ in the
generators of
$\Lambda[\eta_{\alpha}]$, with multi-indices $I,J$.
Then the zeroth-order term $F_{00}$ is a meromorphic
function on the reduced Jacobian, hence a rational expression
in ordinary theta functions.
Using $Z_e$ as the period matrix argument of these theta functions
gives a meromorphic function
on the Jacobian itself, whose reduction agrees with $F_{00}$.
Subtract this expression from $F$ to get a meromorphic function on
the Jacobian
whose zeroth-order term vanishes,
and continue inductively, first in $J$, then in $I$.
For example, $F_{0\alpha}$ is equal, to lowest order in the
$\beta$'s, to a rational expression in theta functions of which one
numerator factor is a $\Theta_\alpha$.
Subtracting this expression removes the corresponding term in $F$
while only modifying other terms of higher order in $\beta$'s.

\section{Super Grassmannian, $\tau$-function and Baker
function.}

\subsection{Super Grassmannians.}
In this subsection we will introduce an infinite super Grassmannian
and related
constructions. The infinite Grassmannian of Sato
(\cite{Sa:KPinfDymGr}) or of
Segal-Wilson (\cite
{SeWi:LpGrpKdV}) consists (essentially) of ``half infinite
dimensional'' vector
subspaces $W$ of an infinite dimensional vector space $H$ such that
the
projection on a fixed subspace $H_-$ has finite dimensional kernel
and
cokernel. In the super category we replace this by the super
Grassmannian of
free ``half infinite rank'' $\Lambda$-modules of an infinite rank
free
$\Lambda$-module $H$ such that the kernel and cokernel of the
projection on
$H_-$ are a submodule respectively a quotient module of a free finite
rank
$\Lambda$-module. In \cite{Schwarz:FermStringModSpa} a similar
construction can
be found, but it seems that there $\Lambda=\mathbb C$ is taken as is
also the
case in \cite{Mu:Jac}. This is too restrictive for our purposes
involving
algebraic super curves over nonreduced base ring  $\Lambda$.

Let $\Linfinf$ be the free $\Lambda$-module $\Lambda[z,z\inv,\theta]$
with $z$
an even and $\theta$ an odd variable. Introduce the notation
\begin{equation}\label{eq:basisLinfinf}
e_i=z^i,\quad e_{i-\frac12}=z^i\theta,\quad i\in \mathbb Z.
\end{equation}
We will think of  an element $h=\sum_{i=-N}^\infty h_i e_i$, $h_i\in
\Lambda$
of $\Linfinf$  not only as a series in $z,\theta$ but also as an
infinite
column vector:
$$
h=(\dots,0,\dots,
h_{-1},h_{-\frac12},h_0,h_{\frac12},h_1,\dots,0,\dots)^t
$$
Introduce on $\Linfinf$ an odd Hermitian product
\begin{multline}
\langle f(z,\theta),g(z,\theta)\rangle=
\frac1{2\pi i}\oint
\frac{dz}{z}d\theta\overline{f(z,\theta)}g(z,\theta)=\\
=\frac1{2\pi i}\oint (\overline{f_{\bar 0}}g_{\bar
1}+\overline{f_{\bar
1}}g_{\bar 0})\frac{dz}{z},
\end{multline}
where $\overline{f(z,\theta)}$ is the extension of the complex
conjugation of
$\Lambda$ (see Appendix \ref{ss:realstrconj}) to $\Linfinf$ by
$\overline{z}=z\inv$ and $\overline{\theta}=\theta$, and
$f(z,\theta)=f_{\bar
0}+\theta f_{\bar 1}$, and similarly for $g$. Let $H$ be the
completion of
$\Linfinf$ with respect to the Hermitian inner product.

We have a decomposition $H=H_-\oplus H_+$, where $H_-$ is the closure
of the
subspace spanned by $e_i$ for $i\le0$, and $H_+$ is the closure of
the space
spanned by $e_i$  with $i>0$, for $i\in \halfz$.

The super Grassmannian $\Sgr$ is the collection of all free closed
$\Lambda$-modules $W\subset H$ such that the projection $\pi_-:W\to
H_-$ is
super Fredholm, i.e., the kernel and cokernel are a submodule
respectively a
quotient module of a free finite rank $\Lambda$-module.

\begin{exmpl}  Let $W$ be the closure of the subspace generated by
$\delta +z, \theta$ and $z^i, z^i\theta$ for $i\le -1$, for $\delta$
a
nilpotent even constant. Let $A\subset \Lambda$ be the ideal of
annihilators of
$\delta$. Then $W$ is free and the kernel of $\pi_-$ is
$A (\delta +z) \subset \Lambda (\delta+z)$ and the cokernel is
isomorphic to
$\Lambda/\Lambda\delta$. \qed
\end{exmpl}

Let $I$ be the subset $\{i\in \halfz\mid i\le 0\}$.
We consider matrices with coefficients in $\Lambda$ of size
$\halfz\times I$:
$$
\mathcal W=(W_{ij}) \quad \text{where } i\in \halfz,\  j\in I.
$$
An even matrix of this type is called an {\it admissible frame} for
$W\in \Sgr$
if the closure of the subspace spanned by the columns of $\mathcal W$
is $W$
and if moreover in the decomposition
$\mathcal W=\begin{pmatrix} W_-\\W_+\end{pmatrix}$ induced by
$H=H_-\oplus H_+$
the operator $W_-:H_-\to H_-$ differs from the identity by an
operator of super
trace class and $W_+:H_-\to H_+$ is compact.

Let $Gl(H_-)$  be the group of invertible maps $1+X:H_-\to H_-$ with
$X$ super
trace class. Then the super frame bundle $\Sfr$, the collection of
all pairs
$(W,\mathcal W)$ with $\mathcal W$ an admissible frame for $W\in
\Sgr$, is a
principal $Gl(H_-)$ bundle over the super Grassmannian. Elements of
$Gl(H_-)$
have a well defined berezinian, see
\cite{ Schwarz:FermStringModSpa} for some details.    This allows us
to define
two associated line bundles $\Bersgr$ and $\Bersgrdual$ on $\Sgr$.
More
explicitly, an element of $\Bersgr$ is an equivalence class of
triples
$(W,\mathcal W, \lambda)$, with $\mathcal W$ a frame for $W$,
$\lambda\in
\Lambda$;  here $(W,\mathcal Wg, \lambda)$ and $(W,\mathcal W,
\ber(g)\lambda)$
are equivalent  for  $g\in Gl(H_-)$. For $\Bersgrdual$ replace
$\ber(g)$ by
$\berdual(g)$. For simplicity we shall write $(\mathcal W,\lambda)$
for
$(W,\mathcal W, \lambda)$, as $\mathcal W$ determines $W$ uniquely.

The two bundles $\Bersgr$ and $\Bersgrdual$ each have a canonical
section.
Let $\mathcal W$ be a frame for $W\in \Sgr$ and  write $\mathcal
W=\begin{pmatrix}W_-\\W_+\end{pmatrix}$ as above. Then
\begin{equation}\label{eq:defsigma*}
\sigma(W)=(\mathcal W, \ber(W_-)),\quad \sigma^*(W)=(\mathcal W,
\berdual(W_-)),
\end{equation}
are sections of $\Bersgrdual$ and $\Bersgr$, respectively. It is a
regrettable
fact of life that neither of these sections is holomorphic; indeed
there are no
global sections to $\Bersgr$ or $\Bersgrdual$ at all, see
\cite{Manin:GaugeFieldTheoryComplexGeom}. This is a major difference
between
classical geometry and super geometry.

\subsection{The Chern class of $\Bersgr$ and the
$gl_{\infty\mid\infty}$
cocycle.}

First we summarize some facts about complex supermanifolds that are
entirely
analogous to similar facts for ordinary complex manifolds. Then we
apply this
to the super Grassmannian, following the treatment in
\cite{PrSe:LpGrps} of the
classical case.

Let $M$ be a complex supermanifold. The Chern class of an invertible
sheaf
$\mathcal L$ on $M$  is an element $c_1(\mathcal L)\in H^2(M,\mathbb
Z)$. By
the sheaf inclusion $\mathbb Z\to \Lambda$ and the de Rham theorem
$H^2(M,\Lambda)\simeq H^2_{\text{dR}}(M)$ we can represent
$c_1(\mathcal L)$ by
a closed two form on $M$. On the other hand, if $\nabla:\mathcal L\to
\mathcal
L\otimes \mathcal A^1$, with $\mathcal A^1$ the sheaf of smooth
1-forms, is a
connection compatible with the complex structure, the curvature $F$
of $\nabla$
is also a two form. By the usual proof (see e.g.,
\cite{GrHa:PrincAlgGeom}) we
find that $c_1(\mathcal L) $ and $F$ are equal, up to a factor of
$i/2\pi$.

We can locally calculate the curvature on an invertible sheaf
$\mathcal L$ by
introducing a Hermitian metric $\langle\,, \rangle$ on it: if $s,t\in
\mathcal
L(U)$ then $\langle s, t\rangle(m)$ is a smooth function in $m\in U$
taking
values in $\Lambda$, linear in $t$ and satisfying $\langle s,
t\rangle(m)=\overline{\langle t, s\rangle(m)}$.  Choose a local
generator $e$
of $\mathcal L$ and let $h=\langle e, e\rangle$. The curvature is
then
$F=\bar\partial \partial \log h$, with $\partial=\sum
dz_i\frac{d}{dz_i}+\sum
d\theta_\alpha \frac{\partial}{\partial \theta_\alpha}$ and
$\bar\partial$
defined by a similar formula.

Now consider the invertible sheaf $\Bersgr$ on $\Sgr$. If
$s=(\mathcal W,
\lambda)$ is a section the square length is defined to be $\langle
s,s\rangle=
\bar \lambda \lambda\ber (\mathcal W^H \mathcal W)$, where
superscript ${}^H$
indicates conjugate transpose. Of course, this metric is not defined
everywhere
on $\Sgr$  because of the rational character of $\ber$,
but we are interested in a neighborhood of the point $W_0$ with
standard frame
$\mathcal W_0=\begin{pmatrix} 1_{H_-}\\ 0\end{pmatrix}$ where there
is no
problem. The tangent space at $W_0$ can be identified with the space
of maps
$H_-\to H_+$, or, more concretely, by matrices with the columns
indexed by
$I=\{i\in \halfz\mid i\le 0\}$ and with rows indexed by the
complement of $I$.
Let $x,y$ be two tangent vectors at $W_0$. Then the curvature at
$W_0$ is
calculated to be
\begin{equation}\label{eq:curvature}
F(x,y)=\bar\partial\partial \log h (x,y)= Str(x^Hy-y^Hx),
\end{equation}
where we take as local generator $e=\sigma$, the section defined by
\eqref{eq:defsigma*}, so that $h=\langle \sigma, \sigma\rangle$.
We can map the tangent space at $W_0$ to the Lie super algebra
$gl_{\infty\mid
\infty}(\Lambda)$ via $x\mapsto \begin{pmatrix}
0&-x^H\\x&0\end{pmatrix}$. Here
$gl_{\infty\mid \infty}(\Lambda)$ is the Lie super algebra
corresponding to the
Lie super group $Gl_{\infty\mid \infty}(\Lambda)$ of infinite even
invertible
matrices $g$ (indexed by $\halfz$) with block decompostion
$\begin{pmatrix}a&b\\c&d\end{pmatrix}$ with $b,c$ compact and $a,d$
super
Fredholm. We see that (\ref{eq:curvature}) is the pullback under this
map of
the cocycle on $gl_{\infty\mid \infty}(\Lambda)$ (see also
\cite{KavdL:SuperBoson}) given by
\begin{equation}\label{eq:cocycle}
\begin{aligned}
c:gl_{\infty\mid \infty}(\Lambda)\times gl_{\infty\mid
\infty}(\Lambda)&\to
\quad\quad\quad\Lambda\\
(X,Y)\quad \quad\quad&\mapsto\quad\frac14 \Str(J[J,X][J,Y]),
\end{aligned}
\end{equation}
where $J=\begin{pmatrix} 1_{H_-}&0\\0&-1_{H_+}\end{pmatrix}$. In
terms of the
block decomposition of $X,Y$ we have
$$
c(X,Y)= \Str(c_Xb_Y-b_Xc_Y).
$$

The natural action of $Gl_{\infty\mid \infty}(\Lambda)$ on $\Sgr$
lifts to
a projective action on $\Bersgr$; the cocycle $c$ describes
infinitesimally the
obstruction for this projective action to be a real action. Indeed,
if
$g_1=\exp(f_1),g_2=\exp(f_2)$ and $g_3=g_1g_2$ are all in the open
set of
$Gl_{\infty\mid \infty}(\Lambda)$ where the ${}_{--}$ blocks $a_i$
are
invertible, the action on a point of $\Bersgr$ is given by
\begin{equation} \label{eq:lift}
g_i\circ(\mathcal W,\lambda)=(g_i\mathcal Wa_i\inv,\lambda).
\end{equation}
(One checks as in \cite{SeWi:LpGrpKdV} that if $\mathcal W$ is an
admissible
basis then so is $g\mathcal Wg_{--}\inv$.)
Then we have
$$
g_1\circ g_2 \circ(\mathcal
W,\lambda)=\exp[c(f_1,f_2)]g_3\circ(\mathcal
W,\lambda).
$$
We can also introduce the {\it projective multiplier } $C(g_1,g_2)$
for
elements $g_1$ and $g_2$ that commute in
$Gl_{\infty\mid\infty}(\Lambda)$:
\begin{equation}\label{eq:projmult}
g_1\circ g_2\circ g_1\inv\circ g_2\inv (\mathcal
W,\lambda)=C(g_1,g_2)(\mathcal
W,\lambda),
\end{equation}
where  $C(g_1,g_2)=\exp[S(f_1,f_2)]$ if $g_i=\exp(f_i)$ and
\begin{equation}\label{eq:logprojmult}
S(f_1,f_2)=\Str([f_1,f_2]).
\end{equation}
We will in subsection \ref{ss:ChernclassBeronPic} use the projective
multiplier
to show that the Chern class of the Berezinian bundle on $\Picxz$ is
trivial.

\subsection{The Jacobian super Heisenberg algebra.}
In the theory of the KP hierarchy an important role is played by a
certain
Abelian subalgebra of the infinite matrix algebra and its universal
central
extension, loosely referred to as the (principal) Heisenberg
subalgebra. In
this subsection we introduce one of the possible analogs of this
algebra in
the super case.

Let the {\it Jacobian super Heisenberg algebra}  be the
$\Lambda$-algebra
$\Jheis=\Lambda[z,z\inv,\theta]$. Of course, this is as a
$\Lambda$-module the
same as $\Linfinf$ but now we allow multiplication of elements. When
convenient
we will identify the two; in particular we will use the basis
$\{e_i\}$ of
(\ref{eq:basisLinfinf}) also for $\Jheis$. We think of elements of
$\Jheis$ as
infinite matrices in $gl_{\infty\mid\infty}(\Lambda)$: if $E_{ij}$ is
the
elementary matrix with all entries zero except for the $ij$th entry
which is 1,
then
$$
e_i=\sum_{n\in\mathbb Z} E_{n+i,n}+E_{n+i-\frac12,n-\frac12},\quad
e_{i-\frac12}=
\sum_{n\in\mathbb Z} E_{n+i-\frac12,n}.
$$
We have a decomposition $\Jheis=\Jheis_-\oplus\Jheis_+$ in
subalgebras
$\Jheis_-=z\inv\Lambda[z\inv,\theta]$ and
$\Jheis_+=\Lambda[z,\theta]$.
Elements of $\Jheis_+$ correspond to lower triangular matrices and
elements of
$\Jheis_-$ to upper triangular ones.  By exponentiation we obtain
from
$\Jheis_-$ and $\Jheis_+$ two subgroups $G_-$ and $G_+$ of
$Gl_{\infty\mid\infty}(\Lambda)$, generated by
$$
g_\pm(t)=\exp(\sum_{i\in\pm I} t_i e_i),
$$
where $t_i\in \Lambda$ is homogeneous of the same parity as $e_i$
(and $t_i$ is
zero for almost all $i$, say).

For an element $g=\begin{pmatrix} a&b\\c&d\end{pmatrix}$ of $G_+$
the  block $b$ vanishes, whereas if $g\in G_-$ the block $c=0$. In
either case
the diagonal block $a$ is invertible and we can lift the action of
either $G_-$
or $G_+$ to a (potentially projective) action on $\Bersgr$ and
$\Bersgrdual$,
via (\ref{eq:lift}). Since the cocycle (\ref{eq:cocycle}) is zero
when
restricted to both $\Jheis_-$ and $\Jheis_+$ we get an honest action
of the
Abelian groups $G_\pm$ on $\Bersgr$ and $\Bersgrdual$, just as in the
classical
case.

In contrast with the classical case, however, as was pointed out in
 \cite{Schwarz:FermStringModSpa},  the actions of $G_-$ and $G_+$ on
the line
bundles $\Bersgr$ and $\Bersgrdual$ mutually commute. This follows
from the
following Lemma.

\begin{lem}\label{lem:comactionJheis}
Let $g_\pm\in G_\pm$ and write $a_\pm=\exp(f_\pm)$, with $f_\pm \in
gl(H_-)$.
Then
$$
\Str_{H_-}([f_-,f_+])=0,
$$
so that the actions of $g_-$ and $g_+$ on $\Bersgr$ and $\Bersgrdual$
commute.
\end{lem}
\begin{proof}
The elements $f_\pm$ act on $H_-$ by multiplication by an element of
$\Jheis_\pm$, followed by projection on $H_-$ if necessary. So write
$f_\pm=\pi_{H_-}\circ \sum_{i>0} c^\pm_iz^{\pm i}
+\gamma^\pm_{i}z^{\pm
i}\theta$. To find the supertrace we need to calculate the projection
on
the rank $1\mid 0$ and $0\mid 1$ submodules of $H_-$ generated by
$z^{-i}$ and
$z^{-i}\theta$:
\begin{align*}
f_+ f_- z^{-k}|_{\Lambda z^{-k}}&=
		f_+(\sum_{i>0}c_i^-z^{-i-k})|_{\Lambda
z^{-k}}=(\sum_{i>0}c^+_ic^-_i)z^{-k},\\
f_+ f_- z^{-k}\theta|_{\Lambda z^{-k}\theta}&=
		(\sum_{i>0}c^+_ic^-_i)z^{-k}\theta,\\
f_- f_+ z^{-k}|_{\Lambda z^{-k}}&=
		f_-(\sum_{i=1}^ka_i^+z^{i-k})|_{\Lambda
			z^{-k}}=(\sum_{i=1}^k c^+_ic^-_i)z^{-k},\\
f_- f_+ z^{-k}\theta|_{\Lambda z^{-k}\theta}&=
		f_-(\sum_{i=1}^kc_i^+z^{i-k}\theta)|_{\Lambda
			z^{-k}\theta}=(\sum_{i=1}^k
c^+_ic^-_i)z^{-k}\theta.
\end{align*}
Since the super trace is the difference of the traces of the
restrictions to
the even and odd submodules we see that $\Str([f_+,f_-])=0$ so that,
by (\ref
{eq:projmult},\ref{eq:logprojmult}), the actions of $G_\pm$ on
$\Bersgr$
commute.
\end{proof}

\subsection{Baker functions, the full super Heisenberg algebra, and
$\tau$-functions.}\label{ss:Bakerf-fullsuperH-tau}
We define $W\in \Sgr$ to be {\it in the big cell} if it has an
admissible frame
$\mathcal W^{(0)}$ of the form
\[\mathcal{W}^{\,(0)}=\begin{pmatrix}\ddots
&\vdots&\vdots&\vdots&\vdots\\
				    \dots &1&0&0&0\\
				    \dots &0&1&0&0\\
				    \dots &0&0&1&0\\
				    \dots &0&0&0&1\\
				    ***&*&*&*&*\\
		    \end{pmatrix},
\]
i.e. $(\mathcal W^{\,(0)})_-$ is the identity matrix. Note that the
canonical
sections $\sigma$ and $\sigma^*$ do not vanish, nor blow up, at a
point in the
big cell.

If $\mathcal{W}$ is any frame of a point $W$ in the big cell we
can calculate the
standard frame $\mathcal{W}^{\,(0)}$ through quotients of
Berezinians of minors
of $\mathcal{W}$. Indeed, if we  put $A=\mathcal W_-$ then the
maximal minor $A$ of $\mathcal W$ is invertible and we have
\begin{equation}\label{eq:connw0w}
\mathcal{W}^{\,(0)}A=\mathcal{W}.
\end{equation}
Write $\mathcal{W}^{\,(0)}=\sum  w_{ij}^{(0)}E_{ij}$. Then we
can solve \thetag{\ref{eq:connw0w}} by Cramer's rule,
\thetag{\ref{eq:supercramer}}, to find for
$i>0,j\le 0$:
\begin{equation*}
w_{ij}^{(0)}=\begin{cases}
\ber\,(A_j(r_i))/
		\ber\,(A)&\text{if $j\in \mathbb Z$},\\
\berdual \,(A_j(r_i))/
		\berdual \,(A)&\text{if $j\in \mathbb Z+\frac
12$}.
	\end{cases}
\end{equation*}
Here $A_j(r_i)$ is the matrix obtained from $A$ by replacing
the $j$th row by $r_i$, the $i$th row of $\mathcal W$. In
particular the even and odd ``Baker vectors'' of $W$, i.e.
the zeroth and $-\frac12$th column of $\mathcal W^{\,(0)}$, are
given by
\begin{equation}\label{eq:bakervector}
\begin{split}
w_{\bar 0}&=e_0+\sum_{\frac{i>0}
{i\in \frac 12 \mathbb Z}}\frac{\ber\,(
		A_0(r_i))}{\ber\,(A)}e_i,\\
w_{\bar 1}&=e_{-\frac 12}+\sum_{\frac{i>0}
{ i\in \frac 12 \mathbb Z}}\frac{\berdual \,(
		A_{-\frac 12}(r_i))}{\berdual \,(A)}e_i
\end{split}
\end{equation}
The corresponding ``Baker functions'' are obtained by using
$e_i=z^{i}$, $e_{i-\frac12}=z^{i}\theta$.
Then  \thetag{
\ref{eq:bakervector}} reads
\begin{equation}\label{eq:bakerfunction}
\begin{split}
w_{\bar 0}(z,\theta)&=1 + \sum_{i>0}
z^{i}\frac	{	\ber\,( A_0(r_i))+
			\ber\,( A_0(r_{i-\frac 12}))\theta
		}
      		{\ber\,(A)},\\
w_{\bar 1}(z,\theta)&=\theta + \sum_{i>0}
	z^{i}\frac{	\berdual \,(A_{-\frac12}(r_i))+
			\berdual
\,(A_{-\frac12}(r_{i-\frac12}))\theta}
     {\berdual \,(A)}.
\end{split}
\end{equation}
Here and henceforth (unless otherwise noted) the summations
run over (subsets of) the integers.

The full super Heisenberg algebra $\Sheis$ is the extension
$\Jheis[\frac{d}{d\theta}]=\Lambda[z,z\inv,\theta][\frac{d}{d\theta}]$

. This
is, just as the Jacobian super Heisenberg algebra, a  possible analog
of the
principal Heisenberg of the infinite matrix algebra used in the
standard KP
hierarchy, see \cite{KavdL:SuperBoson}. $\Sheis$ is non--Abelian and
the
restriction of the cocycle \thetag{\ref{eq:cocycle}} to it is
nontrivial, in
contrast to the subalgebra $\Jheis$.

$\Sheis$ acts in the obvious way on $\Linfinf$ and we can represent
it by
infinite matrices from $gl_{\infty\mid\infty}(\Lambda)$. Introduce a
basis for
$\Sheis$ by
\begin{alignat*}{2}
\lambda(n)&=z^{-n}(1-\theta\frac{d}{d\theta})=\sum_{k\in \mathbb
Z}E_{k,k+n},&\quad
f(n)&=z^{-n}\frac{d}{d\theta}=\sum_{k\in \mathbb Z}E_{k,k+n-\frac
12},\\
\mu(n)&=z^{-n}\theta\frac{d}{d\theta}=\sum_{k\in \mathbb
Z}E_{k-\frac12,k-\frac12+n},&\quad
e(n)&=z^{-n}\theta=\sum_{k\in \mathbb Z}E_{k-\frac12,k+n}.
\end{alignat*}
We can rewrite the Baker functions as  quotients of
Berezinians, using $\Sheis$. To this end define the following even
invertible
matrices (over the ring $\mathbb \Lambda[u,\phi,\frac{\partial}
     {\partial\phi}]$):
\begin{align*}
Q_{\bar 0}(u,\phi)&=1+\sum_{n=1}^\infty
{u^n}[\lambda(n)+f(n)\phi],\\
Q_{\bar 1}(u,\phi)&=1+\sum_{n=1}^\infty
{u^n}[{\mu(n)+e(n)\frac{\partial} {\partial\phi}}],
\end{align*}
where $u$, resp $\phi$, is an even, resp. odd, variable. We can let
these
matrices act on $H$ and obtain in this way infinite vectors over the
ring
$\Lambda[u,\phi,\frac{\partial}{\partial\phi}]$. Also we can let
these matrices
act on an admissible frame and obtain a matrix over
$\Lambda[u,\phi,\frac{\partial}{\partial\phi}]$.

\begin{lem} Let $w_{\bar 0}(u,\phi)$ and $w_{\bar 1}(u,\phi)$ be the
even and
odd Baker functions of a point $W$ in the big cell. For any frame
$\mathcal W$
of $W$ we have:
\begin{equation*}
w_{\bar 0}(u,\phi)=\frac{\ber\,([Q_{\bar
0}(u,\phi)\mathcal W]_-)}
	{\ber\,(A)},\quad
w_{\bar 1}(u,\phi)=\frac{\berdual \,([Q_{\bar
1}(u,\phi)\mathcal W]_-)\phi}
	{\berdual \,(A)},
\end{equation*}
with $A=\mathcal W_-$.
\end{lem}
\begin{proof}
Let $r_i$, $r_{i,\bar 0}$ and $r_{i,\bar 1}$, be respectively the
$i$th
row of $\mathcal W$,  $Q_{\bar 0}(u,\phi)\mathcal W$ and of
$Q_{\bar 1}(u,\phi)\mathcal W$. Then one
calculates that for $i\in \mathbb Z$ we have $\quad r_{i-\frac 12,
\bar 0} =
r_{i-\frac 12}$, and $r_{i, \bar 1}= r_{i}$ and :
\begin{equation}\label{eq:tildewinw}
\begin{aligned}
r_{i,\bar 0}
	&=r_i+\sum_{k\ge 1}{u^{k}}({r_{i+k}+r_{i+k-\frac12}\phi}),
		& & \\
	&=r_i+u( r_{i+1,\bar 0}+ r_{i+\frac 12,\bar 0}\phi),\\
r_{i-\frac12,\bar 1}
	&=r_{i-\frac12}+\sum_{k\ge 1} {u^{k}}(r_{i+k-\frac12}
			+r_{i+k}\frac{\partial}{\partial \phi}),
		& &\quad \\
	&=r_{i-\frac12}+u ( r_{i+1-\frac 12,\bar 1}+
		r_{i+1,\bar 1}\frac{\partial}{\partial \phi}).
\end{aligned}
\end{equation}
Let $X$ be an even matrix. Because of the multiplicative
property of Berezinians we can add multiples of a row to another row
of $X$
without changing $\ber\,(X)$ and
$\berdual(X)$. Using such row operations we see, using
\thetag{\ref{eq:tildewinw}}, that
\begin{align*}
\ber\,([Q_{\bar 0}(u,\phi)\mathcal W]_-)
	&=\ber\,(A_0(r_{0,\bar 0})),\\
\berdual \,([Q_{\bar 1}(u,\phi)\mathcal W]_-)
	&=\berdual\,(A_{-\frac 12}( r_{-\frac 12,\bar 1})).
\end{align*}
Now $\ber$ is linear in even rows, and $\berdual $ in odd
rows, so by \thetag{\ref{eq:tildewinw}} we find
\begin{multline*}
\ber\,([Q_{\bar 0}(u,\phi)\mathcal
W]_-)=\ber\,(A)+
\\\sum_{i>0 }u^{i}[\ber\,(A_0(r_i))+\ber\,(A_0(r_{i-\frac12}))\phi],
\end{multline*}
and
\begin{multline*}
\berdual \,([Q_{\bar1}(u,\phi)\mathcal W]_-) =\berdual \,(A)+\\
	+\sum_{i>0}\,u^{i}[\berdual\,(A_{-\frac12}(r_{i-\frac12}))+
	\berdual\,(A_{-\frac12}(r_i))\frac{\partial}
					{\partial\phi}] .
\end{multline*}
Comparing with \thetag{\ref{eq:bakerfunction}} proves the
lemma.
\end{proof}

We now consider the flow on $\Sgr$ generated by the negative part of
the
Jacobian Heisenberg algebra: define
\begin{equation}\label{expfactor}
\gamma(t)=\exp(\sum_{{i>0}} t_i z^{-i}+ t_{i-\frac12}z^{-i}\theta),
\quad t_i
\in \Lambda_{\text{ev}},t_{i-\frac12}\in\Lambda_{\text{odd}}
\end{equation}
and put for $W\in \Sgr$:
$$W(t)=\gamma(t)\inv W.$$
The $\tau$-functions associated to a point $W$ in the big cell are
then
functions on $\Jheis_{-,\text{ev}}$:
\begin{equation}\label{eq:tau}
\tau_W(t),\tau^*_W(t):\Jheis_{-,\text{ev}}\to \Lambda\cup \{\infty\}
\end{equation}
given by
\begin{align}\label{eq:taudef}
\tau_W(t)
	&=\frac{\sigma(\gamma(t)\inv W)}{\gamma(t)\inv \sigma (W)}=
	\frac{\ber( [\gamma(t)\inv\circ \mathcal
W]_-)}{\ber([\mathcal
		W]_-)},\\
\tau_W^*(t)
	&=\frac{\sigma^*(\gamma(t)\inv W)}{\gamma(t)\inv \sigma^*
(W)}=
	\frac{\berdual( [\gamma(t)\inv\circ \mathcal
W]_-)}{\berdual([\mathcal
		W]_-)}.
\end{align}
Here $\sigma$ and $\sigma^*$ are the sections of $\Bersgrdual$ and
$\Bersgr$
defined in \thetag{\ref{eq:defsigma*}} and $\gamma\inv\in
Gl_{\infty\mid\infty}(\Lambda)$ acts via \thetag{\ref{eq:lift}} on
$\Bersgr$
and $\Bersgrdual$.

The Baker function of $W$ becomes now a function on
$\Jheis_{-,\text{ev}}$, and
we have
an expression in terms of a quotient of (shifted) $\tau$-functions:
\begin{equation} \label{eq:Bakertauquotient}
w_{\bar 0}(t;u,\phi)=\frac{\tau_W(t;Q_{\bar 0})}{\tau_W(t)},\quad
w_{\bar 1}(t;u,\phi)=\frac{\tau_W^*(t;Q_{\bar 1})}{\tau^*_W(t)},
\end{equation}
where
$$
\tau_W(t;Q_{\bar 0})=\frac{\ber([Q_{\bar0}\mathcal
W]_-)}{\ber(\mathcal W_-)},
\quad
\tau^*_W(t;Q_{\bar 1})=\frac{\berdual([Q_{\bar1}\mathcal
W]_-)\phi}{\berdual(\mathcal W_-)}.
$$
Note that even if we are only interested in the Jacobian Heisenberg
flows the
full Heisenberg flows automatically appear in the theory if we
express the
Baker functions in terms of the $\tau$ functions.
In principle we could also consider the flows on $\Sgr$ generated by
the full
super Heisenberg algebra $\Sheis$. However, since $\Sheis$ is
non--Abelian the
interpretation of these flows is less clear and therefore we leave
the
discussion of these matters to another occasion.

\section{The Krichever map and algebro-geometric solutions}
\subsection{The Krichever map.}\label{ss:Krichever}

Consider now a set of geometric data
$(X,P,(z,\theta),\mathcal{L},t)$, where:
\begin{itemize}
\item $X$ is a generic SKP curve as before.
\item $P$ is an irreducible divisor on $X$, so that $P^{\text{red}}$
is a single point of the underlying Riemann surface $X^{\text{red}}$.
\item $(z,\theta)$ are local coordinates on $X$ near $P$, so that
$P$ is defined by the equation $z=0$.
\item $\mathcal{L}$ is an invertible sheaf on $X$.
\item $t$ is a trivialization of $\mathcal{L}$ in a neighborhood of
$P$,
say $U_P = \{| z^{\text{red}}|<1\}$.
\end{itemize}
We will associate to this data a point of the super Grassmannian
$\Sgr$.

For studying meromorphic sections of $\mathcal L$ we have the exact
sequence
\begin{equation}
0 \rightarrow \mathcal{L} \overset{\text{inc}}{\rightarrow}
\mathcal{L}(P)
\overset{\text{res}}{\rightarrow} \mathcal{L}_{P^{\text{red}}} \cong
\Lambda |
\Lambda \rightarrow 0,
\end{equation}
which gives
\begin{equation} \label{polesequence}
H^0(\mathcal{L}) \hookrightarrow H^0(\mathcal{L}(P)) \rightarrow
\Lambda | \Lambda \rightarrow H^1(\mathcal{L}) \rightarrow
H^1(\mathcal{L}(P)) \rightarrow 0,
\end{equation}
where the residue is the pair of coefficients of $z^{-1}$ and $\theta
z^{-1}$
in the Laurent expansion.

Let $\mathcal{L}(*P) = \lim_{n \rightarrow \infty} \mathcal{L}(nP)$
be
the sheaf of sections of $\mathcal{L}$ holomorphic except possibly
for a pole of arbitrary order at $P$.
The {\it Krichever map} associates to a set of geometric data
as above the $\Lambda$-module of formal Laurent series
$W = z\ t[H^0(X,\mathcal{L}(*P))]$, which will be viewed as a
submodule of $H$.

In \cite{MuRa:SupKrich,Ra:GeomSKP} the concern was expressed that $W$
might not be freely generated, and hence not an element of $\Sgr$ as
we have defined it.
However,

\begin{thm} $H^0(X,\mathcal{L}(*P))$ is a freely generated
$\Lambda$-module, and $W \in \Sgr$.
Further, $W$ belongs to the big cell if the geometric data satisfy
$H^0(X,\mathcal{L}) = H^1(X,\mathcal{L}) = 0$, which happens
generically if $\deg \mathcal{L} = g-1$.
\end{thm}

\begin{proof}
Assume first that $H^0(X,\mathcal{L}) = H^1(X,\mathcal{L}) = 0$.
Then the sequence \thetag{\ref{polesequence}} applied to
$\mathcal{L}$ gives
\begin{equation}
0 \rightarrow H^0(\mathcal{L}(P)) \rightarrow \Lambda | \Lambda
\rightarrow 0 \rightarrow H^1(\mathcal{L}(P)) \rightarrow 0,
\end{equation}
so that $H^0(X,\mathcal{L}(P))$ is freely generated by an even and
odd section having principal parts $z^{-1}$ and $\theta z^{-1}$, and
$H^1(X,\mathcal{L}(P))$ is still zero.
Applying the same sequence inductively to $\mathcal{L}(nP)$ shows
that
$H^0(X,\mathcal{L}(*P))$ is freely generated by one even and one odd
section of each positive pole order.
  So $W$ is obtainable from $H_-$ by multiplication by a lower
triangular
invertible matrix, and  $W$ belongs to the big cell of $\Sgr$.
We also have $H^i(\mathcal{L}\spl) = 0, i=0,1$, from Theorem
\ref{thm:freeness}.
And, by the super Riemann-Roch Theorem \thetag{\ref{superRR}},
$\deg \mathcal{L} = \deg \mathcal{L}\spl = g-1$.
Moreover, by semicontinuity, in $\Pic^{g-1}(X)$ the cohomology groups
$H^i(\mathcal{L})$ can only get larger on Zariski
closed subsets, so generically they are zero.

Now consider the general situation in which $H^i(\mathcal{L})$ may
not be zero.
Still, by twisting,
$H^1(\mathcal{L}(nP)) = H^1(\mathcal{L}\spl(nP)) = 0$ for $n$
sufficiently
large.
Then, by the previous argument, $H^0(\mathcal{L}(*P))$ has
non purely nilpotent elements  with poles of order $n+1$ and higher;
the worry
is
that one may only be able to find nilpotent generators for
$H^0(\mathcal{L}(nP))$.
So take $f \in H^0(\mathcal{L}(nP))$ of order $k$ in nilpotents:
its image in $H^0(\mathcal{L}(nP)/\mathfrak m^k)$ is zero, but its
image
$\hat{f}$ in
$H^0(\mathcal{L}(nP)/\mathfrak m^{k+1})$ is nonzero
and also lies in $\Lambda^k H^0(\mathcal{L}\spl(nP))$.
Then $\hat{f}$ can be identified with a sum of elements $f_a$
of $H^0(\mathcal{L}\spl(nP))$ with coefficients from $\Lambda^k$.
By the extension sequence \thetag{\ref{eq:longexactcohom}}, each
$f_a$ can
be extended order by order in nilpotents to an element of
$H^0(\mathcal{L}(nP))$ which is not purely nilpotent.
So we can write the order $k$ element $f$ as a $\Lambda$-linear
combination of
not purely nilpotent elements of $H^0(\mathcal{L}(nP))$, modulo an
element of
order
$k+1$. Induction on $k$ shows then that any element of
$H^0(\mathcal{L}(nP))$
is a
$\Lambda$-linear combination of not purely nilpotent elements of
$H^0(\mathcal{L}(nP))$. So there
exists a set of non-nilpotent
elements which span $H^0(\mathcal{L}(nP))$ over $\Lambda$.
A linearly independent subset of these completes a  basis for
$H^0(\mathcal{L}(*P))$.
\end{proof}

\subsection{The Chern class of the Ber bundle on
$\Picxz$.}\label{ss:ChernclassBeronPic}
By the arguments of the previous subsection we have, in case $W\in
\Sgr$ is
obtained by the Krichever map from geometric data $(X, P, \mathcal L,
(z,\theta), t)$, an exact sequence of $\Lambda$-modules:
\begin{equation}\label{eq:seqcohomW}
0\to H^0(X,\mathcal L)\to W\to H_-\to H^1(X,\mathcal L)\to 0.
\end{equation}
We can interpret the Ber bundle $\Bersgr$ in terms of this sequence
as follows.
Let $M$ be a free $\Lambda$-module, possibly of infinite rank, and
let
$B=\{\mu\}$ be  a collection of bases for $M$ such that any two bases
$\mu,\mu^\prime\in B$ are related by $\mu^\prime=\mu T$ where $T\in
\operatorname{Aut}(M)$ has a well defined Berezinian. Then we
associate to the
pair $(M,B)$ a free rank $(1\mid 0)$ module $\ber(M)$ with generator
$b(\mu)$
for any $\mu\in B$ with identification $b(\mu^\prime)=\ber(T)b(\mu)$.

The fiber of $\Bersgr$ at $W$ can then be interpreted as $\ber(W)$,
using the
collection of admissible bases as $B$ in the above definition.
Similarly we can
construct on $\Sgr$ a line bundle with fiber at $W$ the module
$\ber(H_-)$.
Clearly this bundle is trivial, so we can, even better, think of
$\Bersgr$ as
having fiber $\ber(W)\otimes \berdual(H_-)$. But by the properties of
the
Berezinian we get from \eqref{eq:seqcohomW}
$$
\ber(W)\otimes \berdual(H_-)=\ber(H^0(X,\mathcal L))\otimes
\berdual(H^1(X,\mathcal L)).
$$
Now we have seen in subsection \ref{subs:Berbundles} that the Ber
bundle
$\Ber(\Picxz)$ on $\Picxz$ has  the same fiber, with the difference
that there
we were dealing with bundles of degree 0 and here $\mathcal L$ has
degree
$g-1$.

For fixed $(X,P,(z,\theta))$ the collection  $M$ consisting of
Krichever data
$(X,P,(z,\theta),\mathcal L,t)$ forms a supermanifold and we have two
morphisms
$$
i:M\to \Sgr,\quad p:M\to \Picxz
$$
where $i$ is the Krichever map and $p$ is the projection from
Krichever data to
the line bundle $\mathcal L$. (Here we identify $\Pic^n(X)$ with
$\Picxz$ via
the invertible sheaf $\Ox(-nP)$.) Then we see that
$i^*(\Bersgr)\simeq
p^*(\Ber(\Picxz))$. This fact allows us to prove Theorem \ref
{thm:1ChernBertriv}.

Note first that we have a surjective map
\begin{equation}\label{eq:surjecttoH1}
\Jheis_-\to H^1(X,\Ox).
\end{equation}
Indeed, let $X=U_0\cup U_P$ be an open cover where
$U_0=X-P^{\text{red}}$ and
$U_P$ is a suitable disk around $P^{\text{red}}$. Then if $[a]\in
H^1(X,\Ox)$
is represented by $a\in\Ox(U_0\cap U_P)$  we can write, using the
local
coordinates on $U_P$,
$a=a_P+\sum_{i>0} a_i z^{-i}+\alpha_i z^{-i}\theta$, with $a_P\in
\Ox(U_P)$.
Then
$a-a_P=\sum a_i z^{-i}+\alpha_i z^{-i}\theta\in \Jheis_-$ and
$[a]=[a-a_P]$.
Now the tangent space to any point $\mathcal L\in\Picxz$ can be
identified with
$H^1(X,\Ox)$ and so we have a surjective map from $\Jheis_-$ to the
tangent
space of $\Picxz$. Note secondly that a change of trivialization of
$\mathcal
L$, given by $t\mapsto t^\prime$, corresponds to multiplication of
the point
$W\in \Sgr$ by an element $a_0+\alpha_0\theta +\sum_{i>0}a_i z^{
i}+\alpha_i
z^{i}\theta$ of the group corresponding to $\Jheis_+$. From these two
facts we
conclude that there is a surjective map from $\Jheis$ to the tangent
space  to
the image of the Krichever map $i:M\to \Sgr$ at any point
$W=W(X,P,(z,\theta),\mathcal L,t)$. Now the first Chern class of
$\Bersgr$ is
calculated from the cocycle \eqref{eq:cocycle} on
$gl_{\infty\mid\infty}(\Lambda)$ and it follows from Lemma
\ref{lem:comactionJheis} that the restriction of this cocycle to
$\Jheis$ is
identically zero. This implies that
$$
i^*(c_1[\Bersgr])=p^*(c_1[\Ber(\Picxz)])=0.
$$
But the map $p:M\to \Picxz$ is surjective, so we finally find that
$c_1(\Ber(\Picxz))=0$ and $\Ber(\Picxz)$ is topologically trivial,
proving
Theorem \ref{thm:1ChernBertriv}.

\subsection{Algebro-geometric tau and Baker functions.}
\label{ss:AlgebrogeometrictauBaker}

We consider geometric data mapping to $W$ in the big cell of $\Sgr$,
so that $\deg \mathcal{L} = g-1$.
As discussed in Section 3, we can associate to $W$ both a tau
function
and a Baker function.
A system of super KP flows on $\Sgr$ applied to $W$ produces an orbit
corresponding
to a family of deformations of the original geometric data.
The simplest system of super KP flows, the ``Jacobian" system of
Mulase and
Rabin
\cite{Mu:Jac,Ra:GeomSKP}, deforms the geometric data by
moving $\mathcal{L}$ in $\text{Pic}^{g-1}(X)$.
Solutions to this system for $X$ a super elliptic curve were
obtained in terms of super theta functions in \cite{Ra:SupElliptic}.
On the basis of the ordinary KP theory, cf. \cite{SeWi:LpGrpKdV},
section 9,
we might expect that in general the tau and Baker
functions for this family can be given explicitly as functions of the
flow
parameters by means
of the super theta functions (when these exist) on the Jacobian  of
$X$.
We now discuss the extent to which this is possible.

Recall from \eqref{eq:surjecttoH1}
 that we have a surjection from $\Jheis_{-,\text{ev}}$ to the
cohomology group
$H^1(X,\Oxev)$. By exponentiation we obtain a map from
$\Jheis_{-,\text{ev}}$
to $\Picxz$ and these maps fit together in a diagram
\begin{equation}
\begin{CD}\label{eq:cdtautheta}
{}	@.	{}	@.	{}	@.	0	@.	{}
\\
@.		@.		@.		@VVV		@.
\\
{}	@.	0	@.	{}	@.	H^1(X,\mathbb Z)@.
\\
@.		@VVV		@.		@VVV		@.
\\
0  @>>>	K_0	@>>>\Jheis_{-,\text{ev}}@>>> H^1(X,\Oxev)	@>>>0
\\
@.	@VVV		@\vert		@VVV		@.	\\
0  @>>>	K	@>>>\Jheis_{-,\text{ev}}@>>>\Picxz	@>>>0	\\
@.		@VVV		@.		@VVV		@.
\\
{}	@.  K/K_0	@.	{}	@.0		@.{}\\
@.		@VVV		@.		@.		@.
\\
{}	@.	0	@. 	{}	@.	{}	@.	{}
\\
\end{CD}
\end{equation}
Here $K_0$ is the $\Lambda_{\text{ev}}$-submodule of elements $f$ of
$\Jheis_{-,\text{ev}}$ that split as $f=f_0 + f_P$, with $f_0\in
\Ox(U_0)$ and
$f_P\in \Ox(U_P)$ and $K$ is the Abelian subgroup (not submodule!) of
elements
$k$ of $\Jheis_{-,\text{ev}}$ that after exponentiation factorize:
$e^k=\phi_k
e^{k_p}$, with $\phi_k\in \Ox(U_0)^\times$ and $k_P\in \Ox(U_P)$.
{}From the
Snake Lemma it then follows that $H^1(X,\mathbb Z)\simeq K/K_0$. So a
function
$\hat F$ on $\Jheis_{-,\text{ev}}$ descends to a function $ F$ on
$H^1(X,\Oxev)$ if it is invariant under $K_0$.
The automorphic behavior of such a function  $ F$ with respect to the
lattice
$H^1(X,\mathbb Z)$ translates into behaviour of
$\hat{F}$ under shifts by elements of $K$.
In particular we consider the function $\tau_W$ associated to a point
$W$ in
the big cell of $\Sgr$, see \eqref{eq:tau} and \eqref{eq:taudef}.
This is a function on $\Jheis_{-,\text{ev}}$ and, because of Lemma
\ref{lem:comactionJheis}, we see by an easy adaptation of the proof
of Lemma
9.5 in \cite{SeWi:LpGrpKdV} that
$$
\tau_W(f+k)=\tau_W(f)\tau_W(k),\quad f\in \Jheis_{-,\text{ev}}, k\in
K.
$$
In particular we obtain by restriction a homomorphism
$$
\tau_W:K_0\to \Lambda^\times_{\text{ev}}.
$$
Let $\eta:K_0\to \Lambda_{\text{ev}}$ be a homomorphism such that
$\tau_W(k_0)=e^{\eta(k_0)}$, for all $k_0\in K_0$. Then we can define
a new
function
$$
\hat \tau_1(f)=\tau_W(f)e^{-\eta(f)}.
$$
Then $\hat \tau_1(k_0)=1$, but still we have
\begin{equation}\label{eq:multtau}
\hat \tau_1(f+k)=\hat \tau_1(f)\hat \tau_1(k),
\end{equation}
 so that $\hat \tau_1$ descends to a function $\tau_1$ on
$H^1(X,\Oxev)$. From
 \eqref{eq:multtau} we see that $\tau_1$ corresponds to a
(meromorphic) section
 of a line bundle on $\Picxz$ with trivial Chern class.

A suitable ratio of translated theta functions gives a section of
this same bundle, so that $\tau_1$ is expressed as this ratio times
a meromorphic function, the latter being rationally expressible
in terms of super theta functions.
Then the modified tau function $\tau_1$ is rationally expressed in
terms of
super theta functions.

The even Baker function $w^W_{\bar 0}(z,\theta)$ associated to the
point $W$
is just
the even section of $\mathcal L$ holomorphic except for a pole $1/z$
at $P$.
Such a section can be specified by its restrictions to the charts
$U_0$ and $U_P$.
The Jacobian super KP flows act by multiplying the transition
function of
$\mathcal L$ across
the boundary of $U_P$ by a factor
$\gamma(t)$ as in \eqref{expfactor}.
The corresponding action on the associated point $W$ of $\Sgr$ is
generated by
the matrices $\lambda(n)+\mu(n)$ and
$f(n)$ of Section 3; the remaining matrices generate deformations of
the curve
$X$
and enter the Kac--van de Leur SKP flows.
Then $w^{W(t)}_{\bar 0} / w^W_{\bar 0}$ is a section of the bundle
with
transition function \thetag{\ref{expfactor}}.
Equivalently, it is a meromorphic function on $U_0$ which extends
into $U_P$
except for an
essential singularity of the form \thetag{\ref{expfactor}}, having
zeros at the
divisor of $\mathcal L(t)$
and poles at the divisor of $\mathcal L$.
By analogy with the ``Russian formula" of ordinary KP theory, such a
function
would be expressed in the form
\begin{equation} \label{Russian}
\exp
[
	\sum_{k=1}^{\infty}
	\int_{(0,0)}^{(z,\theta)}
	(t_k \hat{\psi}_k + t_{k-\frac12} \hat{E}_k)
	+ c(t)
]
\end{equation}
times a ratio of theta functions providing the zeros and poles.
Here $\hat{\psi}_k$ and $\hat{E}_k$ are differentials on $\hat{X}$,
with vanishing $a$-periods and holomorphic except for the behavior
near $P$,
\begin{equation}
\hat{\psi}_k \sim \hat{D}(z^{-k}) = -k\rho \hat{z}^{-(k+1)},\;\;\;\;
\hat{E}_k \sim \hat{D}(\theta z^{-k}) = \hat{z}^{-k}.
\end{equation}
The constant $c(t)$ is linear in the flow parameters.
In addition to the symmetry of the period matrix, we have to require
the
existence of these differentials.
This requires that they exist in the split case, and then that these
split
differentials extend through
the sequence \thetag{\ref{eq:longexactcohom}}.
In the split case, the odd differentials $\hat{\psi}_k$ are just
$\theta$ times
the ordinary differentials
on the reduced curve which appear in the Russian formula (and they do
extend).
However, the even differentials $\hat{E}_k$ are sections of $\mathcal
N$, which
is of degree zero
and nontrivial, with $h^1 = g-1$.
Consequently, when $g>1$ there will be Weierstrass gaps in the list
of pole
orders of these differentials.
This means that the odd flow parameters corresponding to the missing
differentials
must be set to zero in order for the Baker function to assume the
``Russian"
form.
Even then, however, the function given by the Russian formula will
generically
behave as $1 + \alpha\theta + \mathcal O(z)$ for
$z \rightarrow 0$, rather than the correct $1 + O(z)$ for
$w^{W(t)}_{\bar 0}$
containing no $\theta/z$ pole.
In \cite{Ra:SupElliptic} this was dealt with by including a term $\xi
\hat{E}_0$ in the exponential, taking
$\partial_\xi$ to construct a section with a pure $\theta/z$ pole,
and
subtracting off the appropriate multiple of this.
In general, however, no such $\hat{E}_0$ will exist.
These difficulties are understandable in view of the relations
(\ref{eq:Bakertauquotient}) which require that the tau function
be known for the full set of K-vdL flows in order to compute
the Baker functions for even the Jacobian flows.
Since the dependence of the tau function on the non-Jacobian flows is
likely to
be far more complicated than
our super theta functions, it is unlikely that the
Baker functions can be expressed in terms of them.

\appendix
\section{Duality and Serre duality}\label{app:dualSerredual}

Let  $\Lambda$ be our usual ground ring $\mathbb
C[\beta_1,\dots,\beta_n]$ with the $\beta_i$ odd
indeterminates. In this Appendix we will discuss duality for
$\Lambda$-modules (cf. chapter 21 in Eisenbud,
\cite{Eis:ComAlg} for the case of commutative rings). Then we
will use this to extend  Serre duality for supermanifolds
over the ground field $\mathbb C$
(\cite{HaskeWells:Serreduality}, \cite{OgPenk:SerreDual}) to
supermanifolds over the ground ring $\Lambda$.

Finally we discuss the more explicit form of Serre duality
that one has in case of super curves. One way of proving Serre
duality for super curves over general $\Lambda$ would be by
using the properties of the supertrace in infinite rank
$\Lambda$-modules to define a residue, generalizing the
method of Tate, \cite{tate:residue}. However, we have
available in our case contour integration which also provides
us with a residue map.

\subsection{Duality of $\Lambda$-modules}\label{appss:dualLmod}
Let $M$ be an object of $\mathcal {M}\text{od}_\Lambda$, the
category of finitely generated $\Lambda$-modules. We give the
$\mathbb Z_2$-graded vector space $E(M)=\text{Hom}_{\mathbb
C}(M,\mathbb C)$ the structure of a $\Lambda$-module via
$$
\lambda \cdot \phi (m)=(-1)^{|\lambda| |\phi|}\phi(\lambda m),
$$
for all $m\in M$ and all homogeneous $\lambda\in \Lambda$ and
$\phi\in E(M)$. Then one checks that $E=E(-)$ is a {\it
dualizing functor} on $\mathcal {M}\text{od}_\Lambda$: it is
contravariant, exact, $\Lambda$-linear and satisfies
$E^2\simeq 1_{\mathcal {M}\text{od}_\Lambda}$. One also
checks that $E(\Lambda)\simeq\Lambda$, up to a possible parity
change. In the sequel we will ignore these parity changes.

An explicit (basis dependent) isomorphism $\Lambda \to
\text{Hom}_{\mathbb C}(\Lambda,\mathbb C)$ can be described as
follows. The monomials
$\beta^I=\beta_1^{i_1}\dots\beta_n^{i_n}$, $i_j=0,1$,  form a
basis of $\Lambda$ as $\mathbb C$-vector space, and we let
$\phi_{\beta^I}$ be the dual basis of
$E(\Lambda)=\text{Hom}_{\mathbb C}(\Lambda,\mathbb C)$, so
that $\phi_{\beta^I}(\beta^J)=\delta_{IJ}$. Then  the
$\Lambda$-homomorphism  that maps $1\in \Lambda$ to the
linear functional $\phi_{\beta_1\dots\beta_n}$ is an
isomorphism, odd in case $n$ is odd.

 From the fact that $E$ is a dualizing functor we see that the
map $\text{Hom}_{\Lambda}(M,N)\to
\text{Hom}_{\Lambda}(E(N),E(M))$ is an isomorphism for all
objects $M,N$ and hence
$$
E(M)=\text{Hom}_{\Lambda}(\Lambda,E(M))\simeq\text{Hom}_{\Lambda}(M,
E(\Lambda))=  \text{Hom}_{\Lambda}(M,\Lambda).
$$
In other words, up to a possible parity switch, we can
identify (functorially) the $\Lambda$-modules of $\mathbb
C$-linear and of $\Lambda$-linear homomorphisms on any
finitely generated $\Lambda$-module $M$:
$$
\text{Hom}_{\mathbb C}(M,\mathbb C)\simeq
\text{Hom}_{\Lambda}(M,\Lambda)
$$
We will identify the two and use for both the symbol $M^*$.

 Note that the above  implies that the functor that maps $M$ to
 $M^*=\text{Hom}_{\Lambda}(M,\Lambda)$ on $\mathcal
{M}\text{od}_\Lambda$ (naturally isomorphic to the functor
$E$) is exact: exact sequences get mapped to exact sequences.
Equivalently: $\Lambda$ is injective as a module over itself.
Also note that the double dual of $M$ is isomorphic to $M$
itself.

\subsection{Serre duality of
supermanifolds}\label{appss:SerredualSupermanifold}
Serre duality
for supermanifolds $(Y,\mathcal O_Y)$ of dimension $(p\mid
q)$ over the complex numbers is discussed in Haske and Wells
\cite{HaskeWells:Serreduality} and in
Ogievetsky and Penkov \cite{OgPenk:SerreDual}.

A {\it dualizing sheaf} on $(Y,\mathcal O_Y)$ is an
invertible sheaf $\omega_Y$ together with a fixed
homomorphism
$$
t: H^p(X,\omega_Y)\to \mathbb C
$$
such that the induced pairing  for all $\mathcal F$
$$
H^i(Y,\mathcal F^*\otimes \omega_Y)\otimes H^{p-i}(Y,\mathcal
F)\overset{t}\to
\mathbb C
$$
gives an isomorphism, called Serre duality,
\begin{equation}
\label{eq:Serredual}
H^i(Y,\mathcal F^*\otimes
\omega_Y)\overset{\sim}\to  H^{p-i}(Y,\mathcal F)^*.
\end{equation} Note that this is an isomorphism of $\mathbb
Z_2$-graded vector spaces over $\mathbb C$. Dualizing sheaves
are unique, up to isomorphism.

If $\mathcal T_{Y}$ is the tangent sheaf of $(Y,\mathcal
O_Y)$ with transition functions $J_{\alpha\beta}:U_\alpha\cap
U_\beta\to Gl(p\mid q,\mathbb C)$, then the {\it Berezinian
sheaf} $\mathcal B\text{er}_Y$ is the invertible sheaf with
transition functions $\ber\,(J_{\alpha\beta}):U_\alpha\cap
U_\beta\to \mathbb C^\times$. In
\cite{HaskeWells:Serreduality}, \cite{OgPenk:SerreDual} it is
proved  that $\mathcal B\text{er}_Y$ is a dualizing sheaf.

We are interested in the relative situation: supermanifolds
$(Z,\mathcal O_Z)\to (\bullet,\Lambda)$ of dimension $(p\mid
q)$ over $\Lambda$. The structure sheaf $\mathcal O_Z$
contains $n$ independent global odd constants. We can reduce
this to the absolute case by thinking of $(Z,\mathcal O_Z)$
as a supermanifold $(Z,\mathcal O_Y)$ of dimension $(p\mid
q+n)$ over  $\mathbb C$: the constants $\beta_i$ generating $\Lambda$
are now interpreted as coordinates.

The tangent sheaf $\mathcal T_Y$ of derivations of $\mathcal
O_Y$ has then global sections $\frac{\partial}{\partial
\beta_I}$, in contrast to the relative tangent sheaf
$\mathcal T_{Z/\Lambda}$.  One sees easily, however, that the
invertible sheaf $\mathcal B\text{er}_Y$, constructed from
the tangent sheaf $\mathcal T_Y$, is the same as $\mathcal
B\text{er}_{Z/\Lambda}$, constructed from the relative
tangent sheaf $\mathcal T_{Z/\Lambda}$.

Now let $\mathcal F$ be a coherent, locally free sheaf  of
$\mathcal O_Z$-modules. Then the associated cohomology groups
are finitely generated $\Lambda$-modules, so we can use the
theory of Appendix \ref{appss:dualLmod}. In particular we see
that there is a homomorphism
\begin{equation}\label{eq:trace}
t: H^p(Z,\mathcal B\text{er}_{Z/\Lambda})\to \Lambda
\end{equation}
inducing an isomorphism of $\Lambda$-modules, also called
Serre duality,
\begin{equation}\label{eq:LSerredual}
H^i(Z,\mathcal F^*\otimes\mathcal
B\text{er}_{Z/\Lambda})\overset{\sim}\to  H^{p-i}(Z,\mathcal
F)^*,
\end{equation}
where now $H^p(Z,\mathcal F)^*$ means the
$\Lambda$-linear dual.

In case of $N=1,2$ super curves $(X,\Ox)\to (\bullet,\Lambda)$
we can be more explicit about the homomorphism $t$ in
\eqref{eq:trace}. Write $\Berx$ for the relative
Berezinian $\mathcal B\text{er}_{X/\Lambda}$.  The cohomology
of $\Berx$ is calculated by the sequence
$$
0\to H^0(X,\Berx)\to \text{Rat}(\Berx)\to \text{Prin}(\Berx)\to
H^1(X,\Berx)\to
0,
$$
where $\text{Rat}(\Berx)$ are the rational sections and
$\text{Prin}(\Berx)$ the principal parts. On
$\text{Prin}(\Berx)$ we have for every $x\in X$ a residue map
$\text{Res}_x\omega\mapsto \frac1{2\pi i}\oint_{C_x}\omega$,
where $C_x$ is a contour around $x$, using the integration on
$N=1,2$ curves introduced in subsection
\ref{ss:IntegraOnSupCurve}. Then $t=\sum_{x\in X}
\text{Res}_x$ and the classical residue theorem holds: if
$\omega \in \text{Rat}(\Berx)$ then $t(\omega)=0$. Now that we
have the residue we can proceed to prove Serre duality as in
the classical case, cf. \cite{Se:AlgGrpsClassFields}.

\section{Real structures and conjugation.}\label{ss:realstrconj}
Let $\Lambda=\mathbb C\,[\beta_1\dots,\beta_n]$ be the Grassmann
algebra
generated by $n$ odd indeterminates. Choose a sign $\epsilon=\pm 1$
and define
a real structure on $\Lambda$ for this choice as a real-linear, even
map
$\omega:\Lambda\to \Lambda$
such that
\begin{enumerate}
\item$\omega(c a)=\bar c \omega(a),\quad c\in \mathbb C,a\in\Lambda
$,
\item $\omega^2=1$.
\item $\omega(\beta_i)=\beta_i$,
\item $\omega(a b)=\epsilon^{|a|\,|b|}\omega(b)\omega(a)$,
\end{enumerate}
where $a,b$ are homogeneous elements of parity $|a|, |b|$.
We will often write $\bar a$ for $\omega(a)$.

A real structure induces a decomposition of $\Lambda$ into
eigenspaces,
$\Lambda=\Lambda_{\Re}\oplus\Lambda_{\Im}$, where
$\Lambda_{\Re}$ is the $+1$ eigenspace and $\Lambda_{\Im}$ the $-1$
eigenspace
for $\omega$. Multiplication by $i$ is an isomorphism, so we get
$\Lambda=\Lambda_{\Re}\oplus i\Lambda_{\Re}$.

\section{Linear equations in the super
category.}\label{app:Lineqsupercat}
\subsection{Introduction.}
Cramer's rule tells us that a system of linear equations over
a commutative ring $R$:
\begin{equation}\label{eq:lineq}
x A=y,
\end{equation}
with $x=(x_1,\dots,x_n)$ and $y=(y_1,\dots,y_n)$
$n$-component row vectors  and $A$ an $n\times n$ matrix with
coefficients in $R$,
can be solved by quotients of determinants:
\begin{equation}\label{eq:Cramer}
x_i=\frac{\det\,(A_i(y))}{\det\,(A)},\quad i=1,\dots,n,
\end{equation}
where $A_i(y)$ is the matrix obtained by replacing in $A$ the
$i$th row by the row vector $y$.

We want to study \thetag{\ref{eq:lineq}} in the super
category. So we fix a decomposition $n=k+l$. Call an index
$i\in\{1,\dots,n\}$ {\it even} if $i\le k$ and {\it odd}
otherwise.  Consider an even  matrix $A=(a_{ij})_{i,j=1}^n$
over $\Lambda$: $a_{ij}$ is an element of the even part
$\Lambda_{\text{ev}}$
if  $i$ and $j$ have the same parity and  of the odd part
$\Lambda_{\text{odd}}$ otherwise. Note that it is not
necessary to specify the parity of  $y$ in
\thetag{\ref{eq:lineq}}: it can be even, odd or
inhomogeneous.

In the theory of linear algebra over a Grassmann algebra
$\Lambda$ in many (but not all) respects a role
analogous to that of the determinant in the commutative case
is played by the {\it Berezinian} and its inverse: for an
even matrix $A$ as above we define
\begin{equation}\label{eq:defberber*}
    \begin{split}\ber\,(A)&=\det\,(X-\alpha
Y\inv\beta)\det\,(Y\inv),\\
    \ber^*(A)&=\det\,(X\inv)\det\,(Y-\beta X\inv
\alpha)=\frac1{\ber(A)}
\end{split}
\end{equation}

We will discuss how Cramer's rule can be generalized in this
setting. The result is as follows: the solution of
\thetag{\ref{eq:lineq}} for $A$ an even super matrix is
given by
\begin{equation}\label{eq:supercramer}
x_i=\begin{cases}{\ber\,(A_i(y))}/{\ber\,(A)}&\text{if $i\le
k$},\\
		  {\ber^*\,(A_i(y))}/{\ber^*\,(A)}&\text{if $i> k$}.
		 \end{cases}
\end{equation}
Note that in general the matrix $A_i(y)$ is not even.
However, in $\ber\,(A)$ for $i\le k$ (resp. $\ber^*\,(A)$
for $i>k$) the entries $a_{ij}$, $j=1,\dots,n$ of row $i$
occur only linearly. So, if $\ber\,(A)$ (resp.
$\ber^*\,(A)$) exists,  we mean by $\ber\,(A_i(y))$ for
$i\le k$ (resp. $\ber^*\,(A_i(y))$ for $i> k$) the element
of $\Lambda$ obtained by replacing in $\ber\,(A)$ (resp.
$\ber^*\,(A)$) all $a_{ij}$ by $y_j$ for $j=1,\dots,n$.

The expression for the entries in even positions, i.e.,
$x_i$, $i\le k$, can be found in \cite{UeYaIk:AlgSKPOspKP},
but we haven't seen the solution for the odd positions in the
literature.
The main ingredient is the Gelfand-Retakh theory of quasi
determinants
\cite{GeRe:DetMatrnoncom,GeRe:Noncomdetchargraphs}.

\subsection{Quasideterminants.}
Let $A=(a_{ij})_{i,j=1}^n$ be an $n\times n$ matrix with entries
 independent (non commuting) variables $a_{ij}$. Let $A^{ij}$ be the
submatrix obtained by deleting in $A$ row $i$ and column
$j$. Following Gelfand-Retakh we introduce $n^2$ {\it
quasideterminants} $|A|_{ij}$, rational expressions in the
variables $a_{pq}$. If $n=1$
we put $|A|_{11}=a_{11}$ and for $n>1$ we define recursively
\begin{equation}\label{eq:defrecurs}
    |A|_{ij}=a_{ij}-\sum_{\frac{p\ne j}
{q\ne i}}a_{ip}|A^{ij}|_{qp}\inv a_{qj}.
\end{equation}
If we assume that the variables $a_{ij}$ commute amongst
themselves, then we have
\begin{equation}\label{eq:comvar}
    |A|_{ij}=(-1)^{i+j}\det\,(A)/\det\,(A^{ij}).
\end{equation}
Returning to the general case, one proves that the inverse
of $A$, as a matrix over the ring of rational functions in the
$a_{ij}$, is the matrix $B=(b_{ij})$, where
$b_{ij}=|A|_{ji}\inv$. This means that we can also define
the quasideterminant
as
\begin{equation}\label{eq:definv}
    |A|_{ij}=a_{ij}-\sum_{\frac{p\ne j}
{q\ne i}}a_{ip}c^{\,(ij)}_{pq} a_{qj},
\end{equation}
where $C^{\,(ij)}=(c^{\,(ij)}_{pq})$ is the inverse to the
matrix $A^{ij}$. This second definition is useful if one
wants to think of the entries of $A$ as elements of a ring
$R$. In that case it is perfectly well possible that the
inverse $C^{\,(ij)}$ of $A^{ij}$ exists,
but that some entry $c^{\,(ij)}_{pq}$ is not invertible in
$R$. In this situation $\thetag{
\ref{eq:definv}}$ allows us to define the quasideterminant
$|A|_{ij}$, whereas
$\thetag{\ref{eq:defrecurs}}$ might make no sense.

Quasideterminants have the following properties:

\begin{description}
\item[P1] If the matrix $B$ is obtained from $A$ by
multiplying row $i$ from the left
by $\lambda$ (a new independent variable) then for all
$j=1,2,\dots,n$ we have
\begin{equation*}
    |B|_{ij}=\lambda|A|_{ij},\quad |B|_{kj}=|A|_{kj}, k\ne  i.
\end{equation*}
Similarly if the matrix $C$ is obtained from $A$ by
multiplying column $j$ from the right by $\mu$ we have for
all $i=1,\dots,n$:
\begin{equation*}
    |C|_{ij}=|A|_{ij}\mu,\quad |C|_{ik}=|A|_{ik}, k\ne j.
\end{equation*}
\item[P2] If $B$ is obtained from $A$ by adding row $k$ to
some other row, then for all
$j=1,\dots,n$:
\begin{equation*}
    |B|_{ij}=|A|_{ij},\quad k\ne i.
\end{equation*}
If $C$ is obtained from $A$ by adding column $k$ to some
other column, then for all
$i=1,\dots,n$:
\begin{equation*}
    |C|_{ij}=|A|_{ij},\quad k\ne j.
\end{equation*}
\end{description}

\subsection{Restriction to super algebra.}
We keep the decomposition $n=k+l$ and we consider an even  matrix
$A=(a_{ij})_{i,j=1}^n$
over supercommuting variables: $a_{ij}$ is an even variable
if $i$ and $j$ have the same parity and an odd
variable otherwise.
\begin{lem}

\label{lem:quasidetqoutber}
Let $A$ be an even matrix as above. Then
\begin{itemize}
\item If $i,j$ are both even then
\begin{equation*}
|A|_{ij}=(-1)^{i+j}\ber\,(A)/\ber\,(A^{ij}).
\end{equation*}
\item If $i,j$ are both odd then
\begin{equation*}
|A|_{ij}=(-1)^{i+j}\ber^*\,(A)/\ber^*\,(A^{ij}).
\end{equation*}
\end{itemize}
\end{lem}

\begin{proof}We can write
\begin{equation*}
A=\begin{pmatrix}X&\alpha\\\beta&Y\end{pmatrix}=
\begin{pmatrix}1&\alpha
    	Y\inv\\0&1\end{pmatrix}
	\begin{pmatrix}X-\alpha Y\inv
\beta&0\\0&Y\end{pmatrix}
	\begin{pmatrix}1&0\\Y\inv\beta&1\end{pmatrix}.
\end{equation*}
This shows that $A$ is obtained from
$\begin{pmatrix}X-\alpha Y\inv \beta&0\\0&Y\end{pmatrix}$ by
operations that don't change the quasideterminant $|A|_{ij}$
in case $i,j$ are both even. Since the last matrix has
entries that commute
with each other we can use \thetag{\ref{eq:comvar}} to find
\begin{equation}\label{eq:qdetasqoutdet}
    |A|_{ij}=(-1)^{i+j}\frac{\det\,(X-\alpha Y\inv\beta)}
	{\det\,([X-\alpha Y\inv\beta]^{ij})}.
\end{equation}
 Now for $i,j$  both even $[X-\alpha
Y\inv\beta]^{ij}=X^{ij}-\alpha^{i\emptyset}Y\inv
\beta^{\emptyset j}$, where $\alpha^{i\emptyset}$ is the
matrix obtained by deleting just the row $i$ and
$\beta^{\emptyset j}$ is obtained by deleting column $j$. So
$\det\,([X-\alpha
Y\inv\beta]^{ij})\det(Y\inv)=\ber\,(A^{ij})$ and the first
part of the lemma follows by multiplying numerator and
denominator of \thetag{\ref{eq:qdetasqoutdet}} by
$\det\,(Y\inv)$. For the second part use the decomposition
\begin{equation*}
    A=\begin{pmatrix}1&0\\\beta X\inv&1\end{pmatrix}
	\begin{pmatrix}X &0\\0&Y-\beta X\inv
\alpha\end{pmatrix}
	\begin{pmatrix}1&X\inv \alpha\\0&1\end{pmatrix}.
\end{equation*}
\end{proof}

\subsection{Linear equations.}
Now we think of the
matrix $A=(a_{ij})_{i,j=1}^n$
as an even matrix over $\Lambda$. Then the quasideterminants
$|A|_{ij}$, thought of as elements of $\Lambda$, can only
exist if $i,j$ have the same parity (since otherwise
$A^{ij}$ is never invertible). Even if $i,j$ have the same
parity $|A|_{ij}$ may or may
not exist (as an element of $\Lambda$), but if it does then the
quasideterminant $|A_i(y)|_{ij} $ exists also. (Recall that
$A_i(y)$ is obtained by replacing the $i$th row of $A$ by
the row vector $y=(y_1,\dots,y_n)$.) Indeed, according to
\thetag{\ref{eq:defrecurs}}, if $y=(y_1,\dots,y_n)$, then,
using that $A_i(y)^{ij}=A^{ij}$,
\[
|A_i(y)|_{ij}=y_i-\sum_{\frac{p\ne i}{p q\ne j}}y_p
|A^{ij}|_{qp}\inv a_{qj}.
\]
Hence, using Lemma \ref{lem:quasidetqoutber}, we find,
whenever $|A|_{ij}$ exists,
\begin{equation}\label{eq:Aiyqdetasquot}
|A_i(y)|_{ij}=\begin{cases}
(-1)^{i+j}\ber\,(A_i(y))/\ber\,(A^{ij})&\text{ if $i\le
k$},\\
(-1)^{i+j}\ber^*\,(A_i(y))/\ber^*\,(A^{ij})&
\text{otherwise}.
		\end{cases}
\end{equation}
Now consider the linear system
\begin{equation}\label{eq:lineq2}
x A=y,
\end{equation}
with $x$, $y$ row vectors, and $A$ an even matrix over
$\Lambda$.
Using the properties of the quasideterminant
\thetag{\ref{eq:lineq2}} implies
\begin{equation*}
x_i|A|_{ij}=|A_i(y)|_{ij}.
\end{equation*}
Combining \thetag{\ref{eq:Aiyqdetasquot}} and Lemma
\ref{lem:quasidetqoutber}
proves   \thetag{\ref{eq:supercramer}}.

Of course, if one deals with infinite systems of equations
\thetag{\ref{eq:lineq}} over a super commutative ring one
has the same expressions \thetag{\ref{eq:supercramer}} for
the solutions, provided that one can define a
Berezinian (with the usual properties) of the infinite
matrices involved.

\section{Calculation of a super tau function}

The Baker functions for arbitrary line bundles over a
super elliptic curve were computed in \cite{Ra:SupElliptic}.
Note that a super elliptic curve, meaning an $N=1$ super Riemann
surface
of genus one with {\it trivial} spin structure $\mathcal N$ is not
a generic SKP curve and its Jacobian is not a supermanifold.
Still, the corresponding super tau function can be computed using
the Baker-tau relations (27,28) of \cite{DoSc:SuperGrass}:
\begin{equation}
\ber [\tilde{B}(W,F,S^*) \cdot S^*] = \left\{
\tau[W,F^T(1-Sz^{-1})^{-1}]/\tau(W,F^T) \right\} ^{-1},
\end{equation}
in the notation of that paper.
We choose the $2 \times 2$ matrix $S$ to be
$\text{diag}(\zeta,\zeta)$,
with $\zeta$ an even variable; then the entries of the
$2 \times 2$ matrix Baker function $\tilde{B}^{ij}(W,F,S^*)$
are simply the coefficients of the even Baker function
$w_{\bar 0}=B^{00}(\zeta) + B^{01}(\zeta)\theta$ and the odd Baker
function
$w_{\bar 1}=B^{10}(\zeta) + B^{11}(\zeta)\theta$ for the subspace
$FW$
of $\Sgr$, where $F$ is multiplication by some invertible formal
Laurent series $f(z)+\phi(z)\theta$, viewed as a $2 \times 2$
matrix
$\left[\begin{smallmatrix} f & 0 \\ \phi & f
\end{smallmatrix}\right]$.
This action of $F$ corresponds via the Krichever map to
deforming $\mathcal{L}$ by tensoring it with a bundle whose
transition function across the circle $|z|=1$ is
$f(z)+\phi(z)\theta$.

The results of \cite{Ra:SupElliptic} give the entries of the
Baker matrix as:
\begin{gather}
B^{00} = 1 + \frac{\alpha\delta}{2\pi i}\left[\frac{\Theta'(a)
\Theta'(\zeta-a)}{\Theta(a)\Theta(\zeta-a)} +
\frac{\Theta'(a)^2}{\Theta(a)^2} \right], \\
B^{01} = \alpha \frac{\Theta'(a)}{\Theta(a)} + \text{terms
proportional to } \delta, \\
B^{10} = \frac{\delta}{2\pi i} \left[
\frac{\Theta'(\zeta-a)}{\Theta(\zeta-a)} +
\frac{\Theta'(a)}{\Theta(a)}
\right], \\
B^{11} = 1 + \frac{\alpha\delta}{2\pi i} \left[
\frac{\Theta''(a)}{\Theta(a)} - \frac{\Theta'(a)^2}{\Theta(a)^2}
-\frac{\Theta''(\zeta-a)}{\Theta(\zeta-a)} +
\frac{\Theta'(\zeta-a)^2}{\Theta(\zeta-a)^2} \right].
\end{gather}

Here $\tau$ and $\delta$ are the even and odd moduli of the
super elliptic curve $X$; all theta functions appearing
are $\Theta_{11}(\bullet;\tau)$.
The parameters $(a,\alpha)$ label the point in the
super Jacobian corresponding to the deformed bundle
$\mathcal{L}$; if $F$ is written in the form
\begin{equation}
F = \exp \sum_{i=1}^{\infty} (t_{i}z^{-i} + t_{i-\frac12}\theta
z^{-i}),
\end{equation}
then they are linear combinations of the $t_i$ plus constants
labeling $\mathcal{L}$ itself.

Computing the superdeterminant, we find
\begin{equation} \label{tauanswer}
\frac{\tau(W,F^T(1-\zeta z^{-1})^{-1})}{\tau(W,F^T)} =
\frac{1 - \frac{\alpha\delta}{2\pi i}[\log \Theta(a-\zeta)]''}
{1 - \frac{\alpha\delta}{2\pi i}[\log \Theta(a)]''}.
\end{equation}
The function
\begin{equation}
(1-\zeta z^{-1})^{-1} = \exp \sum_{k=1}^{\infty} \frac{1}{k}
\frac{\zeta^k}{z^k}
\end{equation}
produces the usual shifts by $\zeta^k/k$ in the even parameters
$t_{k}$.
As a transition function for a line bundle,
$\frac{z}{z-\zeta}$, it corresponds to the
principal divisor $(0,0)-(\zeta,0)$,
so the corresponding deformation adds $-\zeta$ to the even
coordinate $a$ of the Jacobian, as we see on the right side of
\thetag{\ref{tauanswer}}.
We conclude that
\begin{equation}
\tau(W,F^T) = 1 - \frac{\alpha\delta}{2\pi i}[\log \Theta(a)]''.
\end{equation}
Note, however, that if the $2 \times 2$ matrix $F$ represents
multiplication by $f(z)+\phi(z)\theta$ then $F^T$ represents the
action of $f(z)+\phi(z)\partial_\theta$, which in terms of the
Krichever data is no longer just a deformation of $\mathcal{L}$,
but of $X$ as well.
That is, the Baker function for the orbit of $F$ enables us
to calculate the tau function for the ``dual" orbit of $F^T$,
a totally different type of deformation.

We see also that, consistent with the discussion in subsection
\ref{ss:AlgebrogeometrictauBaker}, the super tau function is a
genuine
function,
not a section of a bundle, as the multivaluedness of the
theta function is eliminated by the logarithmic derivatives.
This is again due to the fact that the restricted group
of $F^T$'s considered here acts without central extension in
the Berezinian bundle.

\providecommand{\bysame}{\leavevmode\hbox
to3em{\hrulefill}\thinspace}

\end{document}